\UseRawInputEncoding
\documentclass[12equipped]{article}

\usepackage{amssymb,amsmath}

\input epsf.tex

\title{Analyticity and resurgence in wall-crossing formulas}

\author {Maxim Kontsevich, Yan Soibelman}

\begin{document}
\maketitle

\newcommand{\VV}{{\mathcal V}}
\newcommand{\CC}{{\mathcal C}}
\newcommand{\LL}{{\mathcal L}}
\newcommand{\MM}{{\mathcal M}}
\newcommand{\NN}{{\mathcal N}}
\newcommand{\OO}{{\mathcal O}}
 \newcommand{\ZZ}{{\mathcal Z}}
\renewcommand{\O}{{\mathcal O}}
\newcommand{\XX}{{\mathcal X}}
\newcommand{\Ome}{{\Omega}^{3,0}}
\newcommand{\E}{{\mathcal E}}
\newcommand{\F}{{\mathcal F}}
\newcommand{\TT}{{\mathcal T}}
\newcommand{\g}{{\mathfrak g}}
\newcommand{\GGamma}{{\underline{\Gamma}}}

\renewcommand{\k}{{\bf k}}
\newcommand{\kk}{\overline{\bf k}}

\newcommand{\op}[1]{\operatorname{#1}}

\newtheorem{thm}{Theorem}[subsection]
\newtheorem{defn}[thm]{Definition}
\newtheorem{lmm}[thm]{Lemma}
\newtheorem{rmk}[thm]{Remark}
\newtheorem{prp}[thm]{Proposition}
\newtheorem{conj}[thm]{Conjecture}
\newtheorem{exa}[thm]{Example}
\newtheorem{cor}[thm]{Corollary}
\newtheorem{que}[thm]{Question}

\newtheorem{ack}{Acknowledgements}

\newcommand{\B}{{\bf B}}
\newcommand{\C}{{\bf C}}
\newcommand{\K}{{\bf k}}
\newcommand{\R}{{\bf R}}
\newcommand{\N}{{\bf N}}
\newcommand{\Z}{{\bf Z}}
\newcommand{\Q}{{\bf Q}}
\newcommand{\G}{\Lambda}
\newcommand{\A}{A_{\infty}}
\newcommand{\M}{{\mathsf{M}}}
\newcommand{\epi}{\twoheadrightarrow}
\newcommand{\mono}{\hookrightarrow}
\newcommand\ra{\rightarrow}
\newcommand\uhom{{\underline{Hom}}}
\renewcommand\O{{\cal O}}
\newcommand{\epp}{\varepsilon}


\abstract{We introduce the notion of analytic stability data  on the Lie algebra of vector fields on a torus. We prove that 
the subspace of analytic stability data is open and closed in the topological space of all stability data. 
We formulate a general conjecture which explains how analytic stability data give rise to resurgent series. This conjecture is checked in several examples.}

\tableofcontents

\section{Introduction}\label{Intro}

\subsection{Analytic stability data and wall-crossing structures}

Stability data on a Lie algebra  $\g=\oplus_{\gamma\in \Gamma}\g_\gamma$ graded by the free abelian group $\Gamma\simeq \Z^n$ are defined by a homomorphism of abelian groups $Z: \Gamma\to \C$ (central charge) and a collection of elements $a(\gamma)\in \g_\gamma, \gamma\in \Gamma-\{0\}$ satisfying one axiom called the Support Property (see [KoSo1], Section 2). The set $Stab(\g)$ of stability data on $\g$  carries the natural Hausdorff topology which makes it locally homeomorphic to the space of central charges $Hom(\Gamma,\C)$. This topology encodes the wall-crossing phenomenon, i.e. the change of the collection $(a(\gamma))_{\gamma\in \Gamma-\{0\}}$ along a path when it crosses certain real codimension one subvarieties in $Stab(\g)$ (walls).

Stability data on graded Lie algebras were introduced in the loc.cit. for the purposes of Donaldson-Thomas theory. They found further applications in complex integrable systems and Mirror Symmetry (see [KoSo3]).

In this paper we consider a special case of the general theory, namely the graded Lie  algebra 
$\g=Vect_\Gamma$ of  algebraic vector fields on the complex algebraic torus ${\bf T}_\Gamma:=Hom(\Gamma, {\bf G}_m)$.
Our main object of study is a subset $Stab^{an}(Vect_\Gamma)\subset Stab(Vect_\Gamma)$ of {\it analytic stability data}. The definition of the latter  is quite involved. It is based on the new approach to the notion of stability data proposed in Section \ref{stab data}.  Hypothetically, analyticity of stability data is equivalent to exponential bounds on the growth of $a(\gamma)$ as $|\gamma|\to \infty$.

Main result concerning analytic stability data   is proved in Section \ref{open and closed}. It says that the subset $Stab^{an}(Vect_\Gamma)$ is open and closed in the topological space $Stab(Vect_\Gamma)$.

Recall that the notion of stability data on a graded Lie algebra was generalized in [KoSo3] to the notion of {\it wall-crossing structure}. The latter is roughly speaking given by a family of stability data on a local system of graded Lie algebras on a topological space. We generalize the notion of analytic stability data  to this framework thus arriving to the notion of {\it analytic wall-crossing structure}. Donaldson-Thomas theory, Gromov-Witten theory, complex integrable systems, $4d$ $\mathcal {N}=2$ gauge theories, exponential integrals (possibly infinite-dimensional, e.g. Feynman integrals) are just few examples in which analytic wall-crossing structures appear or expected to appear.

Another motivation for this paper comes from the classical subject of  resurgence (see [Ec]). 
In this paper we propose a general conjecture which says how to construct resurgent series in a simple explicit way starting with analytic stability data. We show that our conjecture agrees with known results which include   \'Ecalle-Voronin theory, Voros resurgence of WKB asymptotics as well as exponential integrals.

In general the notion of wall-crossing structure (WCS for short) gives a useful tool for study  properties of generating series arising in wall-crossing formulas (see [KoSo1,3]). It looks plausible that one can  prove algebraicity, analyticity or resurgence of those generating series in a uniform way, based on the  properties of the corresponding wall-crossing structures (see e.g. [BarbSt], [CorSh], [GaMoNe1], [GuMaPu], [GarGuMa] for some examples of situations where such theory can be useful). This paper is devoted to analyticity and resurgence. We will discuss algebraicity in a separate paper.

\subsection{Contents of the paper}

In Sections \ref{cones}-\ref{alt topology} we review main definitions and facts about stability data and propose a reformulation of this notion in terms of quotient sets associated with {\it cyclic covers} of $\R^2-\{0\}$ by sectors. We also discuss an upper bound on the support of stability data in continuous families in Section \ref{families and support property}.
In the remaining part of Section \ref{stab data} we restrict ourselves to the case $Vect_\Gamma$. E.g. we recall  an interpretation of stability data on $Vect_\Gamma$ in terms of certain formal schemes.

 In Section \ref{analytic stab data} we define analytic stability data on $Vect_\Gamma$ in two different ways:
 
 a) via the language of quotient sets;
 
 b) via analytifications of  formal schemes.

The key result is Theorem \ref{analyticity of factors} about analyticity of factors. 
Using it we prove in Section \ref{open and closed} that the subset  of analytic stability data is open and closed.  Another topic discussed in Section \ref{analytic stab data} is a generalization of analytic stability data to analytic wall-crossing structures.

In Section \ref{analytic stab data and resurgence} we introduce an extra complex parameter to the story. It plays a role of Planck constant.  An important situation in which the new parameter appears naturally is the relation between analytic stability data and resurgence. We formulate the general Conjecture \ref {analyticity and resurgence} which underlies resurgence  of many series in mathematics and physics. The Conjecture is illustrated  in several examples, including   \'Ecalle-Voronin theory about normal forms of  germs of analytic automorphisms of the line tangent to the identity. 
 In Section \ref{variation of central charge} we consider  the central charge  as an additional parameter. Arising mathematical structure resembles the one in non-commutative Hodge theory (see [KaKoPa]). After motivational results about the formal classification of certain non-linear connections, we present the new structure in axiomatic way. 
 
 In the last Section \ref{final remarks} we briefly discuss other classes of stability data like e.g. algebraic stability data. The latter are responsible for {\it algebraicity} of many series appearing in Donaldson-Thomas theory. Also we expand the list of examples illustrating the  Conjecture \ref {analyticity and resurgence}. Those include Gromov-Witten theory, Voros resurgence of WKB solutions of differential equations with small parameter as well as exponential integrals. In the Appendix we give a proof of an important result about independence of the set of stability data of the cyclic cover.

{\it Acknowledgements.} The work  of Y.S. is partially supported by the Munson-Simu Star Excellence Award of KSU. He also thanks IHES for excellent research conditions.  This paper is partly a result of the ERC-SyG project, Recursive and Exact New Quantum Theory (ReNewQuantum) which received funding from the European Research Council (ERC) under the European Union's Horizon 2020 research and innovation program under grant agreement No 810573.

\section{Stability data on graded Lie algebra revisited, and the case of Lie algebra of  vector fields}\label{stab data}

\subsection{Graded Lie algebras and cones}\label{cones}

Let $\Gamma$ be a free abelian group of rank $n$, i.e. $\Gamma\simeq \Z^n$ (but the isomorphism is not part of the structure). Suppose we are given a $\Gamma$-graded Lie algebra $\g$ over the field $\Q$, i.e. 
$\g=\oplus_{\gamma\in \Gamma}\g_\gamma, [\g_\gamma, \g_\mu]\subset \g_{\gamma+\mu}$. Let $\Gamma_\R=\Gamma\otimes \R$. Then $\Gamma\subset \Gamma_\R$ is a lattice. The dual vector space is $\Gamma_\R^\ast=\Gamma^\vee\otimes \R$, where $\Gamma^\vee=Hom(\Gamma, \Z)$ is the dual lattice.

Recall that a {\it strict cone} in $\Gamma_\R$ is a cone in $C\subset \Gamma_\R$, which does not contain a straight line.
If $rk\,\Gamma=2$ instead of strict cones we will speak about {\it strict sectors}.

For any strict convex cone $C\subset \Gamma_\R$ we define the pronilpotent Lie algebra $\g_C:=\prod_{\gamma\in C\cap \Gamma-\{0\}}\g_\gamma$ as well as the corresponding pronilpotent Lie group $G_C:=exp(\g_C)$. Clearly if $C^\prime\subset C$ is a convex subcone then $G_{C^\prime}\subset G_C$. In the $2$-dimensional case we will use the following terminology.

\begin{defn} An  admissible sector  is a sector of one of the following types: 

a) a strict sector $V$ which is of the form $V=\{z=re^{i\theta}|r>0, \theta \in S\}$, where $S\subset \R$ is a non-empty connected subset of  length smaller than $\pi$;

b) open half-plane, i.e. we allow in  a) the set $S$ to be an open interval of  length $\pi$.
\end{defn}

Sometimes we will call admissible sectors of the type a) {\it strictly admissible}.

Accordingly,  we say that $V$ is  {\it open, closed or semiclosed}  if $S$ is open, closed or semiclosed. In particular,  $V$ can be a ray with the vertex removed (in this case $S$ is a point and the removed vertex is $0\in \R^2$). In what follows we orient $\R^2$ in the clockwise order. Then each sector has its {\it left} and {\it right boundary} rays, so that left precedes right. Similarly, if  we split a sector by an admissible ray into two, then they are called {\it left} and {\it right} depending on which of their left boundary rays precedes the other one.

Let $Z:\Gamma\to \C$ be a homomorphism of abelian groups which we will call a {\it central charge}. We denote by $Z_\R:\Gamma_\R\to \C=\R^2$ its $\R$-linear extension or simply by $Z$, if it does not lead to a confusion. Let us assume that $rk\, \Z_\R=2$. 
Let $C$ be a closed strict convex cone  such that the restriction of $Z_\R$ to $C$ is proper. Then $V=Z_\R(C)$ is a closed strict sector. For any  ray $l\subset V$ we have the corresponding cones $C(l)=Z_\R^{-1}(l)\cap C$, $C_\pm=Z_\R^{-1}(V_\pm)$, where $V_\pm$ are strict   sectors which are connected components of the complement $V-l$ ordered  in such a way that $V_+$ is on the left of  $V_-$. Then $G_C$ contains $G_{C(l)}, G_{C_\pm}$ and moreover $G_C=G_{C_+}G_{C(l)}G_{C_-}$ understood as the product of subsets of a group.

Let now $V$ be an {\it open} admissible sector. Let us fix the central charge $Z$. Consider the set $Cones(V, Z)$ consisting of strict convex cones $C\subset \Gamma_\R$ which satisfy the properties: 

a) $Z(C-\{0\})\subset V$;

b) $\overline{C}\cap Ker(Z_\R)=\{0\}$. 

The set $Cones(V, Z)$ carries a partial order (inclusion
of cones). 

Consider also sets  $Cones_{\pm}(V, Z)$ consisting of such {\it closed} strict convex cones $C$ that they satisfy the condition b) above
as well as the following conditions

$c)_\pm$ The closed set ${Z_\R(C)}$ belongs to the union of $V$ and its left (for $+$) or right (for $-$) boundary ray.

We denote by $G_V=G_{V,Z}=G_V^{new}$ the inductive limit $\varinjlim_{C\in Cones(V,Z)}G_C$. Notice that $G_V$ is {\it always} a group, even if $V$ is not strict (e.g. it can be an open half-plane). Since $V$ is open, in the definition of $G_V$ we can take the inductive limit over the smaller set of {\it closed} strict convex cones $V$.  

\begin{rmk}
The reader should not confuse the group $G_V=G_V^{new}$ introduced above with  the group $G_V$ introduced in [KoSo1]. In order to distinguish them we will denote the group defined in the loc.cit. by $G_V^{old}$. In particular, the elements $A_V$ introduced in [KoSo 1]  belong to $G_V^{old}$, but not necessarily to $G_V$. Notice that one can  define as above
the group $G_V$ for a closed admissible $V$. In this case $G_V$ and $G_V^{old}$ coincide.
In particular we have the group $G_l=G_l^{old}$ for an admissible ray $l$. Furthermore $G_V=\varinjlim_{W\subset V}G_W^{old}=
\varinjlim_{W\subset V}G_W$, where the inductive limit is taken over all closed strict admissible subsectors of $V$.
\end{rmk}

Let us also define groups $G_V^{\pm}$ as above by taking inductive limits over $Cones_{\pm}(V, Z)$. Then $G_V=G_V^+\cap G_V^-$.

Recall (see [KoSo1], Section 2) that if the ray $l$ divides a strict open admissible sector $V$ into the disjoint union $V=V_0\cup l\cup V_1$ where $V_0, V_1$ are clockwise ordered open admissible subsectors of $V$ then $G_V^{old}=G_{V_0}^{old}G_lG_{V_1}^{old}$ (product of subsets in the group $G_V$). This observation implies the following result which will be useful later in the paper.

\begin{prp}\label{group factorization} Let $V$ be an open admissible sector, $l_1,..., l_k\subset V$ are admissible rays, clockwise ordered, and  $V=V_0\sqcup l_1\sqcup V_1\sqcup l_2\sqcup V_2\sqcup...\sqcup l_k\sqcup V_k$ is the clockwise ordered disjoint decomposition. Then we have $G_V=G_{V_0}^+G_{l_1}G_{V_1}^{old}G_{l_2}G_{V_2}^{old}....G_{l_k}G_{V_k}^-$ (product of subsets in the group $G_V$).

\end{prp}

\subsection{Reminder on stability data}\label{stability reminder}

We are going to recall following  [KoSo1], Section 2 the notion of stability data on a graded Lie algebra.

\begin{defn}
A stability data $\sigma$ on $\g$ consists of the following:

1)  A central charge $Z:\Gamma\to \C$,

2) Elements $a(\gamma)\in \g_\gamma, \gamma\in \Gamma-\{0\}$,

Let now $a=(a(\gamma))_{\gamma\in \Gamma-\{0\}}$ and $Supp(\sigma):=Supp(a)$ denotes the set of $\gamma$ such that $a(\gamma)\ne 0$.
Then the only axiom for the data 1) and 2) says that there exists a  non-trivial  quadratic form $Q$ on $\Gamma_\R$ such that $Q(\gamma)\ge 0$ for
$\gamma\in Supp(a)$ and $Q_{|Ker\, Z}<0$.
\end{defn}

Let  now $V$ be a strict admissible sector and $Z,Q$ mean  the same as in the definition.
 Consider a  strict convex cone
$C(V):=C(V,Q,Z)$ which is the convex hull of the set
of such $u\in \Gamma_\R$ that $Q(u)>0$ and $Z(u)\in V$.  Then we have the group $G_{C(V)}=G_{C(V,Q,Z)}$ associated with this cone.

An equivalent way to describe stability data $\sigma$  on $\g$ is via the 
collection of elements $A_V:=A_V^\sigma\in G_{C(V)}$, where $V$ runs through the set of all strict admissible sectors,
satisfying the following property (called Factorization Property in the loc. cit.):

If $V=V_1\sqcup V_2$ is a decomposition of  an admissible sector into the disjoint union of two admissible sectors such that $V_1$ precedes $V_2$
in the clockwise order on $\R^2$ then $A_V=A_{V_1}A_{V_2}$ (this makes sense since $G_{C(V_i)}\subset G_{C(V)}, i=1,2)$.

\begin{rmk}\label{quadratic form}

a) Factorization Property is equivalent to the fact that $A_V$ is equal to the clockwise product $\prod_{l\subset V}A_l$ over all rays $l$ in $V$ (notice that rays are admissible sectors). All but countably many factors $A_l$ are equal to $1$, hence the product makes sense.

b) Two descriptions given above are related in the following way: the elements $a(\gamma)\in \g_\gamma, \gamma\in \Gamma-\{0\}$ are $\gamma$-components of the elements $log(A_l)$
for $l=\R_{>0}\cdot Z(\gamma)$. 

c)  The fact that the quadratic form $Q$ does exist is called the Support Property in [KoSo1].

\end{rmk}

The group $G_{V}^{old}$  mentioned in the previous subsection coincides with the inductive limit $G_{C(V,Z)}$ of the groups  $G_{C(V, Z, Q)}$  in the case when $V$ is  admissible.  More precisely,  we define the pronilpotent group $G_{V,Q, Z}$ for fixed $Q$ and $Z$ as in this subsection,
and then take the inductive limit of $G_{V,Q, Z}$ over the set of all quadratic forms $Q$ (notice that the set of quadratic forms  is partially ordered with respect to the natural order).  Sometimes will skip $Q$ and $Z$ from the notation.

Let $\sigma\in Stab(\g)$.
Having in mind the  terminology of  [KoSo1] the element $a(\gamma)\in \g_\gamma$ will be sometimes called {\it numerical rational Donaldson-Thomas invariant of $\gamma$}.  The word {\it rational} reflects the difference with {\it numerical Donaldson-Thomas  invariants $\Omega(\gamma)$} consider in loc. cit. which are inverse M\"obius  transforms of $a(\gamma)$.

\begin{defn}
1) We say that the central charge is rational if $Z(\Gamma)\subset \Z^2$. Stability data are called rational if the central charge is rational.

2) We call $dim\,Z_\R(\Gamma_\R)$ the rank of the central charge (or the corresponding stability data).

\end{defn}

Thus we can have three possibilities:

0) $Z=0$. It follows that all $A_V=1$. This is the rank zero case.

1) $dim_\R\,Z_\R(\Gamma_\R)=1$. This means that $Z_\R(\Gamma_\R)=l_-\sqcup\{0\}\sqcup l_+$ for some rays $l_\pm$.
Hence stability data are determined by $A_{l_\pm}$. This is the  rank one case.

2) $dim_\R\,Z_\R(\Gamma_\R)=2$. Therefore $Z_\R(\Gamma_\R)=\C$. This is the   rank two case, which corresponds to generic stability data.

\subsection{Stability data as a quotient set}\label{alternative def}

Here we propose a description of the set $Stab(\g)$ in terms of what we call cyclic covers of $\R^2-\{0\}$. Later we will explain the topology in these terms. We assume the set up of Subsection \ref{cones}. 

\begin{defn} Let $J$ be a finite cyclically ordered set. A closed (resp. open, semiclosed) cyclic cover of $\R^2-\{0\}$ parametrized by $J$  is given by a finite clockwise ordered set $\mathcal{V}=(V_i)_{i\in J}$ of 
closed (resp. open , semiclosed) set of sectors $V_j=\{z=re^{i\theta}|r>0, \theta_j^-\le \theta\le \theta_j^+\}$ (resp. obvious modifications of the inequalities for the angles in the open and semiclosed cases) such that $ \theta_j^-\le \theta_{j+1}^-\le \theta_{j}^+, \theta_{j+1}^+\ge \theta_j^+$. We require that the sectors of the cyclic cover are admissible, i.e. they are strict or they are open half-planes. We call the cyclic cover rational if the sectors are bounded by rational rays.
\end{defn}

In other words, a cyclic cover  is a cover by a finite cyclically ordered set of  closed (resp. open, semiclosed) sectors such that two consecutive sectors have non-empty intersection. The latter can be e.g. a common boundary  ray without vertex. An example of cyclic cover is given by the wheel of admissible sectors (see  [KoSo3]), which we are going to use later in the paper.

Let $Z:\Gamma\to \C$ be a central charge. Then for any open sector $V$ which is either strictly admissible or an open half-plane we have defined the group
$G_V=G_{V, Z}$.  Recall that it is given by $\varinjlim G_C$, where the inductive limit is taken over all strict convex  cones $C\subset \Gamma_\R$ such that $Z(C-\{0\})\subset V$. 

For any open cyclic cover  $(V_i)_{i\in J}$ by strictly admissible sectors or open half-planes let us consider the group $\prod_{i\in J}G_{V_i}$. It is acted by the group $\prod_{i\in J}G_{V_j\cap V_{j+1}}$ in the natural way, namely
$G_{V_j\cap V_{j+1}}$ acts on $G_{V_j}\times G_{V_{j+1}}$ via the natural embeddings $G_{V_j\cap V_{j+1}}\to G_{V_j}, G_{V_j\cap V_{j+1}}\to G_{V_{j+1}}$.
Since this action is free we have a well-defined quotient set  $\MM_{\mathcal V}=\MM_{(V_i)_{i\in J}}$. We warn the reader in advance that the topology on $\MM_{\mathcal V}$ which we will discuss later  is {\it not} the naturally  induced topology on the quotient set.

We would like to prove that the quotient space does not depend on the admissible cover. We start with some preliminaries. First, we remark that cyclic covers form a category $Cov$. Objects of $Cov$ are open cyclic covers. Morphism  $(V_i)_{i\in J}\to (V_{i^\prime})_{i^\prime\in J^\prime}$ is given by an {\it embedding} $\psi: J\to J^\prime$ which satisfy the following conditions:

a1) $V_i\subset V_{\psi(i)}^\prime$;

a2) $V_{\psi(i+1)}^\prime$ is position on the right from $V_{\psi(i)}^\prime$ in the clockwise order (i.e. left boundary ray of $V_{\psi(i+1)}^\prime$ is obtained from the left boundary ray of $V_{\psi(i)}^\prime$ by a clockwise rotation by a non-zero angle.

Equivalently we can reformulate a2) in the following way:

a2)' The map $\psi$ respects the natural relation on the open sectors $(V_i)_{i\in J}$ ``to be positioned on the right" (i.e. $\psi$ preserves the cyclic order). 

Sometimes instead of ``positioned on the right" we will simply say ``follows" and write $i_1<i_2$ if $i_2\ne i_1$ follows $i_1$ in the cyclic order. 

Then it is clear that we can compose to morphisms and get a morphism. 
\begin{prp}\label{functor}
We have a well-defined functor $\Phi: Cov\to Sets$ given by $\Phi(\mathcal V)=\MM_{\mathcal V}$ which  is induced by the natural map $\prod_{i\in J}G_i\to \prod_{i^\prime\in J^{\prime}}G_{i^\prime}$ coming from  the embedding $\psi: J\to J^\prime$.
\end{prp}
{\it Proof.} Notice that $i<i+1$ implies $\psi(i)<\psi(i+1)$. Then $V^\prime_{\psi(i+1)}$ follows $V^\prime_{\psi(i)}$. Consider the interval $\psi(i)=i_1^\prime<i_1^\prime+1:=i_2^\prime<...<i_{k-1}^\prime+1:=i_k^\prime=\psi(i+1)$. Notice that $V_{\psi(i)}\cap V_{\psi(i+1)}\subset V_{i_p^\prime}, 1\le p\le k$.  Then there is a well-defined map 

$$G_{V_{\psi(i)}}^{old}\times_{G^{old}_{V_{\psi(i)}\cap V_{\psi(i+1)}}}G_{\psi(i+1)}^{old}\to G_{V_{i_1^\prime}}^{old}\times_{G_{V_{i_1^\prime}\cap V_{i_2^\prime}}^{old}}G_{V_{i_2^\prime}^{old}}\times_{G_{V_{i_2^\prime}\cap V_{i_3^\prime}}^{old}}G_{V_{i_3^\prime}}^{old}\times...\times G_{V_{i_k^\prime}}^{old}$$
The result follows. $\blacksquare$.

\begin{prp} \label{connected nerve} The category $Cov$ is connected (i.e. its nerve is connected).

\end{prp}
{\it Proof.} The category $Cov$ contains a subcategory $Cov_{1/2}$ consisting of open covers by open half-planes. This category is filtered and hence connected. Furthermore there is a map on objects (not a functor) $\phi: Ob(Cov)\to Ob(Cov_{1/2})$ defined such as follows. To each open admissible sector $V$ we associate an open half-plane $\alpha_V$ bounded by the straight line, which contains {\it left} boundary ray of $V$ (and hence $V\subset \alpha_V$). Hence to an open cyclic cover $(V_i)_{i\in J}$ we associate an open cyclic cover by half-planes $(\alpha_i=\alpha_{V_i})_{i\in J}$. This gives our map $\phi$. Notice that we have a natural morphism in $Cov$ which assigns to the open cyclic over $(V_i)_{i\in J}$ the cover $(\phi(V_i))_{i\in J}$.
 The result follows.
$\blacksquare$

\begin{thm}\label{independence of cyclic cover}
The quotient set $\MM_\VV$ does not depend on the open cyclic cover. More precisely:
there exists a set $\MM$ and a collection of bijections $iso_\VV: \MM_\VV \to \MM$ one for each open cyclic cover $\VV$ such that if $\psi:\VV_1\to \VV_2$ is a morphism in $Cov$ then $iso_{\VV_1}=iso_{\VV_2}\circ \Phi(\psi)$.
\end{thm}

Proof of the Theorem will be given in Appendix. For that we will develop a purely combinatorial language for  open cyclic covers and corresponding groups.
\begin{rmk}  Let us define the sets $\MM^{old}_{\VV}$ for closed and semiclosed cyclic covers similarly to $\MM_\VV$, but using the groups $G_{V_i}^{old}$ instead of $G_{V_i}$. Then the Theorem \ref{independence of cyclic cover} still holds.

\end{rmk}

We will use this Remark in Section 6 without further comments.

\subsection{Reminder on the topology on the space of stability data}\label{stability topology}

There is  a Hausdorff topology on the space $Stab(\g)$ of stability data on $\g$ such that the natural projection of $Stab(\g)$ to the space of central charges $Z$ is a local homeomorphism. This topology was defined in Section 2 of [KoSo1]. We recall it here (see [KoSo1], Section 2.3 for details).  In order to define the topology it suffices to say what means to have a continuous family of stability data parametrized by a topological space $X$. According to loc. cit., the family of stability data $\sigma_x, x\in X$ is said to be continuous at $x_0\in X$ if the following conditions are satisfied:

a) the corresponding family of central charges $Z_x, x\in X$ is continuous at $x_0$;

b) there is an open neighborhood $U_{x_0}$ of $x_0$ and a quadratic form $Q_{x_0}$ which satisfies the condition c) from the Remark \ref{quadratic form} for all $x\in U_{x_0}$;

c) most important condition consists of the requirement that for any closed admissible sector $V\subset \R^2$ such that its boundary does not contain points $Z_{x_0}(Supp(\sigma_{x_0}))$ the map $x\mapsto log(A_{V}^{\sigma_x})\in \prod_{\gamma \in \Gamma}\g_\gamma$ is continuous at $x=x_0$.

Here $A_{V}^{\sigma_x}$ denotes the element $A_V$ constructed for the stability data $\sigma_x$, and the product topology is discrete.

Equivalently, the  condition c) means that if we cut out a triangle $\Delta$ from $V$ by a generic straight line, then for all $\gamma\in \Gamma\cap \Delta$ the element $a_x(\gamma)\in \oplus_{\gamma\in \Gamma\cap \Delta}\g_\gamma$ of  $log(A_{V}^{\sigma_x})$ is {\it locally constant} in a neighborhood of $x_0$. It is not true in general that there is a neighborhood of $x_0$ which ``serves" in this way all $\gamma\in \Gamma\cap V$.

The above-described Hausdorff topology underlies the {\it wall-crossing formulas} from the loc.cit. {\it Walls of first kind} according to [KoSo1] are subsets in $Stab(\g)$ consisting of those stability data for which the central charge $Z$ maps a sublattice of $\Gamma$ of the rank at least $2$ into a line. Generically walls have real codimension $1$, so they are unavoidable in $1$-dimensional continuous families, e.g. in paths in $Stab(\g)$.

\subsection{Topology from the point of view of  cyclic covers}\label{alt topology}

Let us discuss now the topology in terms of the approach of  Subsection \ref{alternative def}. This allows us to propose a much more explicit definition of the topology on $Stab(\g)$. 
 
\begin{defn} \label{open basis} Let $\sigma\in Stab(\g)=\MM_\VV$ where $\VV=(V_i)_{i\in J}$ is an open cyclic cover. Then for any representative $(g^\sigma_{V_i})_{i\in J}\in \prod_{i\in J}G_{V_i}$ the basis of open neighborhoods of $\sigma$ consists of stability data $\sigma^\prime$ with the same representatives, i.e. $(g^{\sigma^\prime}_{V_i})_{i\in J}=(g^\sigma_{V_i})_{i\in J}$ and central charges $Z^\prime$ which belong to a basis of open neighborhoods of the central charge $Z$ of $\sigma$.

\end{defn}

One can prove the following result.

\begin{prp} The Definition \ref{open basis} is equivalent to the definition from [KoSo1].

\end{prp}

In particular, in  the case of a cyclic cover by open half-planes, each group $G_{V_i}$ is  an inductive limit of the groups  $G_{W_i^{(m)}}$ for closed subsectors $W_i^{(m)}$ such that $V_i=\varinjlim_m W_i^{(m)}$ (see Subsection \ref{alternative  def}). For a fixed $m\ge 1$ we can define an open neighborhood $U_m$ of $\sigma$. It consists of  those stability data for which the central charge belongs to such an open neighborhood of $Z\in Hom(\Gamma, \C)$ that the representatives are the same as those for $\sigma$, e.g.
$(A^\sigma_{W_i^{(m)}})_{i\in J}$. To bigger values of $m$ correspond smaller open subsets $U_m$.

\begin{prp}\label{limit}
The stability data $\sigma$ with central charge $Z$ is a limit as $n\to \infty$ of stability data $\sigma_n$ with central charges $Z_n$ iff for any $\epsilon>0$  there exists $n_0$ such that for all $n\ge n_0$ all $\sigma_n$ have the same representative $(g^\sigma_{V_i})_{i\in J}$ and $|Z-Z_n|<\epsilon$ (here we can use any norm on $\Gamma_\R$).

\end{prp}
{\it Proof.} Follows from the Definition \ref{open basis}. $\blacksquare$

New definition of the topology allows us to give a simple example of the continuous path in the space $Stab(\g)$. We will use it in Section 3. Namely, let $\sigma$ be stability data with central charge $Z$ and $(g_{V_i})_{i\in J}$ be its representative associated with an open cyclic cover  $V_i, i\in J$. Let us choose closed strict convex cones $C_i\in \Gamma_\R, i\in J$ such that
$Z(C_i-\{0\})\subset V_i$. Consider now a continuous path $t\mapsto Z_t, t\in [0,1]$ in $Hom(\Gamma,\C)$ such that $Z_0=Z$. Let us define the stability data $\sigma_t$ by the requirement that it have representative $(g_{V_i}^t)_{i\in J}=g_{V_i}^{t=0}:=g_{V_i}\in G_{C_i}$ for all $t\in [0,1]$ and all $ i\in J$. Then $\sigma_t$ is a continuous path containing $\sigma$ as its endpoint.

\subsection{Families of stability data and Support Property}\label{families and support property}

Let $\Gamma$ be a free abelian group of rank $n$. Let us fix $Z_0\in Hom(\Gamma, \C)$  as well as an open neighborhood $U_0\subset Hom(\Gamma,\C)$ of $Z_0$. Main purpose of this subsection is to prove the following result.

\begin{prp}\label{support for families}
There exists a quadratic form $Q=Q(Z_0, U_0)$ on $\Gamma_\R$ satisfying the following properties:

1) $Q_{|Ker\, Z_0}<0.$

2) For any graded Lie algebra $\g=\oplus_{\gamma\in \Gamma}\g_\gamma$ and a continuous family of stability data $\sigma_u\in Stab(\g), u\in U_0$ the support of $\sigma_{Z_0}$ is contained in the set $\{\gamma\in \Gamma-\{0\}| Q(\gamma)>0\}$.

\end{prp}

\begin{rmk} Notice that the ``upper bound" on $Supp(\sigma_{Z_0})$ which follows from the above Proposition does not depend  on the graded Lie algebra. Furthermore, for given $\g$ by definition of the stability data $\sigma_{Z_0}$ we have the quadratic form $Q_0$ which satisfies  the property 1) and which is positive on $Supp(\sigma_{Z_0})$. We remark that the form $Q$ from the Proposition does not have to coincide with $Q_0$. The latter depends on the stability data $\sigma_{Z_0}$, while the former depends only on the central charge $Z_0$ and its neighborhood. 
\end{rmk}

Proof of the Proposition \ref{support for families} will consist of several steps and occupy the rest of this subsection.
Main idea is given in the following lemma.

\begin{lmm} \label{central charges on the line} Suppose we are given a continuous family $(\sigma_b)_{b\in B}$ of stability data on $\g=\oplus_{\gamma\in \Gamma}\g_\gamma$ parametrized by a connected topological space $B$. Assume that for all $b\in B$ the rank of $Z_{\R,b}$ (i.e. of the $\R$-linear extension of $Z_b$) does not depend on $b$ (i.e. $rk\,Z_{\R,b}\in \{0,1,2\}$). Let $\kappa\subset \R^2$ be a straight line passing through the origin, and $\alpha_b=Z_b^{-1}(\kappa)\subset \Gamma_\R$ the corresponding subspace.  If $\alpha=\alpha_b$ does not depend on $b\in B$ then the elements $a_b(\gamma)$ do not depend on $b$ as long as $\gamma\in \alpha\cap (\Gamma-\{0\})$.
\end{lmm}

{\it Proof of the Lemma.} We will prove the Lemma in the rank $2$ case. Rank $1$ case is similar, while in the rank $0$ case all $a_b(\gamma)=0, b\in B$, and the Lemma is obvious.
Since $B$ is connected, it suffices to show that $a_b(\gamma)$ is locally constant for each $\gamma\in \alpha\cap (\Gamma-\{0\})$, where $\alpha=\alpha_b$ is independent on $b$. Let $b_0\in B$ and $Q_{b_0}$ is a quadratic form on $\Gamma_\R$ such that $Q_{b_0}$ is negative on $Ker\, Z_{b_0}$ and positive on $Supp(\sigma_{b_0})$. Then $Q_{b_0}$ enjoys the same property for all $\sigma_b$ with $b$ belonging to a small open neighborhood of $b_0$.

Let us denote by $\kappa_\pm$ two open opposite rays forming the line $\kappa$ such that they have a common vertex $0\in \R^2$. Consider an open cyclic cover of $\R^2-\{0\}$ by open admissible sectors $V_1, V_2, V_3, V_4$ such that $\kappa_+\subset V_1, \kappa_-\subset V_3$. Recall the strict convex open cones $C(V_i), i=1,2,3,4$ which are convex halls of $Z_{b_0}^{-1}(V_i)\cap \{x\in \Gamma_\R|Q_{b_0}(x)>0\}$. Recall also the elements $A_{V_i}\in G_{V_i}, i=1,2,3,4$ which give a representative for $\sigma_b$ for all $b$ sufficiently close to $b_0$ (see Definition \ref{open basis}). Recall the factorization  formula for the sector $V_1$, i.e. $A_{V_1}=A_{V_1}^{>,b}A_{\kappa_{+,b}}A_{V_1}^{<,b}$. The LHS does not depend on $b$, but a priori this is not true for the factors. The assumption that $\alpha=\alpha_b$ does not depend on $b$ implies  that the factors also do not depend on $b$. In particular $A_{\kappa_+}:=A_{\kappa_{+,b}}$ does not depend on $b$, and hence all $a_b(\gamma)$ with $\gamma\in C(V_1)$ do not depend on $b$. Similarly for $A_{\kappa_-}:=A_{\kappa_{-,b}}$ and $C(V_3)$.  This finishes the proof of the Lemma.

\begin{cor} Under the assumptions of the Lemma \ref{central charges on the line} let us fix $b_0\in B$ and consider the set $S_{b_0}$ consisting of $\gamma\in \alpha\cap (\Gamma-\{0\})$ such that  $a_{b_0}(\gamma)\ne 0$. Then $Z_b(\gamma)\ne 0$ for any $b\in B$ as long as $\gamma\in S_{b_0}$.
\end{cor}
{\it Proof of the Corollary.} Indeed, if $Z_b(\gamma)=0$ for some $b\in B$ then $\gamma$ cannot belong to  $Supp(\sigma_b)$ and hence it cannot satisfy the condition $a_b(\gamma)\ne 0$. But  the Lemma \ref{central charges on the line} implies that $a_b(\gamma)=a_{b_0}(\gamma)\ne 0$,  since $\gamma\in S_{b_0}$. This contradiction proves that $Z_b(\gamma)\ne 0$. The Corollary is proven.

Let us prove the Proposition \ref{support for families}. As above we will give the proof in the most difficult case when $rk\, Z_0=2$.

By choosing a basis and Euclidean norm $|\bullet|$ on $\Gamma_\R$ we may assume that $\Gamma_\R=\R^n$ with coordinates $(x_1,...,x_n)$ and $Z_0(x_1,...,x_n)=x_1+\sqrt{-1}x_2$. For any $\epsilon>0$ the quadratic form $Q_\epsilon(x_1,...,x_n)=x_1^2+x_2^2-\epsilon^2 (\sum_{3\le i\le n}x_i^2)$ is negative on $Ker\, Z_0$. Let us choose $\epsilon>0$ such that the $\epsilon$-ball with the center of $Z_0$ belongs to $U_0$.
We claim for any $\g$  as in the Proposition \ref{support for families}, the support of $\sigma_{Z_0}$ belongs to the subset $\{\gamma\in \Gamma-\{0\}|Q_\epsilon(\gamma)>0\}$. 
We will prove this fact assuming the converse and arriving to a contradiction.

Let $x\in \Gamma_\R=\R^n$ satisfies the condition $x\ne 0$ and  $Q_\epsilon(x)<0$. Let $x=x_h+x_v$ be the orthogonal decomposition into the ``horizontal" and ``vertical" parts, so that the projection $Z_0$ satisfies the property $Z_0(x_v)=0$.  Notice that $\epsilon|x_v|>|x_h|$.

We can identify $x_h$ with $Z_0(x_h)$. Let us assume that both $x_h$ and $x_v$ are non-zero vectors. We denote by $F$ the real plane spanned by $x_h$ and $x_v$ and by $F^\perp\subset \Gamma_\R$ the orthogonal subspace.  Consider a continuous family of linear maps $Z_t: \R^n\to \R^2\simeq \C, t\in [0,1]$ defined by the conditions that:

a) $Z_{t=0}=Z_0$;

b) the restriction of $Z_t$ to $F^\perp$ coincides with $Z_0$; 

c) $Z_t(x_h)=x_h, Z_t(x_v)=-tx_h$. 

Setting $t=1$ we see that $Z_1(x)=0$. 
Furthermore, $Z_t(F^\perp)=(\R\cdot x_h)^\perp\subset \R^2, t\in [0,1]$. In the induced norm $||\bullet||$ on $Hom(\R^n,\R^2)$ we have $||Z_t-Z_0||\le \epsilon$.  Hence $Z_t\in U_0$ for all $t\in [0,1]$.

Notice that the subspace $Z_t^{-1}(\R\cdot x_h)$ does not depend on $t$, so we can apply the Lemma \ref{central charges on the line}.
Let $\gamma\in Supp(\sigma_{Z_0})$ satisfy the condition $Q_\epsilon(\gamma)<0$. Then $\gamma,\gamma_h,\gamma_v$ are all non-zero (the non-triviality of $\gamma_v$ follows from non-triviality of $\gamma_h$ and the inequality $\epsilon |\gamma_v|>|\gamma_h|>0$). By Lemma \ref{central charges on the line} we see that $a_{\sigma_{Z_t}}(\gamma)$ does not depend on $t\in [0,1]$. Since $\gamma\in Supp(a_{\sigma_{Z_0}})$ we conclude that $a_{\sigma_{Z_t}}(\gamma)$ is non-zero for all $t\in [0,1]$. But $Z_1(\gamma)=0$, which contradicts to the Corollary of the Lemma \ref{central charges on the line}. This proves our Proposition \ref{support for families} in the rank $2$ case. Similar considerations give the proof in the rank $1$ case. The rank zero case is obvious. The Proposition \ref{support for families} is proven.

\subsection{Lie algebra of vector fields}\label{vector fields}

Let $\Gamma$ be a free abelian group of rank $n$ and $\g=Vect_{\Gamma}:=Vect_{\Gamma,\Q}=Der(\Q[\Gamma])=\oplus_{\gamma\in \Gamma}\g_{\gamma}$ be the graded Lie algebra of vector fields on the torus  $Hom(\Gamma,{\bf G}_m)$. The graded  component $\g_{\gamma}=Hom(\Gamma,\k)$ can be identified with the vector space spanned by $x^{\gamma}x_i\partial/\partial x_i, 1\le i\le n$ as long as we choose an isomorphism $\Gamma\simeq \Z^n$.
A choice of such an isomorphism allows us also to identify
the Lie algebra $Vect_{\Gamma}$ with  the $\Gamma$-graded Lie algebra $Vect_n$, consisting of derivations of the ring of Laurent polynomials $\Q[x_1^{\pm 1},...,x_n^{\pm 1}]$, where $n=rk\, \Gamma$.

\begin{rmk} If $\,$   $\Gamma$ carries a skew-symmetric integer bilinear form $\langle\bullet,\bullet\rangle:\bigwedge^2\Gamma\to \Z$.
then the Lie algebra $Vect_{\Gamma}$ contains a Lie subalgebra ${Ham}_{\Gamma}$ of Hamiltonian vector fields,
i.e.  vector fields of the form $\{f,\bullet\}$, where $f\in {\mathcal O}(Hom(\Gamma,{\bf G}_m))$
is a regular function on the Poisson torus of characters of $\Gamma$. After a choice of isomorphism $\Gamma\simeq \Z^n$
we can identify
$Ham_{\Gamma}$ with the Lie algebra $Ham_n\subset Vect_n$ of Hamiltonian vector fields on the Poisson torus ${\bf G}_m^n$.
The Poisson algebra $\OO({\bf T}_\Gamma)$ considered as a Lie algebra is isomorphic to the direct sum to its center (which is abelian Lie algebra) and its image in ${Ham}_{\Gamma}$.

Notice that for the Lie algebra of  Hamiltonian vector fields
all graded components are at most one-dimensional
(they correspond to the Hamiltonian vector fields  $\{x^{\gamma},\bullet\}, \gamma\in \Gamma$).
In the case of the Lie algebra of {\it all} vector fields we  have $dim\,\g_{\gamma}=n$.

\end{rmk}

For any field $\k$ of characteristic zero we will denote by $Vect_{\Gamma,\k}=Vect_\Gamma\otimes \k$ the corresponding Lie algebra of $\k$-valued vector fields on the torus. Sometimes we will skip $\k$ from the notation, if the ground field is clear.

\subsection{Wheels of cones and formal schemes arising from stability data}\label{wheels and schemes}

In this subsection work over an arbitrary ground field $\k$ of characteristic zero. Also we will use the notation for cones which is dual to those in the rest of the paper, but in agreement with those in [KoSo3].

Let us recall some definitions and results from [KoSo3], slightly changing the notation with respect to the loc.cit. We are going to consider the rank two case of stability data, i.e. we assume that the image of $Z_\R$ is the whole  plane $\R^2=\C$.
Let $J$ be a cyclically ordered finite set, i.e. $J\simeq \Z/N\Z$. We will assume that $N\ge 3$. Let $\Gamma$ be a free abelian group of rank $n$ as above.

\begin{defn} A wheel of closed cones $(C_j)_{j\in J}$ is a collection of closed rational strict convex cones $C_j\subset (\Gamma_\R)^\ast \simeq (\R^n)^\ast$ such that the interior $int(C_j)$ of each $C_j$ is non-empty and pairwise do not intersect, and the following conditions are satisfied:

1) $dim(C_j\cap C_{j+1})=n-1, j\in J$;

2) for any $j_1\ne j_2\in J$ the set $\{j|C_j\subset Conv(C_{j_1}\cup C_{j_2})\}$ is either the whole set $J$ and a connected interval of cyclically ordered integers in one of the two forms: $\{j_1,j_1+1,...,j_2\}$ or $\{j_2, j_2+1,...,j_1\}$.

Here $Conv$ denotes the convex hull of the set.

Dropping the condition that $J$ is cyclically ordered we arrive to the notion of the chain of closed cones.

\end{defn}

We impose the following version of the {\it Connectedness Assumption} from [KoSo3], Section 5.1:

{\it for an $j\in J$ the union $C_j\cup C_{j+1}$ is again a closed strict convex cone.}

It follows that the intersection $C_j\cap C_{j+1}$ is a face of both $C_j$ and $C_{j+1}$.

The {\it Non-degeneracy Assumption} from the loc.cit. says that the intersection of all dual cones $C_i^\vee$ consists of one point $0\in \Gamma_{\R}$.

We also impose the {\it Wheel Assumption}:

{\it $C_i\cap C_j=0$ unless $i\in \{j,j\pm 1\}$}.

Then considering the wheel of cones $C_j, j\in J$ including the faces of the cones as a fan we obtain a toric variety $X:=X_{(C_j)_{j\in J}}$ over $\k$. It is naturally stratified by the toric strata $X_s, dim\,X_s=s, 1\le s\le n$.  Also we have the toric divisor $D\subset X$ which is a normal crossing divisor given by the complement to the open orbit of the $n$-dimensional torus.

Consider  $1$-dimensional torus orbits $F_i\simeq {\bf G}_m$ as well as their closures $\overline{F}_i\simeq {\bf P}^1$. We denote by $X^{form}:=X^{form}_{(C_j)_{j\in J}}$ the formal scheme obtained as the completion of $X$ along $\overline{F}:=\cup_{i\in J}\overline{F}_i$. Let $U_i$ denotes the formal neighborhood of $F_i\cup \{p_i\}\cup F_{i+1}$,
where $p_i$ is the only intersection point of $\overline{F}_i$ and $\overline{F}_{i+1}$.
Then $U_{i,i+1}=U_i\cap U_{i+1}$ is the formal neighborhood of $F_{i+1}$. 
Thus we have the wheel of projective lines
$\overline{F}=\cup_{i\in J}\overline{F}_{i}$ both in $X$ and $X^{form}$ (cf. loc.cit. Def. 5.1.1). 

In term of cones, the algebra if functions on the formal neighborhood of $p_i$ consists of series in the monomials $x^\gamma$ with $\gamma\in \Gamma\cap C_i^\vee$.

\begin{defn} \label{compatible wheel} Given stability data $\sigma$ of rank $2$ on $Vect_\Gamma$ we say that the wheel of closed cones $\CC=(C_i)_{i\in J}$ is compatible with $\sigma$ if:

a) for any $i\in J$ the intersection $t_{i,i+1}^0=(C_i\cap C_{i+1}-\partial (C_i\cap C_{i+1}))\cap Z_\R^\ast(\R^2)$ is an open ray;

b) the cyclic ordering of the rays $l_{i,i+1}, i\in J$  agrees with the clockwise order with respect to the orientation on the plane $Z_\R^\ast(\R^2)$;

c)  support of the stability data belongs to $\cup_{i\in J}(C_i\cap C_{i+1})^\vee-\{0\}\subset \Gamma_\R-Ker\, Z_\R$.

\end{defn}

If the conditions a) and b) are satisfied for the  linear map $Z_\R=Z\otimes \R: \Gamma_\R\to \R^2$ we will say that the wheel of closed cones is {\it compatible with $Z$ (or with the stability data)}. 
By the Proposition 5.3.2 from loc. cit. for given stability data on $Vect_\Gamma$ there exists a compatible wheel of closed cones in $\Gamma_\R^\ast$.

Let $C_{i,i+1}=C_i\cap C_{i+1}$ and let $t_{i,i+1}$ denote the closure of $t_{i,i+1}^0$. Then we have a cyclically ordered set of closed rays $l_{i,i+1}=(Z_{\R}^\ast)^{-1}(t_{i,i+1})\subset \R^2$. They are boundaries of the cyclically ordered set of admissible sectors 
$V_i, i\in J$. Let $G_i$ be the proalgebraic group of automorphisms of the formal neighborhood $U_i$. Similarly we have groups $G_{i,i+1}$ of automorphisms of $U_{i,i+1}:=U_i\cap U_{i+1}$.  Clearly $G_{i}$ (resp. $G_{i+1}$ acts from the left (resp. from the right) on $G_{i,i+1}$ via the natural embedding. 

We remark that if the central charge $Z$ is rational then all rays $l_{i,i+1}$ and sectors $V_i$ are rational.

\begin{rmk}
Notice that the above definition can be inverted in the  following way. Suppose we are given the wheel of closed cones $\CC=(C_i)_{i\in J}$. We say that the stability data $\sigma\in Stab(Vect_\Gamma)$ of rank $2$ with the central charge $Z$ are compatible with $\CC$ if $\CC$ is compatible with $Z$ and the conditions a)-c) of the Definition \ref{compatible wheel} are satisfied.
We are going to use this terminology without further comments. We will denote the set of stability data compatible with the wheel of closed cones $\CC$ by $Stab_{\CC}(Vect_\Gamma)$ or simply by $Stab_\CC$.

\end{rmk} 

Recall the  $\k$-formal scheme ${X}_{\overline{F}}^{form}$ corresponding to the wheel of projective lines $\overline{F}$. Notice that for the formal completion $X_{p_i}^{form}$ of the point $p_i$ we have:
$Aut({X}_{p_i}^{form})\supset G_{i, i+1}$ and similarly for the completion at $\overline{F}_i$ we have $Aut({X}_{\overline{F}_i})^{form}\supset G_{i}$. In fact $G_{i, i+1}$ is a kernel of the natural epimorphism of $Aut({X}_{p_i}^{form})$
onto the torus $Hom(\Gamma,{\bf G}_m)\simeq ({\bf G}_m)^n$. The torus is interpreted as the group of rescaling of coordinates.

\begin{rmk} Our version of Connectedness Assumption means that we can contract any projective line $\overline{F}_{i}$ to a point obtaining a new  formal scheme with a new closed point which is a toric singularity
(i.e. the completed local ring of the point is a completion of the polynomial algebra in the closed rational strict convex cone). Notice that the $\widehat{X}_{p_i}$ are already in this form.

\end{rmk}

It follows from the construction of the formal scheme that there is  a natural action of the group 
$(\k^\times)^n$ on the formal scheme given by rescaling of variables. The {\it decoration of the formal scheme} is a trivialization of the corresponding $(\k^\times)^n$-torsor.  Arising {\it decorated formal scheme} will be sometimes called {\it decorated formal wheel}.  An equivalent description of the decoration can be found in [KoSo3], Def. 5.2.2.

We also recall the following relation between the description of stability data in terms of
wheels of sectors (see [KoSo3], Proposition 5.3.3).

\begin{prp} \label{sectors and cones} For given stability data on $Vect_{\Gamma, \k}$  there exist an admissible wheel of cones $(C_i)_{i\in \Z/m\Z}$ compatible with the stability data as well as the following:

1) a cyclic decomposition $\R^2=V_{1,1}\cup...\cup V_{1,k_1}\cup V_{2,1}\cup...\cup V_{2,k_2}\cup...\cup V_{m,1}\cup...\cup V_{m,k_m}, m\ge 3, k_i\ge 1$, where $V_{i,j}$ are closed strict sectors such that two consecutive sectors have a common edge, $Z^{-1}(\partial V_{i,j}-\{0\})\cap \Gamma=\emptyset$;

2) a cyclically ordered collection of closed strict convex cones $C(V_{i,j})\subset \Gamma_{\R}$ compatible with $Z$ and such that $Z(C(V_{i,j}))\subset V_{i,j}$, the set $Supp(a)$ belongs to $\cup_{i,j}C(V_{i,j})$, and  for any $i,j$ the set $C(V_{i,j})-\{0\}$ belongs to $int(C_{i,i+1}^{\vee})$.
\end{prp}

Proof of the following result is similar to the one of Proposition 5.2.7 from [KoSo3]. It is also a special case of the approach to stability data proposed in Subsection \ref{alternative def}

\begin{prp}
\label{formal coset}

Suppose we are given stability data $\sigma$ on $Vect_{\Gamma,\k}$ of rank $2$ such that $Z_\R(\Gamma_\R)=\C$.  Let $(C_j)_{j\in J}$ be a wheel of closed cones such that $(Z_\R^\ast)^{-1}(int(C_j))$ are disjoint non-empty cyclically ordered admissible sectors in $(\R^2)^\ast$. Let us consider the subset $Stab_{\sigma,(C_j)_{j\in J}}(Vect_\Gamma)\subset Stab(Vect_\Gamma)$ satisfying the following property:

$\bullet$ the set of $\gamma\in \Gamma$ such that $Z(\gamma)\ne 0$ and $\gamma$-component of $log(A_{\R_{>0}\cdot Z(\gamma)})\ne 0$ belongs to $\cup_{j\in J}C_j^\vee$.

Then $Stab_{\sigma,(C_j)_{j\in J}}(Vect_\Gamma)$ is in one-to-one correspondence with the homogeneous space $\MM_{(C_i)_{i\in J}}:=...\times G_{i,i+1}\times_{G_{i+1}}G_{i+1,i+2}\times_{G_{i+2}}G_{i+2,i+3}\times...:=\prod_{i\in J}G_{i,i+1}\times_{G_{i+1}}G_{i+1,i+2}=\prod_{i\in J}G_{i, i+1}/\prod_{i\in J}G_i$.

\end{prp}
Notice that the natural actions of the group $G_{i+1}$ on $G_{i,i+1}$ and  $G_{i+1,i+2}$ is free for each $i\in J$.

We will call $\MM_{(C_i)_{i\in J}}$ the {\it moduli space of stability data compatible with the stability data $\sigma$ and the wheel of cones $(C_i)_{i\in J}$}.

\begin{rmk} In terms of the above-mentioned quadratic forms $Q$ one can say that in the rank two case for any quadratic form $Q$ on $\Gamma_\R$ of signature $(2,n-2)$ whose
restriction to $Ker(Z_\R)-\{0\}$ is strictly negative, there exists of wheel of closed cones $(C_j)_{j\in J}$ such that $\{u\in \Gamma_\R|Q(u)>0\}$ belongs to $\cup_{j\in J}C_j^\vee$.

\end{rmk}
Similarly to [KoSo3] one can show that the moduli space $\MM_{(C_i)_{i\in J}}$
is an affine scheme.

In the notation of the Proposition \ref{formal  coset} we will denote the formal scheme corresponding to an element $g=(g_i)_{i\in J}\in G:=\prod_{i\in J}G_{i,i+1}$ (cyclic product) by $X_{g, (C_i)_{i\in J}}^{form}$ or simply $X_g^{form}$ if it will not lead to a confusion.  If $g$ and  $h$ give rise to the different points of $\MM_{(C_i)_{i\in J}}$ then $X_g^{form}$ is non-isomorphic to $X_h^{form}$. At the same time there might be formal schemes  containing the same wheel of projective lines $\overline{F}$ as the set of closed points, but which are not of the type $X_g^{form}$.

Recall Remark 5.2.8 from [KoSo3]. Adopted to our case of the Lie algebra $Vect_{\Gamma, \k}$ it says that if we have another admissible wheel of cones $(C_i^\prime)_{i\in J^\prime}$ such that for any $j\in J^\prime$ there exists $i\in J$ such that $(C_j^\prime)\subset C_i$  then there is a natural embedding $\MM_{(C_i)_{i\in J}}\to \MM_{(C_i^\prime)_{i\in J^\prime}}$.  
At the level of formal schemes $X_g^{form}$ such subdivisions of cones correspond to two operations at the level of formal schemes (see  loc.cit):

a)  one makes a finite sequence of blow-ups of $X^{form}_g$ with centers at some toric strata;

b) in the formal scheme resulting from the blow-ups in a) one chooses a wheel of projective lines formed by $1$-dimensional strata and takes the completion along the wheel.

In this way one gets  an inductive system with respect to the above subdivisions  (equivalently, operations a) and b) on formal toric varieties).
The following result gives  an alternative description of the set of stability data as in Subsection \ref{alternative def}. It  is also the analog of the Theorem 5.3.4 from [KoSo3]. The proof is similar.

\begin{thm}\label{indlim}
The set of all stability data of rank $2$ on $Vect_\Gamma$ can be identified with the inductive limit $\varinjlim \MM_{(C_i)_{i\in J}}$ taken with respect to the inductive structure on the wheels of cones given by the above-mentioned subdivisions.
\end{thm}

Finally let us describe the set of stability data of rank $1$.
We remark that in the rank $1$ case stability data give rise to a pair of  elements $A_{l_\pm}\in G_{l_\pm}$ corresponding to two admissible opposite rays $l_\pm\subset \R^2$. They determine the stability data uniquely.

\begin{rmk}\label{local system of data} The above considerations can be generalized to the case when one has a local system $\underline{\Gamma}$ of lattices over $\R^2-\{0\}$, local system of $\underline{\Gamma}$-graded Lie algebras $\g_{\underline{\Gamma}}$ and a locally constant family of central charges $\underline{Z}:\underline{\Gamma}\to \C$. Then in the case of the Lie algebra of vector fields one can allow the local system $\underline{\Gamma}$ to have non-trivial monodromy on the wheel of projective lines. This ``twisted' version of stability data can be thought of as a special case of wall-crossing structure. It  appears in practice, e.g. in the case of complex integrable systems.

\end{rmk}

\section{Analytic stability data }\label{analytic stab data}

\subsection{Analytic transformations}\label{analytic transformations}

From now on we will assume that  $\k=\C$ unless we say otherwise.\footnote{ The definition and some of the results below make sense also for $K=\C((t))$, or an algebraic extension of the latter.}

Let us denote by $x^\gamma, \gamma\in \Gamma$ the monomial function on the torus ${\bf T}_\Gamma=Hom(\Gamma, \C^\ast)$, such that $x^\gamma(\phi)=\phi(\gamma)\in \C^\ast$. Clearly $x^\gamma x^\mu=x^{\gamma+\mu}$.

If $C \subset \Gamma_\R$ is a strict convex cone then it gives rise to the completion 
$\OO_{C}({\bf T}_\Gamma)$ of the algebra of regular functions on ${\bf T}_\Gamma$. Namely $\OO_{C}({\bf T}_\Gamma)$ consists of series 
$f=\sum_{\gamma\in \Gamma\cap C}c_\gamma x^\gamma, c_\gamma\in \C$. If $\phi\in {\bf T}_\Gamma$ then we have the series $f(\phi)=\sum_{\gamma\in \Gamma\cap C}c_\gamma x^\gamma(\phi)$. 
\begin{defn}\label{analytic series}
In case if the series is absolutely convergent for some $\phi$ (and hence convergent on an open set) we will call $f$ {\it analytic}. This property is equivalent to the exponential bound on the coefficients $c_\gamma$, i.e. there exists $R>0$ such that $|c_\gamma|\le R^{||\gamma||}$. Here $||\bullet||$ is any fixed norm on $\Gamma$.
\end{defn}

Recall that given a graded Lie algebra $\g=\oplus_{\gamma\in \Gamma}\g_\gamma$ and a strict convex cone $C\subset \Gamma_\R$ with non-empty interior, we can construct a pronilpotent group $G_{C}$ associated with the pronilpotent Lie algebra $\g_{C}=\prod_{\gamma\in \Gamma\cap C}g_\gamma$. 

\begin{defn} Let $g\in G_{C} $. We say that $g$ is analytic if  $g(x^\gamma)=x^\gamma(1+...)$ is
analytic for all $\gamma\in \Gamma$
\end{defn}

\begin{prp} It suffices to check the analyticity of $g$ for $g(x^\gamma), \gamma\in C$.

\end{prp}

{\it Proof.} Notice that for any $\gamma\in \Gamma$ we have: $g(x^\gamma)=x^\gamma(1+\sum_{\mu\in \Gamma\cap C}c_{\gamma\mu}x^\mu)$.  Since $C$ has non-empty interior, we can find a basis $\gamma_1,...,\gamma_n$ of $\Gamma$ which belongs to $\Gamma\cap C$.  If $g(x^{\gamma_i})$ is analytic then for any
$\gamma=\sum_{1\le i\le n}a_i\gamma_i$ we have $g(x^\gamma)=g(x^{\gamma_1})^{a_1}...g(x^{\gamma_n})^{a_n}$ is analytic as the product of analytic series. This proves the Proposition.
$\blacksquare$
 
 Notice that the property to be analytic does not depend on the choice of a basis in $\Gamma$, which can give different explicit formulas for the series $g(x^\gamma)/x^\gamma$ if written in coordinates.
 Thus the analyticity of $g(x^\gamma)$ can be effectively verified by checking exponential bounds on coefficients of the series in the RHS.

\begin{rmk}
For each admissible open sector $V$  the group $G_V$ contains a subgroup $G_V^{an}$ consisting of analytic elements. An element $g\in G_V$ belongs to $G_V^{an}$ iff there exists a closed strict convex cone $C\in \Gamma_\R, Z(C)\subset V$ and an element $g_C\in G_C$ which is  analytic, and which is mapped to $g$ under the natural map $G_C\to \varinjlim_{C}G_C=G_V$.

Notice that in the same way one can define $G_V^{an, old}$ for the version of $G_V^{old}$ from [KoSo1]. In all the above definitions the sector $V$ in fact does not have to be open, It can be closed or semiclosed. 
\end{rmk}

\subsection{Analyticity of factors: two-dimensional case}

We are going to consider only rational central charges.  Our aim is to prove the following result.

\begin{thm}\label{analyticity of factors}
Suppose we are given  stability data of rank $2$ with rational central charge.
If $V$ is a closed rational admissible sector and $V=V_1\cup V_2$ is its decomposition into the union of closed rational admissible sectors with the common boundary ray $l_{12}$, then $G_V^{an}=G_{V_1}^{an}G_{V_2}^{an}$.

\end{thm}

In this subsection we consider the two-dimensional case only. We will finish the proof in the next subsection.

In the two-dimensional case we may assume that $Z=id$, hence rationality condition is empty.

After several reductions the result will follow from its special case, when $V$ is the first quadrant $\{x\ge 0, y\ge 0\}$, and $l_{12}=l$ is the diagonal ray $\{x=y\ge 0\}$. We order $V_1:=V_+$ and $V_2:=V_-$ clockwise, i.e. $V_+$ is the subsector of $V$ above $l$ and $V_-$ is the one below $l$.

Let $X={\bf P}^1\times {\bf P}^1$ be the toric surface and $\widehat{X}$ be its blow-up at the points $x_+=(0,\infty)$ and $x_-=(\infty, 0)$.  The toric surface $\widehat{X}$ contains six smooth components of the toric divisor: the exceptional divisors $D_\pm$ corresponding to $x_\pm$, as well as proper transforms $D_x, D_y, D_1, D_2$ of the projective lines ${\bf P}^1\times \{0\}, \{0\}\times {\bf P}^1, {\bf P}^1\times \{1\}, \{1\}\times {\bf P}^1$ respectively. We keep same notation for proper transforms. Let $U_x\subset \widehat{X}$ and $U_y\subset \widehat{X}$ denote    formal neighborhoods of $D_x$ and $D_y$ respectively. Then the union $U=U_x\cup U_y$ is a formal surface stratified by exceptional divisors $D_\pm$ old divisors $D_x, D_y$ and their intersection point $(0,0)=D_x\cap D_y$. 

The group element $g_V\in G_V$ gives rise  the formal automorphism of $U_x\cap U_y$  such that  $g_V(x)=x(1+...), g_V(y)=y(1+...)$. Hence we can use $g_V$ in order to glue the new formal scheme $U^{new}$ which is still stratified by 
$D_x, D_y, (0,0)$, since $g_V$ preserves these strata.  The formal scheme $U^{new}$ contains also strata $D_\pm$ (here we abuse notation for formal schemes corresponding to algebraic divisors). Let $V_\pm^0=V_\pm-l$ denote the semi-open subsectors of $V_\pm$. Then $g_V=g_{V_+^0}g_lg_{V_-^0}$. Furthermore $g_{V_+}=g_{V_+^0}g_l, g_{V_-}=g_lg_{V_-^0}$. Here $g_{V_\pm}\in G_{V_\pm}^{old}$.

Consider the groups $Aut(U_x)$ and $Aut(U_y)$ of automorphisms of formal neighborhoods  which  preserves the induced toric stratifications by the toric divisors in the neighborhood. There are homorphisms $Aut(U_x)\to \C^{\ast}\times \C^{\ast}$ and $Aut(U_y)\to \C^{\ast}\times \C^{\ast}$ which at the level of Lie algebras, in coordinates, are projections to the Lie algebras spanned by the vector fields $x\partial_x, y\partial_y$ (at the level of groups this Lie algebra corresponds to the group of dilations). Let $Aut^{(1)}(U_x)$ and  $Aut^{(1)}(U_y)$ denote the kernels of the above homomorphisms of groups. We have similarly defined group $Aut^{(1)}(U)$.

\begin{lmm} \label{groups}We have the following isomorphisms of groups 
$$G_{V_+}^{old}\simeq Aut^{(1)}(U_y), G_{V_-}^{old}\simeq Aut^{(1)}(U_x).$$ 

Also, we have an isomorphism of groups $G_l\simeq Aut^{(1)}(U)$. The 
latter group  is also isomorphic to the semidirect product $Aut(D_\pm)\ltimes H_\pm$, where $H_\pm$ is the group of changing trivializations of the normal bundle to $D_\pm$ (its Lie algebra can be called Atiyah algebra of the normal bundle).

\end{lmm}

{\it Proof.} Notice that elements of $G_{V_+}^{old}$ are exactly  automorphisms preserving the formal series  $\sum_{n,k\ge 0}a_{nk}x^n(xy)^k$ (up to dilations) and hence are automorphisms of the stratified formal toric surface $U_y$. Similarly for $U_x$. 
The result about $G_l$ follows from these two observations.
$\blacksquare$

It follows that $U^{new}$ is isomorphic to $U$ as a stratified formal scheme.



Let us now consider  family of hyperbolas $C_\epsilon$ given by the equations $\{xy=\epsilon\}\subset {X}, \epsilon\ne 0$. After the blow-up we obtain the corresponding family of curves on $\widetilde{X}$, for which we will keep the same notation. Their natural compactifications $\overline{C}_\epsilon$ are isomorphic to ${\bf P}^1$. Then $\overline{C}_\epsilon$ intersects $D_+$ transversally at the point $0_\epsilon$ and $D_-$ at the point $\infty_\epsilon$. As $\epsilon\to 0$ the curves degenerate to $D_x\cup D_y$. Furthermore the family of smooth projective curves $\overline{C}_\epsilon$ gives rise to a versal deformation of $D_x\cup D_y$. Notice that by construction the normal bundle to each curve is trivial. Same is true for each of the divisors $D_\pm$. 

The group $G_{V_\pm^0}$ is the group of changes of analytic coordinates which preserves the trivialization of normal bundle to $D_\pm$. Indeed, for, say $g\in G_{V_+^0}$ we have $g: x\mapsto x(1+r(x,y))$, where $r(x,y)$ is a series in $x^iy^j, i>j$. The coordinates on the blow-up in the formal neighborhood of $D_+$ are $(xy, y)$. Since $r(x,y)$ does not have terms $(xy)^i$ the automorphism $g$ preserve the trivialization of the normal bundle to $D_\pm$ which agrees with the trivialization of the tangent bundle to $\overline{C}_\epsilon$ at $0_\epsilon=\overline{C}_\epsilon \cap D_+$. Of course similar considerations apply to $D_-$ and $G_{V_-^0}$ but we don't need them. 

Choice of trivializations as above gives us the coordinates $(\epsilon, t)$, where $\epsilon$ can be now interpreted as the formal coordinate along the formal divisor $D_+$, while $t$ is the coordinate along the algebraic curve $\overline{C}_\epsilon\simeq {\bf P}^1$.  It is easy to see that these coordinates are defined up to a transformation from the group $G_l$.  

Now we are ready to finish the proof of our special case of the Theorem \ref{analyticity of factors}. If $g_V$ is analytic as a transformation on $(x,y)$ then it follows that the formal neighborhoods $U^{new}$ is in fact a formalization of the analytic neighborhood $U^{new, an}$ of $D_x\cup D_y$. The formal family of curves $\overline{C}_\epsilon$ is in fact an analytic family (this is a special case of Douady theorem). Thus $(\epsilon, t)$ is an analytic coordinate system in $U^{new,an}$. Since $G_{V_+}$ consists of changes of these analytic coordinates preserving the trivialization of the normal bundle to $D_+$, all $g\in G_{V_+}$ are in fact analytic transformations. Similarly for $G_{V_-}$. In order to prove analyticity of the transformations of $g\in G_l$ we notice the following general result.

\begin{lmm} \label{analyticity of diagonal} If $l$ is the boundary ray for the closed admissible sector $V_+$, then analyticity of $g_{V_+}$ implies analyticity of $g_l$.
\end{lmm} 

{\it Proof.} Indeed, $g_l(x^\gamma)=x^\gamma(1+...)$, where the dot terms have monomial $x^\mu, \mu\in \Gamma\cap C(l)$. But in the convergent series $g_{V_+}(x^\gamma)/x^\gamma$ the terms with monomials $x^\mu, \mu\in \Gamma\cap C(l)$ can come from $g_l$ only. This proves the Lemma and the special case of the Theorem \ref{analyticity of factors} in the case when $V$ is the first quadrant in $\R^2$. $\blacksquare$

Let us return to the case of more general ray. Namely, consider the case when $V$ is the first quadrant as before, but $l_1$ is the ray in $V$ which has the slope $b/a$ for a pair of relatively prime positive integers. Notice that if $l_1$ bisects $V_+$ or $V_-$ the analyticity can be deduced from the special case $b/a=1$ considered above. Indeed, choosing an integer basis $\gamma_1, \gamma_2$ in, say, $V_-$ such that $\gamma_i, i=1,2$ are primitive vectors $(1,1), (1,0)$, we can identify $V_-$ with the first quadrant and $l$ with the bisecting ray. Now we have two rays $l, l_1$  in $V$ with analyticity of $g_{V_i}$ for each of the three arising sectors $V_i$ as well as $g_V$. We can continue this process organizing the  slopes as elements of the Farey sequence. This gives the proof for an arbitrary slope $b/a$. Furthermore, since any closed rational sector can be identified with the first quadrant by choosing primitive basis of $\Z^2$ belonging to the boundary rays, the result holds for any such sector and any rational ray dissecting it. This concludes the proof of the Theorem \ref{analyticity of factors} in the case when $rk\, \Gamma=2$. $\blacksquare$

\subsection{Higher-dimensional case}
We will finish the proof of the Theorem \ref{analyticity of factors} in this subsection. We start with the general result.

Let $Z: \Gamma\to \C$ be a central charge such that $Z(\Gamma)=\Z^2$, and $C$ be a closed rational strict convex cone $C\subset \Gamma_\R$  such that 
$Z_\R(C)=V$ is a closed rational strict sector. 

By choosing a primitive basis $\gamma_1,...,\gamma_n$ of $\Gamma$ we may 
assume  that  $\Gamma=\Z^n$, the cone $C$ is given by the inequalities $x_i\ge 0, 1\le i\le n$.
Let $l\subset V$ be a rational ray of the slope $b/a$ with the vertex at $0$. We denote by $V_+$ and $V_-$ the corresponding subsectors of $V$ (see previous subsection) and by $\pi=Z^{-1}(l)$ the rational cooriented hyperplane dissecting $C$ into two closed cones $C_\pm=Z^{-1}(V_\pm)$.
Since the central charge is rational, we may assume that  $Z=\sum_{1\le i\le n}a_ix_i$, where $a_i\in Z$. We may assume that $a_1>0, a_2<0$. Let $d=g.c.d. (a_1,a_2)$.

\begin{prp}
 Using transformations from $Aut(\Gamma)$ we can bring $Z$ to the form  $ax_1-bx_2$, where $a,b>0$ are integers.
\end{prp}
{\it Proof.} If $d=1$ then any $a_i, i\ge 3$ is an integer linear combination of $a_1$ and $a_2$ and hence can be made equal to zero. Otherwise we can make all $a_{\ge 3}$ non-negative integers strictly smaller than $d$.  In particular $a_i<a_1, i\ge 3$. Hence by an integer change of variables we can replace $a_1$ by $a_1^\prime=a_1-a_3>0$ without changing all other $a_i$'s. Thus we can make $a_1$ smaller. Applying the above considerations to $a_1^\prime, a_2$ we see that we can make all $0\le  a_{\ge 3}<d$ smaller than any positive integer, hence equal to zero. The Proposition is proved. $\blacksquare$

Let us now finish the proof of Theorem \ref{analyticity of factors}.

Let  $\Delta$ be the formal disc in the subspace with coordinates $x_3,...,x_n$ and $\widehat{Y}=\widehat{X}\times \Delta$, where $\widehat{X}$ is the toric surface from the previous subsection. If $0_\Delta$ is the center of $\Delta$ then we have formal neighborhoods $U_{D_x}=D_x\times \Delta$ and $U_{D_y}=D_y\times \Delta$ of the projective lines $D_x\times 0_\Delta$ and $D_y\times 0_\Delta$ respectively. Similarly $U_{D_x\cup D_y}=(D_x\cup D_y)\times \Delta$ is the formal neighborhood of the $(D_x\cup D_y)\times 0_\Delta$. All these formal neighborhoods are stratified by the toric strata (toric divisors and their intersections).
For each formal neighborhood $U$ we denote by $Vect_U$ the Lie algebra of vector fields which are tangent to the toric strata with order of at least $1$. Identifying a toric stratum with the coordinate stratum given by the equations $x_{i_1}=...=x_{i_k}=0$ we see that the vector fields tangent to it has the form $\sum_{j\in \{i_1,...,i_k\}}a_j(x)x_j\partial_{x_j}$, where $a_j(x)=\sum_{m\in \Z_{>0}^n}a_{j,m}x^m$ are formal series.

Then we have a family of curves $C_\epsilon\times \Delta\subset \overline{C}_\epsilon$. This gives a versal deformation in the formal neighborhoods of $C_\epsilon\times 0_\Delta$ and $\overline{C}_\epsilon\times 0_\Delta$ respectively. In other words the deformation theory of the $2$-dimensional case lifts trivially along the rest of the directions. Hence the factorization result from the previous subsection says that the series $g_{V_\pm}(x^\gamma)/x^\gamma$ is tautologically analytic in the variables $x_3,...,x_n$ and analytic in the variables $x_1, x_2$. This concludes the proof. $\blacksquare$

\begin{cor} \label{analyticity and rays}Suppose in the Proposition \ref{group factorization} the boundary rays of $V$ and all rays $l_i, 1\le i\le k$ are rational and $Z(\Gamma)=\Z^2$. If $g\in G_V$ is analytic then the factors $g_{V_0^+},g_{l_1}..., g_{V_k^-}$ in the factorization formula are also analytic.

\end{cor}
{\it Proof.} We can replace $V$ by a closed strict admissible rational subsector $V^\prime$. Then the result follows from Theorem \ref{analyticity of factors}. $\blacksquare$

\subsection{Analytic stability data}\label{analytic stab}

The  easiest way to  define the notion of analytic stability data is to the approach of Subsection \ref{alternative def}.

\begin{defn} \label{def analyticity} Stability data $\sigma\in Stab(Vect_\Gamma)$ with the central charge $Z$ are called analytic if there exists an open cyclic cover $(V_i)_{i\in J}$ and a representative $(g_{V_i})_{i\in J}\in \prod_{i\in J}G_{V_i}$ of $\sigma$ such that all elements $g_{V_i}$ are analytic.

\end{defn}

We denote by $Stab^{an}(Vect_\Gamma)$ the subset of analytic stability data.
The above Definition says that $\sigma$ is analytic if there exists a representative of $\sigma$ such that all $g_{V_i}$ belong to the corresponding subgroups $G_{V_i}^{an}$. The above definition covers stability data of all ranks.
The following result gives a description of rational analytic stability data.

\begin{prp} \label{analyticity rational}Let $\sigma \in Stab(Vect_\Gamma)$ be rational and $l_i, i\in J$ be a cyclically ordered finite collection of rational admissible rays such that any consecutive pair $l_i, l_{i+1}$ bounds an open strict admissible sector $V_i$. Then $\sigma$ is analytic iff all $A_{V_i}$ and all $A_{l_i}$ are analytic.

\end{prp}
{\it Proof.} The proof uses the Corollary \ref{analyticity and rays}   instead  of \ref{group factorization}.
Otherwise it is completely parallel to the proof of Proposition \ref{independence of cyclic cover}, where we replace the groups $G_{V_i}$ and $G_{V_i\cap V_j}$ by their analytic subgroups. $\blacksquare$

\begin{rmk} \label{independence of cover} Assume that $\sigma$ is rational. If $\sigma$ admits a representative consisting of analytic elements for an open rational cyclic cover $(V_i)_{i\in J}$  then the same is true for any other rational cyclic cover $(V_i^\prime)_{i\in J}$ . This follows immediately from the definition of analyticity.

\end{rmk}

\begin{conj} Stability data on $Vect_\Gamma$  are  analytic iff for every admissible sector $V\subset \R^2$ the corresponding group element $A_V\in G_V$ is analytic.

\end{conj}
As we have seen this Conjecture holds for rational stability data.


\subsection{Open and closed property of analytic stability data}\label{open and closed}

Here we will use the discussion of Subsection \ref{alt topology} in order to prove that the subspace of analytic stability data is open and closed (and hence is a disjoint union of connected components).

\begin{prp}\label{open} For any analytic stability data $\sigma$ there exists an open neighborhood of $\sigma$ in the space $Stab(Vect_\Gamma)$ consisting of analytic stability data only.

\end{prp}
{\it Proof.} Let us fix  a cyclic cover $(V_i)_{i\in J}$ such that the corresponding representative elements $g_{V_i}$ are analytic. Then there is an open neighborhood consisting with stability data represented by the same collection $(g_{V_i})_{i\in J}$ and sufficiently closed central charges. By definition all these stability data are analytic. $\blacksquare$

\begin{prp}\label{closed} If $\sigma\in Stab(Vect_\Gamma)$ belongs to the closure of $Stab^{an}(Vect_\Gamma)$ then it is analytic.
\end{prp}
{\it Proof.} Since $Stab^{an}(Vect_\Gamma)$ is open the stability data $\sigma$ with the central charge $Z$ is the limit of a sequence of rational analytic stability data $\sigma_n$ with the central charges $Z_n$ as $n\to \infty$. 
Let now $Q$ be a non-zero quadratic form such that $Q\ge 0$ on $Supp(\sigma)$ and $Q$ is negative on $Ker\,Z$. Then for all sufficiently large $n\ge n_0$ the form $Q$  satisfies the same property for $\sigma_n$.  

Let us choose  a rational cyclic cover $(V_i^\prime)_{i\in J}$ by four semiclosed quadrants in $\R^2$ with the origin $0\in \R^2$ and the right boundary ray removed from each $V_i^\prime$.  We may assume by Proposition \ref{analyticity rational} (which holds for semiclosed cyclic cover as well) that  all corresponding elements $A_{V_i}^{(n)}$ are analytic.

Let us choose $n_0$ in such a way that $Z(Z_{n_0}^{-1}(\overline{V}_i)\cap \{Q\ge 0\})$ belongs to an open half-plane. Here $\overline{V}_i$ denotes the corresponding closed quadrant.
Each of the four closed quadrants belongs to the interior of an open half-plane bounded by the line orthogonal to the  bisector of the quadrant. E.g. for the quadrant $\{x\ge 0, y\ge 0\}$ it is the open half-plane bounded by the line $x+y=0$. Let $\alpha_i, i\in J$ be these cyclically ordered open half-planes. Then the collection $(A_{V_i}^{(n_0)})_{i\in J}$ of the elements of $G_{V_i}^{old,an}$ gives a representative  of the rational analytic stability data associated with an open cyclic cover
$(\alpha_i)_{i\in J}$. By definition of the topology this collection gives a  representative of  stability data $\sigma$ associated with the open cyclic cover $(\alpha_i)_{i\in J}$. Hence $\sigma$ is analytic.
$\blacksquare$


\subsection{Analytic stability data via analytic germs}\label{analytic germs}

We start with a general definition.

\begin{defn}\label{analytic germ}
Let $\widehat{X}_Y$ be a formal scheme over $\C$, which is obtained by the completion of an algebraic variety $X$ along a closed algebraic subset $Y$.
Then  an {\it analytification of $\widehat{X}_Y$ }  is a pair consisting of an analytic germ $X^{an}_Y$ at the closed subset $Y$ and an isomorphism $\phi$ of the formal completion of $X^{an}_Y$ along $Y$ and $\widehat{X}_Y$. Here we abuse notation and denote by $Y$ both algebraic and analytic closed sets. 
\end{defn}

 In general the analytification is non-unique. E.g. taking $Y=\{pt\}$ and $\widehat{X}_Y=Spf(\C[[z]])$ we see that the analytifications form a set isomorphic to $Aut(\C[[z]])/Aut(\C\{z\})$. In general the isomorphism between two analytifications is unique, hence they form a groupoid which is a set (no automorphisms of points).

We will use the notation of Subsection \ref{wheels and schemes} which are dual to the notation in other sections. We consider stability data of rank $2$.
Let us fix a wheel of rational closed cones $\CC=(C_i)_{i\in J}$.
An element $g\in G$ (see loc.cit.) gives rise to the formal scheme  $X^{form}_{g, (C_i)_{i\in J}}$  derived from  $(C_i)_{i\in J}$ and $g\in G$, called formal wheel. If we fix just $\CC$ the set of formal schemes of the type $X^{form}_{g, (C_i)_{i\in J}}, g\in G$ is in bijection with the set of stability data $\sigma\in Stab(Vect_\Gamma)$ such that $Supp(\sigma)\subset \cup_{i\in J}C_i^\vee$.
If $X^{form}_{g, (C_i)_{i\in J}}$ admits analytification in the above sense, then the corresponding analytic germ at the wheel of projective lines $\overline{F}$ will be called {\it  analytic wheel}. The complete proof of the following result will take too much space, so we give only main steps. The interested reader can fill the details.

\begin{thm}\label{components vs germs} Rational stability data $\sigma$ of rank $2$ bounded by the wheel of rational closed cones $\CC$ are analytic iff  there exists $g\in G$ corresponding to $\sigma$ such that the  formal scheme $X^{form}_{g, (C_i)_{i\in J}}$ admits analytification. Moreover this analytification is unique.
\end{thm}

{\it Sketch of the proof.} If $\sigma$ is analytic then we take the collection $A_{V_i}, i\in J$ with $V_i=C_i^\vee\cap Z^{\ast}(\R^2)$ as the representative $g$. Using analytic automorphisms $A_{V_i^0}, i\in J$ we can glue  the desired analytic space. 

Conversely, assume that there exists an analytification. We want to prove that the stability data are analytic, i.e. the automorphisms $A_{V_i}$ are analytic. Here one uses the analytic coordinates coming from the two families of rational analytic curves as in the proof of the Theorem \ref{analyticity of factors}. The  difference with the latter is that we consider degenerating families of compact rational curves in the neighborhood of a toric singularity. 
In order  to show that a particular $A_{V_{i_0}}, i_0\in J$ is analytic we consider two families of straight lines parallel to the boundary rays of $V_{i_0}$ and ``lift" them to the above-mentioned two formal families of rational  curves. By Douady theorem the formal families of compact analytic curves are in fact analytic families. Finally, we interpret the automorphism $A_{V_{i_0}^0}$ corresponding to the interior $V_{i_0}^0\subset V_{i_0}$ as the change of these analytic coordinates preserving the transversal algebraic divisors determined by the directions of boundary rays for $V_{i_0}$ . 

Uniqueness of analytification follows from the fact that having an arbitrary analytification we interpret it as a result of gluing by analytic automorphisms $A_{V_i^0}$. $\blacksquare$

Recall the notions of decorated formal scheme and decorated formal wheel. We have the corresponding notions in the analytic setting. Since the group of automorphisms of the decorated formal wheel is trivial, the corresponding decorated analytic wheel is defined {\it canonically}.

\subsection{Stability data with gaps}

Let us discuss analyticity for  {\it stability data with gaps.} This is a generalization analytic stability data of rank $1$.

\begin{defn} We say that stability data of a graded Lie algebra $\g=\oplus_{\gamma\in \Gamma}\g_\gamma$ is the stability data with gaps
if there is a disjoint union $\sqcup_{1\le i\le m} V_i$ of closed strict sectors in $\R^2$ such that its complement is a disjoint union of strict open sectors and such that for each $\gamma\in \Gamma-\{0\}$ with $a(\gamma)\ne 0$ we have: $Z(\gamma)$ belongs to one of $V_i$'s.

\end{defn}

For stability data with gaps we have $A_V=id$ provided $V$ does not intersect  $V_i$'s. 
We consider  this notion for the Lie algebra of vector fields $Vect_\Gamma$.
The rank one case corresponds to  $m=2$, when we have only two sectors $V_\pm$.  
The situation with several gaps is in a sense similar to the rank one case because of the following result.

\begin{prp} Stability data with gaps are analytic if and only if each $A_{V_i}(x^\gamma), \gamma\in C(V_i), 1\le i\le m$ is a series which defines an analytic germ at $x=0$ (i.e. each element $A_{V_i}$ is analytic).
\end{prp}

{\it Proof.} If stability data are analytic then all $A_{V_i}$ are analytic. Indeed we can find a cyclic open cover by half-planes such that each half-plane contains a single sector $V_i$.

Conversely, suppose we have stability data with rational central charge and gaps and such that each $A_{V_i}(x^\gamma)$ gives an analytic germ. Using the construction with coordinate systems coming from immersed rational curves as we did before we can prove analyticity of the corresponding germ, since the only non-trivial gluing automorphisms are $A_{V_i}, 1\le i\le m$. $\blacksquare$

\subsection{Stability data with exponential bound}

Let us fix a Euclidean norm $|\bullet|$ on $\Gamma_\R=\Gamma\otimes \R$.
Suppose that the Lie algebra $Vect_\Gamma$ carries (as a graded vector space) a Banach norm, which we denote by the same notation $|\bullet|$ (i.e. each graded component carries such a norm) .
\begin{defn} Stability data on $Vect_\Gamma$ are said to satisfy  exponential bound if there exist positive constants $c_1, c_2$ such that we have 
$|a(\gamma)|\le c_1e^{c_2|\gamma|}$ for all $\gamma\in \Gamma-\{0\}.$

\end{defn}

Exponential bound condition is equivalent to the following {\bf Uniform Convergence Assumption} (UCA for short) : 
$$\sum_{l\subset V}log\left(A_l(x^{\gamma_i})\over{x^{\gamma_i}}\right)$$
converges for each strict sector $V$ and a choice of primitive generators $\gamma_i, 1\le i\le n$ of $C(V)\cap \Gamma$.

\begin{prp} Under the Uniform Convergence Assumption the following holds:
for any admissible sector $V$ the element $A_V$ is analytic.
\end{prp}

{\it Proof.} Consider a pair of tube domains $U_1\subset U_2\subset {\bf T}^{an}$ which satisfy the following property.  Let us a choose a norm $|\bullet|$ on $\Gamma$. It gives rise to the norm $|\bullet|_\R$ on $\Gamma_\R$. We have also the ``tropical map" $log|\bullet|: {\bf T}^{an}\to \Gamma_\R$. The UCA implies that the expression
$$\sum_{l\subset V}sup_{x\in U_2}|log|A_l(x)/x||_\R:=C<\infty,$$
where $log|A_l(x)/x|$ is the image of the series $A_l(x)/x=1+...$ under the tropical map.
We can choose $U_1$ and $U_2$ in such a way that the constant $C$ is smaller then the distance between $U_1$ and $\Gamma_\R-U_2$.

Then the  clockwise ordered product $\prod_{l\subset V}A_l$ gives rise to  an analytic map $U_1\to U_2$.  But this product is equal to $A_V$. $\blacksquare$

\begin{rmk}

a) The exponential bound is equivalent to the following one: there exist $c_1,c_2$ as in the Definition such that $\sum_{\gamma\in C(l), ||\gamma||\le const}|a(\gamma)|\le c_1e^{c_2\cdot const}$ for each ray $l\in \R^2$ and the corresponding strict convex cone $C(l)$ associated with $l$ ( recall that the ray $l$ is a strict sector in $\R^2$).

b) Let $l$ be an admissible ray in $\R^2$ and $A_l$ the corresponding element of the group $G_l$. Then we can consider $log(A_l)$ either as an element of the Lie algebra $\g_l$ or as a linear operator on the vector space of monomials $x^\mu, Z(\mu)\in l$, as above. The exponential bound on stability data is {\bf not equivalent} to the property that $log(A_l)(x^\mu)$ is convergent.

c) The set of stability data with exponential bound does not coincide with the set of analytic stability data. On the other hand, assume that  $\Gamma$ carries an integer skew-symmetric form and pick up a generic central charge. Then the Lie algebra of Hamiltonan vector fields corresponding to a rational ray is commutative. The  UCA is equivalent to the exponential bound on DT-invariants: $|\Omega(\gamma)|\le c^{|\gamma|}$ for some $c>0$.

\end{rmk}

\begin{que} 
Consider rational stability data on $Vect_\Gamma$. Is it true that  analyticity of the stability data is equivalent to the following property: there are exponential bounds as above such that the constants $c_1,c_2$ 
can be chosen in such a way that the corresponding inequalities hold with the same $c_1,c_2$ for all $\gamma\in C(V)$ and all strict sectors $V\subset \R^2$?
\end{que}

\begin{prp} Suppose we are given rational  stability data on $Vect_\Gamma$. Then we have exponential
bounds on $|a(\gamma)|$  for $\gamma\in C(l)$ for a strict convex cone $C(l)$ associated with a rational ray $l$. 

\end{prp}

{\it Proof.} Let us choose a decomposition $\R^2=\cup_{i\in I}V_i$ into the union of the wheel of strict rational sectors $V_i$.
Let us pick one of the sectors, say $V_{i_0}$ and
cut it into three subsets by a rational ray $l$, i.e. $V_{i_0}=V_{i_0}^+\sqcup l\sqcup V_{i_0}^-$. Then we have $A_{V_{i_0}}=A_{V_{i_0}^+}A_lA_{V_{i_0}^-}$. If $A_{V_{i_0}}$ is analytic then all three factors in the RHS are also analytic.  In particular, the map $A_l$ is analytic,
hence the exponential bound for $a(\gamma)$ with $\gamma\in C(l)$ holds. On the other hand any rational ray belongs to some $V_i$. $\blacksquare$

\begin{conj} The above Proposition holds for any, not necessarily rational, central charge.

\end{conj}

\subsection{Examples of analytic stability data}\label{examples}

There are several geometrically defined stability data (more general, wall-crossing structures) {\it which are expected} to be analytic and hence satisfy the exponential bound. 

1) Let $C$ be a complex smooth projective curve, $\Sigma\subset T^\ast C$ be a spectral curve (i.e. finite ramified cover of $C$), and
$X=X_\Sigma$ be the non-compact Calabi-Yau $3$-fold associated with these data (see e.g. [KoSo3]). It is expected that semistable objects of the  compact Fukaya category $\mathcal{F}_c(X_\Sigma)$ bijectively correspond to spectral networks on $C$ (see e.g. [GaMoNe2]).
The central charge is 
$Z(\gamma)=\int_{\gamma}\Omega^{3,0}, \gamma\in \Gamma$, where $\Omega^{3,0}$ is the holomorphic volume form, and $\Gamma=H_3(X_\Sigma, \Z)=H_1(\Sigma, \Z)$. In this way we obtain stability data on the graded Lie algebra of Hamiltonian vector fields acting on the Poisson torus ${\bf T}_\Gamma$. It is expected that these stability data are analytic. This should come from  exponential bounds on the growth of Donaldson-Thomas invariants $\Omega(\gamma)$. One can generalize the story to the case of marked points.

Notice that although in the case of compact Calabi-Yau $3$-fold $X$ we can define the stability data in a similar way,
the corresponding stability data are not expected to be exponential or analytic. In particular, Ooguri-Strominger-Vafa conjecture predicts that the Donaldson-Thomas invariants $\Omega(\gamma)$ grow as $e^{c|\gamma|^2}$.

Instead of the cotangent bundle $T^\ast C$ one can consider its ``trigonometric" analog $(\C^\ast)^2$, and the corresponding spectral curve, Calabi-Yau $3$-folds, etc. The  estimates for $\Omega(\gamma)$ are not known. Nevertheless one can hope for analytic stability data in this case as well.

2) The more general source  of analytic stability data and analytic wall-crossing structures (the latter will be briefly discussed later) is the theory of complex integrable systems (see [KoSo3]). There, instead of a single lattice $\Gamma$ one has a local system with the generic fiber isomorphic to the lattice $\Gamma$, which is the first integer cohomology of the corresponding complex torus (and hence it is symplectic if the torus is compact). Calabi-Yau threefolds (more generally and conjecturally, $3CY$ categories) give rise to complex integrable systems with fibers which are intermediate Jacobians of those $3$-folds (resp. categories). Let us  consider the case of local Calabi-Yau threefolds. In that case fibers are {\it polarized} abelian varieties. The base $B$ of a complex integrable system (or rather its open dense part $B^{sm}$ parametrizing non-degenerate tori) can be thought of as a middle-dimensional  complex  submanifold in the space of stability data on the  graded Lie algebra $\Gamma$. In case if the complex symplectic form is exact, one has a family of
homomorphisms $Z_b:\Gamma_b\to \C, b\in B^{sm}$. In order to define the  stability data on a given $\g_{\Gamma_b}$  it is necessary to have integers $\Omega_b(\gamma), \gamma\in \Gamma_b$. They are also called Donaldson-Thomas invariants  of the corresponding data ($3CY$ category, complex integrable system, etc.). We explained in the loc.cit. how to construct (under certain assumptions) DT-invariants for complex integrable systems. By construction they satisfy the WCF from [KoSo1].  For more details see [KoSo3]. Conjecturally arising wall-crossing structures  satisfy the exponential bounds which are uniform. On the contrary, non-polarized complex integrable systems arising from compact Calabi-Yau $3$-folds should give wall-crossing structures which do not satisfy the  exponential bounds.

3) (Higgs bundles)
This example is in a sense a special case of the previous one.
Let $M=T^\ast C$ where $C$ is a smooth complex projective curve of genus $g$. Let $\Sigma\subset M$ be a smooth spectral curve. Then the constructions of the above two subsections have the following interpretation. Let $\pi: M_{Higgs}\to B$ be the Hitchin map, and $A_\Sigma=\pi^{-1}(\Sigma)$ be the abelian variety over $\Sigma\in B$. We have a holomorphic family $Conn_t$ of $t$-connections, which coincides with $M_{Higgs}$ as $t=0$. Let $\MM^{tot}$ be the total space of this family. Let us consider the blow-up $Bl_{A_\Sigma}(\MM^{tot}):=\MM^{tot}_{A_\Sigma}$. The exceptional divisor is a ${\bf P}^g$-bundle over $A_\Sigma$. It contains as an open subset a $\C^g$-bundle over $A_{\Sigma}$. The latter is analytically isomorphic to $(\C^\ast)^{2g}$ (rank one local systems on $\Sigma$). Let us consider the analytic neighborhood $\MM^\prime$ of this $(\C^\ast)^{2g}$ in $\MM^{tot}_{A_\Sigma}$. There is a natural map $\MM^\prime\to \C_t$. Over formal series $\C[[t]]$ the corresponding formal neighborhood can be trivialized by a result of D. Arinkin (see [Ar]). 
Conjecturally the analytic germ $\MM^\prime$ together with the above formal trivialization is isomorphic to ${\mathcal X}$ defined in the  next section and is an analytic germ corresponding to some analytic stability data.

Furthermore, one can consider the ``Betti side", where $M_{Higgs}$ is replaced by the corresponding $M_{Betti}$. There is a compactification of $M_{Betti}$ by log-divisors, and the analytic neighborhood of the boundary divisor is isomorphic to the analytic neighborhood of a wheel of projective lines, as explained previously. The corresponding analytic stability data should be compared with those in Example 1.

4) ($2d-4d$ wall-crossing formulas)

Assume that the lattice $\Gamma$ carries an integer skew-symmetric form. Then we have the corresponding Poisson torus ${\bf T}_\Gamma$. Let us fix also $m\ge 1$. Consider the cross-product $\g_{\Gamma,m}=\OO({\bf T}_{\Gamma})\ltimes Mat(m, \OO({\bf T}_{\Gamma}))$.

The fiber product structure is a corollary  of a general construction of $\g_1\ltimes \g_2$ for Lie algebras $\g_i, i=1,2$ provided there is a homomorphism of Lie algebras $\g_1\to Der(\g_2)$. 
In our case we have a decomposition $\OO({\bf T}_\Gamma)\simeq Der(\OO({\bf T}_\Gamma))\oplus Center$, where $Center$ means the Poisson center.  Then we have the action of the Lie algebra of derivatives on the matrix elements. The cross-product is then analogous to the Atiyah algebra of a vector bundle.
Notice that there is an action of the torus $(\C^\ast)^n\times (\C^\ast)^m$ on $\g_{\Gamma,m}$ , which gives a grading on the latter.

In the end we have an embedding $\g_{\Gamma,m}\to Vect((\C^\ast)^{n+m})$.  Then the previous considerations which can be thought of as the case $m=0$ can be generalized. The resulting wall-crossing structure underlies the so-called $2d-4d$ wall-crossing formulas.  Our definition of analytic stability data easily extends to this case.

5) Stability data arising from complexified Chern-Simons functional are expected to be analytic. On the contrary,
in the case of holomorphic Chern-Simons functional  for a compact Calabi-Yau $3$-fold the stability data are determined by the virtual number $\Omega(\gamma)$ of special Lagrangian submanifolds. Thus it is essentially the same story as in the case of compact Calabi-Yau $3$-folds mentioned in Example 1. Hence the corresponding stability data are not expected to be analytic.

\subsection{Non-archimedean stability data}

Let now $K$ be a non-archimedean field, i.e a field endowed with non-trivial complete non-archimedean Banach norm. Then all the above results remain true. Proofs are basically the same, as long as we assume a non-archimedean version of the Douady theorem. Construction of the $(n-1)$-dimensional family of smooth projective curves is in a sense easier, since we can lift the family of straight lines in $\R^2-\{0\}$ to the corresponding non-archimedean tube domains. The difference of the arising notion of non-archimedean analytic stability data with the one of complex analytic stability data is that the former notion is equivalent to the notion of stability data with exponential bound. We will discuss the details elsewhere.

\subsection{Analytic wall-crossing structures}
In this subsection in order to save space, we will heavily use the notation of [KoSo3], Section 2. 

According to [KoSo3], Section 2.2 wall-crossing structures can be defined as sections of certain sheaf of sets.
Hence it suffices to consider wall-crossing structures on an open subset $U\subset \Gamma_\R^\ast$. There a WCS is determined by a collection of elements $g_{y_1,y_2}$ belonging to a certain pronilpotent group $G_C$ associated with our graded Lie algebra $\g$ and a strict closed convex cone $C\subset \Gamma_\R$. The points $y_1,y_2$ are assumed to be {\it $C$-irrational} in the sense that they do not belong to the union of walls $\gamma^\perp, \gamma\in C\cap (\Gamma-\{0\})$.

The elements $g_{y_1,y_2}$ satisfy some conditions. They include the cocycle condition $g_{y_1,y_2}g_{y_2,y_3}=g_{y_1,y_3}$ as well as the condition that the adic limit of $g_{y_1,y_2}$ as long as $y_1,y_2$ approach the same generic point of the wall $\gamma^\perp$ belongs to the group $G_{\R_{>0}\cdot \gamma}$ (see [KoSo3], Section 2 for details).

We extend the definition of $g_{y_1,y_2}$ to the case when $y_i, i=1,2$ are arbitrary points of $U-\{0\}$, i.e. they can belong to walls.
More precisely, for each $y\ne 0$ we define its {\it decoration} $\epsilon$ as an element of the set $\{+, -\}$. Informally decoration is responsible for the ``direction of the infinitesimal shift" of $y$.

For $y\ne 0$ let  $\mathcal{A}_y$ be the set defined such as follows:

a)  $\mathcal{A}_y=\emptyset$ if $y$ is $C$-irrational;

b) $\mathcal{A}_y=\{\gamma\in C\cap(\Gamma-\{0\})|y(\gamma)=0\}$ otherwise.

Then we define two open convex cones $C_y^\pm\subset \Gamma_\R^\ast$ such that $C_y^\pm=\{y^\prime\in \Gamma_\R^\ast-\{0\}| \pm y^\prime(\gamma)>0, \gamma\in \mathcal{A}_y\}$. Notice that $y$ belongs to the closures $\overline{C}_y^\pm$ of both opposite cones.

For a pair of  points $y_i\ne 0, i=1,2$ and a choice of their decorations $\epsilon_i, i=1,2$ we 
will  define an element $g_{y_1,y_2}^{\epsilon_1, \epsilon_2}\in G_C$ as the limit in the adic topology of the elements $g_{y_1^\prime,y_2^\prime}$ such that $y_i^\prime\in C_y^{\epsilon_i}, i=1,2$ are $C$-irrational points which approach to the corresponding $y_i, i=1,2$.
The new collection of elements still satisfy the cocyle condition. Notice that if a point $y_i, i=1,2$ is $C$-irrational then the element $g_{y_1,y_2}^{\epsilon_1, \epsilon_2}$ does not depend on $\epsilon_i, i=1,2$.

\begin{defn} Let $\g=Vect_\Gamma$. In the above notation we call the wall-crossing structure analytic, if for all sufficiently close pairs of points $y_i\in \Gamma^\vee\otimes \Q-\{0\}$ endowed with arbitrary decorations $\epsilon_1, \epsilon_2$ the elements  $g_{y_1,y_2}^{\epsilon_1, \epsilon_2}$ are analytic.

\end{defn}
In other words the property of the wall-crossing structure to be analytic can be checked on sufficiently close pairs of rational points in $\Gamma_\R-\{0\}$.

Suppose that we are given analytic stability data in the sense of Section 3. Assume that $Z\ne 0$. Then we have a linear map $Z^\ast:\R^2\to \Gamma_\R^\ast$. Recall (see [KoSo3]) that in this case we can define a WCS as a section of a sheaf of sets in a neighborhood of $Z^\ast(\R^2-\{0\})$. Then this WCS is analytic iff the stability data are analytic. In order to see that one observes that the above-described group elements $g_{y_1,y_2}^{\epsilon_1, \epsilon_2}$ can be expressed via the transformations $A_{V_i}$ and $A_{l_i}$ from the Proposition 3.4.2 for a certain choice of rational rays $l_i$ and rational central charges sufficiently close to $Z$.

\section{Analytic stability data and resurgence}\label{analytic stab data and resurgence}

Let $c(t)=\sum_{n\ge 0}c_nt^n$ be a formal power series. Consider its Borel transform (equivalently, Laplace transform along the ray) 
${\mathcal B}(c)=
\sum_{n\ge 0}{c_n\over{n!}}\lambda^n$. The series $c$ is called resurgent if ${\mathcal B}(c)$ defines a germ of analytic 
function at $\lambda=0$ which admits endless analytic continuation. The latter means that there is a countable subset $S\subset \C$
which is an inductive limit of finite subsets $S_L$ such that ${\mathcal B}(c)$ can be analytically continued along any path in $\C-S_L$
of the length less or equal than $L$. This notion goes back to \'Ecalle (see [Ec]) and provides a useful tool for working with divergent series.
We will discuss below the relationship of resurgence series with the notion of analytic stability data. Our main conjecture claims that the latter gives rise to the former.

\subsection{Introducing  parameter to analytic stability data}\label{introducing parameter}

Given $\C$-analytic stability data on $Vect_\Gamma$ we will construct a germ of $(n+1)$-dimensional complex manifold ${\mathcal X}$ such that:

a) $\XX\supset \XX_0:=(\C^{\ast})^n$ and $\XX$ is an analytic germ at $\XX_0$.

b) Holomorphic map $\pi: \XX\to \C$ which is submersion s.t. $\pi_{|\XX_0}=0$. 

c) The formal completion $\widehat{\XX}_{\XX_0}$ is indentified with $\XX_0\times Spf(\C[[t]])$ in such  a way that $\pi$-pull-back of
the standard coordinate on $\C$ is equal to $t\in \OO(\XX_0\times Spf(\C[[t]]))$.

But first let us recall the description of the set of stability data of rank $2$ on $Vect_\Gamma$ given in  Section \ref{wheels and schemes}.
According to loc. cit. as long as we fixed the central charge, this set
can be identified with the quotient set $\MM_{(C_i)_{i\in J}}$, where $(C_i)_{i\in J}$ is a compatible wheel of cones.

Suppose we are given analytic stability data $\sigma$ with the central charge $Z^\sigma$.
Let us  discuss  the construction of  $\XX$. Consider the $C^\infty$ manifold $N_0:=S^1_\theta\times Hom(\Gamma, \C^{\ast})$ and fix its diffeomorphism with $S^1\times (\C^{\ast})^n$.  Let $p_1$ be the projection to the first factor.  
Then we can define a constructible sheaf of groups $\underline{\mathcal G}$ on the circle $S^1_\theta$ (the subscript denotes the coordinate). We may assume that we are given a wheel  of cones $(C_\nu)_{\nu\in \Omega}$ in $\Gamma_\R^\ast$  such that $V_\nu^\sigma:=Z^\sigma(C_\nu^\vee)$ are strict admissible and they form a wheel.  Let us choose admissible rays $l_s, s\in \Omega$ in such a way that for any sector $V_s$ bounded by $l_s$ and $l_{s+1}$,
an element $\gamma\in C_s^\vee$ and $t\in V_s$ we have $Re(Z^\sigma(\gamma)/t)>0$. Denote the intersection  of $l_s$ with the standard circle $S^1=\{|t|=1\}$ by $p_s$. Then we have a stratification of $S^1$ by points $p_s, s\in \Omega$ and the complementary arcs. Using the groups $G_s^{old}$ and $G_{s,s+1}^{old}$ corresponding to $V_s$ and $V_s\cap V_{s+1}$ respectively we obtain the constructible sheaf of groups $\underline{\mathcal G}^{old}$ on $S^1$ such that the stack at $p_s$ is $G_{s,s+1}^{old}$ while the stack at the adjacent to $p_s$ ( in the clockwise order) arc is $G_s$. We have the natural maps $G_{s}^{old}\leftarrow G_{s,s+1}^{old}\rightarrow G_{s+1}^{old}$ (specialization maps for the constructible sheaf). Then the quotient set describing the set of stability data can be interpreted as an element of $H^1(S^1,\underline{\mathcal G}^{old})$ (which is trivial). The corresponding analytic version can be interpreted as $H^1(S^1,\underline{\mathcal G}^{an})$, where we use subgroups of analytic elements  $G_s^{an}$ and $G_{s,s+1}^{an}$ in the definition of the sheaf of groups $\underline{\mathcal G}^{an}$.

Then we have the sheaf of groups $p_1^\ast(\underline{\mathcal G}^{an})$ on $N_0$.  There is another sheaf of groups $\underline{\mathcal H}$
on $N_0$ consisting  of diffeomorphisms of the $C^\infty$ germ at $N_0$ of the manifold $N=S^1_\theta\times [0,\varepsilon]_r\times (\C^\ast)^n$ (here $\varepsilon>0$ is small) which are equal to $id_{N_0}$ on $N_0$ with all derivatives.  We have a homomorphism of sheaves of groups $p_1^\ast(\underline{\mathcal G}^{an})\to \underline{\mathcal H}$.
We should think of the pair $(r,\theta)$ as of polar coordinates corresponding to the point $t$. We should think of the circle $S^1_\theta$ as the oriented real blow-up (circle of directions) of $\C_t$ at $t=0$. Consider diffeomorphism of the torus over the point $(r,\theta)$ defined by $\psi_{(r,\theta)}(z^\beta)=e^{-Z(\beta)/re^{i\theta}}z^\beta$ and similarly for $\overline{z}^\beta$, which is the complex conjugate monomial on the torus $(\C^{\ast})^n$ considered as a $C^\infty$ manifold. Then we conjugate each $A_{V_i}^\sigma$ with  $\psi_{(r,\theta)}$,
considering each $A_{V_i}^\sigma$ as a map of the space of monomials.  Then we get a collection of diffeomorphisms $A_{V_i, (r,\theta)}^\sigma$ whose fixed set of points is $N_0$.  Using these diffeomorphisms we glue a new $C^\infty$ manifold $N^{new}$ with the boundary $N_0$. Here are some properties of $N^{new}$:

a) It is a germ of $C^\infty$ manifold at $N_0$, which is trivialized in the formal neighborhood of $N_0$. This follows from the exponential decay of the function $e^{-a/r}, a>0$ as
$r\to 0$. Thus all diffeomorphisms are not only equal to $id_{N_0}$ on $N_0$ but the equality holds with all higher derivatives.

b) There is a complex structure on $N^{new}-N_0$ which continuously extends to $N_0=\partial N^{new}$.

c)  There is a germ of continuous map of the manifolds with boundary $N^{new}\to \C_t$ which is holomorphic outside of the boundary.
In coordinates it is given by the projection $(r,\theta, z)\mapsto t=re^{i\theta}$.

\begin{prp} Consider the topological space obtained from $N^{new}$ by contractions of all  $S^1$-orbits on $N_0$ with respect to the natural action of $S^1$. 
Then the resulting topological space carries a structure of $C^\infty$ manifold.

\end{prp}
{\it Proof.} Consider the following subsheaf in the sheaf of functions on $N^{new}$. Its sections are smooth functions on $N^{new}-N_0$ which admit continuous $S^1$-invariant  extension to $N_0$, and which admit asymptotic expansion as $r\to 0$ (i.e. at $N_0$) of the form $ \sum_{m\ge 0}r^m\sum_{|l|\le m, l-m\in 2\Z}c_{l,m}e^{i\theta l}$. One can check that after taking the quotient by $S^1$-action this sheaf gives a sheaf of $C^\infty$ functions on the quotient topological space. $\blacksquare$

Furthermore, one can see that the complex structure from b) descends to the quotient by the $S^1$-action. This gives us a germ of complex manifold together with an analytic projection to the germ of complex line $\C_t$ at $t=0$.

\subsection{Conjecture about the resurgence}\label{resurgence conjecture}

In this subsection we are going to spell out the construction of the previous subsection in a different way. This will allow us to formulate a general
conjecture which relates analytic stability data and resurgent series in the parameter $t$ mentioned above.

Let $t$ denote the standard coordinate on $\R^2=\C$  and $\C-\{0\}=\cup_{i\in J}V_i$ be a cyclic cover by open half-planes. Choose a basis of monomials $z^{\gamma}\in \OO({\bf T})$, where ${\bf T}=Hom(\Gamma, \C^\ast)$.
Let $l_i$ denote the admissible open ray perpendicular to the boundary line of $V_i$. Then for any $\beta\in \Gamma$ such that
$Z(\beta)\in l_i$ we have
 $Re({Z({\beta})\over {t}})>0$ 
as long as  $t$ belongs to the interior of  $V_i$. Also, we have a $1$-parameter subgroup $t\mapsto e^{{-Z(\beta)\over{t}}}$ of the torus ${\bf T}\simeq (\C^{\ast})^n$ associated with each $\beta\in \Gamma$.

Suppose we are given analytic stability data $\sigma$ on $Vect_\Gamma$. For each $V_i$ we can choose a representative $g_{V_i}\in G_{V_i}^{an}$.  By definition such a representative is given by a collection of analytic elements $g_{W_{i,\epsilon}}$ for the inductive limit as $\epsilon\to 0$ of closed strict admissible subsectors $W_{i,\epsilon}\subset V_i, i\in J$ such that the left (resp. right) boundary rays of 
each $W_{i,\epsilon}$ forms the angle $\epsilon$ with the left (resp. right) boundary ray of $V_i$.
Then 
$g_{W_{i,\epsilon}}(z^\gamma)=z^\gamma(1+\sum_{\beta\in Z^{-1}(W_{i,\epsilon})\cap \Gamma} c_{i\gamma\beta}z^\beta)$ is a convergent series, because
of the exponential bounds on  $c_{i\gamma\beta}$. 

Consider now the trivial analytic fiber bundle over $\C-\{0\}$ with the fiber ${\bf T}^{an}$.
We modify this trivial bundle  along all the rays $l_i$ via the analytic automorphism $g_{i,t}$ of ${\bf T}^{an}$ given at the level of monomials
by the formula
$$g_{i,t}(z^\gamma)=z^{\gamma}(1+\sum_{\beta\in Z^{-1}(l_i)\cap \Gamma}c_{i\gamma\beta}e^{-{Z({\beta})\over {t}}}z^\beta),$$
where $t\in V_i$.
The series in the RHS converges by the assumption of analyticity of $\sigma$.

In this way we obtain  an analytic fiber bundle ${\mathcal F}(\sigma)$ over $\C-\{0\}$.  Since  the Taylor series expansion at $t=0$ of all $g_{i,t}(z^\gamma)$ do  not depend on  $V_i, i\in J$ (in fact it is equal to $z^\gamma$), the bundle can be extended analytically to $t=0$. Moreover, it has a canonical trivialization over the formal neighborhood of the fiber at $t=0$ (which is ${\bf T}^{an}$).

\begin{conj}\label{analyticity and resurgence}Fix $\gamma\in \Gamma$ and denote by $\chi_\gamma: {\bf T}\to \C^\ast$ the corresponding character.
Let $s(t)$ be a germ at $t=0$ of an analytic section of ${\mathcal F}(\sigma)$. Then in the above-mentioned canonical formal trivialization the Taylor series of $\chi_\gamma(s(t))$ is resurgent.
Same is true for  the $\chi_\gamma(log(s(t))$.
\end{conj}

Clearly there are plenty of germs of analytic sections of ${\mathcal F}(\sigma)$, so the Conjecture \ref{analyticity and resurgence} will be a source of resurgent series arising from wall-crossing structures.

In practice such fiber bundles as well as their sections appear as a result of gluing of trivial fiber bundles over sectors together with solving
a collection of  Riemann-Hilbert problems for the common boundary rays. An archetypical example is the one of exponential integrals discussed later in the paper. In the simplest examples singularities of the Borel transform belong to the set $Z(\Gamma)$, but we do not have enough evidence
to make a more general prediction.

The above conjecture covers a big class of examples of  resurgent behavior of formal series (one can optimistically say  ``most of them"). Some examples will be discussed later in the paper. Physics and mathematics  literature devoted to perturbative expansions and generating functions of enumerative invariants is the rich source of such examples. There are other attempts to approach axiomatically their analytic properties. The interested reader can try to translate them into the language of our formalism. An example of such an approach is given in [Br], where exponential bounds on the growth of DT-invariants are taken as a part of the axiomatics.

\subsection{Chains of projective lines and twistor family}\label{bundle and chains}

Recall that for given analytic stability data $\sigma$ of rank $2$ we can construct the open analytic space $X^{an}$ which is the analytic neighborhood of the wheel of projective lines corresponding to a wheel of closed cones in $\Gamma_\R^\ast$ (see Section \ref{wheels and schemes} and Section \ref{analytic germs}). The corresponding formal scheme  $X^{form}$ is obtained by gluing trivial formal neighborhoods of projective lines in the toric variety associated with the wheel of cones $(C_i)_{i\in J}$  via automorphisms $g_{C_i^\vee}\in G_{C_i^\vee}$ corresponding to the interiors of the cones. Analyticity of the stability data ensures that the above procedure can be upgraded to the level of analytic spaces.
The corresponding analytic space $X^{an}$ is an open analytic space of dimension $n=rk\, \Gamma$, with the underlying ``tropical variety" given by the wheel of cones  $C_i\subset \Gamma_\R^\ast, i\in J$. As a topological space $X^{an}$ has the homotopy type of the torus $(S^1)^{n+1}$.
Then we have a trivial analytic fiber bundle over $\C^\ast_t$ with the fiber $X^{an}$. It carries the  trivial  non-linear flat connection. Furthermore the bundle can be analytically extended to the point $t=0$, with the fiber over $t=0$ being $(\C^\ast)^n$.

On the other hand, the construction of Section \ref{resurgence conjecture} (and even an earlier Section \ref{introducing parameter}) gives rise to an analytic fiber bundle over $\C_t^\ast$ with a non-linear connection, coming from the maps $z^\gamma\mapsto g_{i,t}(z^\gamma)$. Then we have an open embedding of this bundle with non-linear connection to the trivial one.

\begin{rmk}\label{analogy with Higgs bundles}

The reader can compare the above considerations with  the construction of the analytic bundle with non-linear flat connection over $\C_t$ with the fiber at $t=0$ given by  moduli space $\MM_{Higgs}$ of possibly irregular Higgs bundles on a curve (see [KoSo3]). The non-linear analytic fiber bundle over $\C^\ast_t$ has  the fiber given by the cluster variety $\MM_{Betti}$. The chain of projective lines ``at infinity" for each fiber  can be visualized by means of  the WKB expansion as $t\to 0$ of the  flat sections of the corresponding (by the Riemann-Hilbert equivalence) $t$-connection. The cluster variety is ``glued" from tori ${\bf T}$ corresponding to the lattice $\Gamma$, which is essentially the first homology group of the spectral curve. In this example we have a wall-crossing structure rather than a stability data on $Vect_\Gamma$, as at the beginning of this subsection. Fiber at $t=0$ is obtained by the blow-up of $\MM^{tot}$ (see Section \ref{examples}, example 3)) at an abelian variety, which is a fiber of the Hitchin fibration $\MM_{Higgs}$ placed at $t=0$. All that will be discussed in more detail in our subsequent papers on our big project ``Holomorphic Floer Theory" (see [KoSo4]).
\end{rmk}

The above considerations  can be upgraded further giving a  full twistor family of analytic manifolds over ${\bf P}^1$ as long as we impose some ``reality conditions" on the stability data. The twistor family (equivalently hyperk\"ahler metric, maybe incomplete on the fiber at $t=0$) should be obtained via the techniques analogous to [GaMoNe1].

\subsection{\'Ecalle-Voronin theory and analytic stability data}\label{Ecalle and analytic stability}

In this subsection  we discuss a non-trivial example in which the Conjecture \ref{analyticity and resurgence} can be verified.
More precisely, we are talking about resurgence of series arising in the \'Ecalle-Voronin theory (see e.g. Section 7 of [Mit Sau]).
Let us describe the corresponding analytic stability data.

Let $\Gamma=\Z$, and the central charge $Z: \Gamma\to \C$ is given by $Z(n)=2\pi in$.
Then the Lie algebra of vector fields on the torus ${\bf T}_\Gamma=\C^{\ast}$ with the chosen coordinate $z$ is $Vect_\Gamma=\g=\oplus_{n\in \Z}\g_n$,
where $\g_n=\C\cdot z^n(z\partial_z)$. 

For strict admissible sectors $V_\pm$ containing $\pm i\R_{>0}$
we define $A_{V_\pm}^{\pm 1}: z\mapsto z\cdot exp(\sum_{n\ge 1}a_{\pm n}z^{\pm n})$, where the series
$f_\pm=\sum_{n\ge 1}a_{\pm n}z^{\pm n}$ converges as $|z|^{\pm 1}$ belongs to a disc of sufficiently small radius.
Then $A_{V_\pm}$ give rise to the analytic stability data.

For any $t$ which belongs to a  small sector containing $i\R_{>0}$ there is an automorphism $\psi_t^+: \C^{\ast}\to \C^{\ast}$ given by $\psi_t(z)=e^{-2\pi i/t}z$. We can twist $f_+$ using $\psi_t$ such as follows:
$$\psi_t(f_+): z\mapsto z\cdot exp(\sum_{n\ge 1}a_ne^{-2\pi in/t}z^n).$$

For any $t$ which belongs to a small sector containing the ray $-i\R_{>0}$ we define $\psi_t^-(z)=e^{2\pi i/t}z$. It gives the twist of $f_-$ such as follows:
$$\psi_t(f_-): z\mapsto z\cdot exp(\sum_{n\ge 1}a_{-n}e^{2\pi in/t}z^{-n}).$$

Next we obtain via the gluing procedure described below a germ of complex surface $Z$ at $\C^\ast$ together with the natural map $p: Z\to \C_{t=0}$, where $\C_{t=0}\subset \C_t$ denotes the analytic germ of $\C_t$ at $t=0$. 

The surface $Z$ is constructed such as follows (cf. Sections \ref{introducing parameter}, \ref{resurgence conjecture}). Let us consider two subsets of the torus $\C^\ast$ endowed with the coordinate $z$ given by
$\Delta_\epsilon=\{|z|<\epsilon\}$ and $\Delta_\epsilon^{-1}=\{|z|> \epsilon^{-1}\}$ such that $f_+$ converges in the first disc, while $f_-$ converges in the second one.
In the plane $\C_t=\R^2$ consider two open domains $U_L$ and $U_R$ defined such as follows. The domain $U_L$ is obtained by removing
from a small open disc $|t|<\delta$ the closed sector $\pi-\theta\le Arg(t)\le \pi+\theta$, where $\theta$ is sufficiently close to $\pi/2$.
Similarly we define $U_R$ as the domain obtained by removing from the same small disc the closed sector $-\theta\le Arg(t)\le \theta$ with the same
condition on $\theta$. Then $U_L\cap U_R$ is a disjoint union of the opposite sectors $S_+\cup S_-$. Consider now the set of holomorphic maps $z_L:U_L\to \C^\ast_z$ and $z_R:U_R\to \C^\ast_z$ such that $z_L(0)=z_R(0)=1$. We can choose $\delta>0$ to be such small that $1-\nu<|z_L|<1+\nu$ for sufficiently small $\nu$.
Namely, we can satisfy the condition $|e^{-2\pi i/t}z_L|<\epsilon$ if $t\in S_+$ and $|(e^{-2\pi i/t}z_L)^{-1}|<\epsilon$ if $t\in S_-$. It follows that the expressions
$f_{\pm}(e^{-2\pi i/t}z_L)$ are well-defined.

Now we treat the above pair $(z_L, z_R)$ as a holomorphic section of the projection $p: Z\to \C_t$  over $\{|t|<\delta\}\subset \C^\ast_t$ 
provided on $S_\pm$ they satisfy the relations
$$z_R=z_L\cdot exp(f_\pm(e^{2\pi/t} z_L)).$$
This defines the germ $Z$ which satisfies the desired properties. 
Clearly $Z$ carries a holomorphic foliation defined by the holomorphic vector field $v=t^2\partial_t$. The image of a small disc (or a germ at $t=0$) under the above section is then endowed with the Poincar\'e map $P$ induced by the return map for the foliation.

Notice as $t\to 0$ in $S_+\cup S_-$ (i.e. where both $z_L$ and $z_R$ are defined) then the asymptotic expansions of both series coincide.
Thus $log(z_L)$ and $log(z_R)$ has Taylor expansion $\sum_{n\ge 0}c_nt^n$, where $|c_n|\le C^nn!$ for some $C>0$.

The data $\Phi=(f_+,f_-, z_L, z_R)$ are equivalent to those in the \'Ecalle-Voronin theory (see e.g. [MitSau]).
We will formulate the relationship between two types of data such as follows.

\begin{thm}\label{equivalence} Let us assign to the data $\Phi$ a germ $g: \C_{t=0}\to \C_{t=0}$ given by the Poincar\'e map $P$, provided we take points in the image of the section $p: Z\to \C_{t=0}$ sufficiently close to the circle $|z|=1$. Then

a) $g(t)=t+t^2+...$, where dots denote higher terms;

b) the above assignment is a bijection with the set of germs of the form a);

c) keeping $f_\pm$ fixed but changing $(z_L, z_R)$ is equivalent to the conjugation $g\mapsto hgh^{-1}$, where $h$ is a germ of the automorphism of $\C_{t=0}$ of the form $h(t)=t+...$, where dots denote higher order terms.

d) consider Taylor expansions of $log(z_L)$ and $log(z_R)$ at $t=0$ which in fact coincide and have the form $t+\sum_{n\ge 2}c_nt^n$. 
Then this  series transforms $g$ to the normal form over formal series. More precisely,
let $w=w(t)=t\cdot exp(\sum_{n\ge 2}(c_nt^n))$. Then $g(w(t))={w(t)\over{1-w(t)}}$.

\end{thm}

Resurgence of the series $\sum_nc_nt^n$ follows from  \'Ecalle-Voronin theory (see loc.cit.)

Let us comment on the Theorem \ref{equivalence}.

First we make the  substitution $x=1/t$
and consider the corresponding disc $\Delta_\delta^{-1}=\{\delta^{-1}<|x|<\infty\}$ about $x=\infty$. Let us cover $\Delta_\delta^{-1}$
 by the union of two open domains
$W_L$ and $W_R$ which ``looks'' as complements to the interior of parabolas on the plane. This means that $W_L=\{(u,v)|u>v^2+r\}$ and
$W_R=\{(u,v)|u<v^2-r\}$ for  an appropriately chosen $r>0$. Here $x=u+iv$. Notice that $W_L\cap W_R=J_+\sqcup J_-$ similarly to the story with $U_L$ and $U_R$ discussed above.

 It is known that the classification of the changes of coordinates transforming
our automorphism to the standard form $t\mapsto t+t^2$ is controlled  by a pair of holomorphic  germs (see [MitSau]). In the variable $x$ the normal form is $x\mapsto x+1$. We can always transform the germ of an analytic map at $x=\infty$ to this normal form {\it over formal series in $1/x$}. 
Analytic transformation depends on a pair of convergent series, which we will identify below with $f_\pm(1/x)$. This  can be seen such as follows.
Let $y_L(x)=x+{1\over{2\pi i}}log(z_L(1/x)), y_R(x)=x+{1\over{2\pi i}}log(z_R(1/x))$, be the analytic functions defined in $W_L$ and $W_R$ respectively (for them
to be well-defined we choose an appropriate $r$. Then the above formulas connecting $z_L$ and $z_R$ on $S_+\sqcup S_-$ lead to the following connection formulas for $y_R$ and $y_L$:
$$(y_R-y_L)(x)={1\over{2\pi i}}(f_\pm(e^{2\pi i y_L(x)}), x\in J_\pm.$$

Consider now the transformation $T: y_L\mapsto y_L+1, y_R\mapsto y_R+1$. Notice that it is well-defined on $J_+\sqcup J_-$, since if $y_L\mapsto y_L+1$ then $y_R\mapsto y_R+1$ by the above connection formulas. Returning to the variable $x$ we see that $T(x)=x+1+\sum_{n\ge 1}{b_n\over{x^n}}$, where the series converges in $\Delta_\delta^{-1}$.  As we have seen above, transforming $T(x)$ to its normal form depends on two convergent series $f_\pm$.  Furthermore, the common  formal expansion of $z_L$ and $z_R$ as $t\to 0$ should be a resurgent series by the Conjecture \ref{analyticity and resurgence}.  This agrees with
the  \'Ecalle-Voronin theory.

\begin{rmk} The geometry underlying above considerations is the one of a complex surface endowed with a complex one-dimension foliation.
It is given by the vector field $v=t^2\partial_t+2\pi i\partial_z$ in the coordinates $(t,z)$ on the total space of the trivial $\C^\ast_z$ bundles over
the domains $U_L$ and $U_R$. This vector field is preserved under the gluing formulas, hence it gives a complex foliation on the total
space of the corresponding surface. The automorphism $T$ arises after a choice of a curve transversal to the foliation. Then it becomes the Poincar\'e return map.

This remark gives a hint to a higher-dimensional generalization of the above story. The fiber of the trivial bundle over the domains in $\C^\ast_t$
will be ${\bf T}_{\Gamma}=Hom(\Gamma, \C^\ast)\simeq (\C^{\ast})^n$. The vector field defining the foliation can be chosen as
$v=t^2\partial_t+\sum_{1\le i\le n}z_i\partial_{z_i}$, where $rk\,\Gamma=n$ and $z_i=Z(e_i)$ are the coordinate functions on ${\bf T}_{\Gamma}$.
Here $Z$ is the central charge and $e_1,...,e_n$ is a chosen basis of $\Gamma$. The automorphism generalizing $T$ will be the return map
after a choice of a hypersurface transversal to the foliation.

\end{rmk}

\section{Variation of the central charge and non-commutative Hodge theory}\label{variation of central charge}

In this section we generalize the previous discussion by considering the central charge $Z\in Hom(\Gamma,\C)$ as a variable.  This makes the story different from the one for $\MM_{Higgs}$. \footnote{Probably the latter can be generalized in such a way that the analogy persists, e.g. by allowing to vary the complex structure of the underlying curve.}
We will present the new  structure in axiomatic way.  The results will be formulated in terms of germs of complex manifolds. 
Recall that  a germ of a complex analytic space  is defined as a ringed locally compact topological space $(X_0,\OO_{X_0})$ which is in a neighborhood of each point is isomorphic as a ringed space to the pair $(K, (\OO^{an}_Y)_{|K})$, where $Y$ is a complex analytic space, $K\subset Y$ is a locally closed complex analytic subset and $(\OO^{an}_Y)_{|K}$ denotes the restriction of the sheaf of complex analytic functions.  

We will interpret  the trivial stability data on $Vect_\Gamma$ in terms of a non-linear flat connection of certain type on the trivial analytic fiber bundle over $U\times Spf(\C[[t]])$, where $U\subset Hom(\Gamma,\C)$ is an open domain. 
We call such connections {\it standard}. Then we will introduce a more general class of analytic fiber bundles with non-linear flat connections which we call  {\it almost standard}. Then we will prove that every almost standard connection is gauge equivalent over $\C[[t]]$ to a standard one. This can be thought of as an analog of the Hukuhara-Levelt-Turrittin theorem (HLT theorem for short) which gives a formal classification of vector bundles with meromorphic connections on the formal disc. 
Analytic stability data appear in the story as Stokes automorphisms associated with  analytic fiber bundles endowed with non-linear connections
which are almost standard in the formal with respect to $t$ sense. 

\subsection{Trivial stability data and standard connections}\label{trivial WCS}
We start with an example of the trivial stability data, in which case all data and axioms below is easy to verify. 

Let $\Gamma=\Z^n, X=Hom(\Gamma, \C)=\C^n$. We endow $X$ with natural coordinates $Z_1,...,Z_n$, which we can think of coordinates of the central charge $Z: \Gamma\to \C$. Let $\pi: E\to X\times \C$ denote the trivial fiber bundle with the fiber $(\C^\ast)^n$. We endow the factor $\C$ with the coordinate $t$ and the torus fiber with coordinates $z_1,...,z_n$. Let  $E_0=X\times (\C^\ast)^n\to X$ be the trivial fiber bundle.

We will endow the restriction of $E$ to $X\times (\C-\{0\})$ with a non-linear connection $\nabla^{nl,st}$ which extends to $E$ as a meromorphic connection having a pole of order $2$ at $t=0$. For that consider the collection of sections $w_i=e^{Z_i/t}z_i, i=1,...,n$ and declare them to be  flat sections of $\nabla^{nl,st}$. It follows that 
$\nabla^{nl,st}_{\partial_t}=\partial_t+\sum_{1\le i\le n}{Z_i\over{t^2}}z_i\partial_{z_i}$ and $\nabla^{nl,st}_{\partial_{Z_i}}=\partial_{Z_i}-{1\over{t}}z_i\partial_{z_i}, 1\le i\le n$. Here to simplify the notation we denote the vector field $(0,\partial_t)$ by $\partial_t$ and $(\partial_{Z_i}, 0)$ by $\partial_{Z_i}$.

Let now $\xi=(\xi_1,...,\xi_n)\in Hom(\Gamma, \C)$, and $\partial_\xi=\sum_{1\le i\le n}\xi_i\partial_{\xi_i}$  denote the corresponding constant vector field on $X$. Then $\nabla^{nl,st}_{t\partial_\xi}=\sum_{1\le i\le n}\xi_iz_i\partial_{z_i}$ is a {\it holomorphic} vector field on $E$. Furthermore, the lift $\nabla^{nl,st}_{t^2\partial_t}$ of the vector field $(0, t^2\partial_t)$ is equal to $t^2\partial_t+\sum_{1\le i\le n}Z_iz_i\partial_{z_i}$ is holomorphic and its restriction to the fiber of $E_0$ to $X\times \{0\}$ is constant in logarithmic coordinates.

We also notice the $\C^\ast$-action on $X$ corresponding to the rescaling $Z_i\mapsto \lambda Z_i, 1\le i\le n$, which can be naturally extended to $X\times \C$ by the rescaling $t\mapsto \lambda t$ of the coordinate $t$. The non-linear connection $\nabla^{nl,st}$ lifts the corresponding vector field, and the lift preserves the space of flat sections. Geometrically we can say that the lift of the $\C^\ast$-action preserves the foliation defined by $\nabla^{nl,st}$.

\begin{defn}
We call the non-linear connection $\nabla^{nl,st}$ the {\it standard connection}.
\end{defn}

Later we will prove that standard connections are rigid (i.e. undeformable) in a bigger class of non-linear connections which we will discuss in the next subsection.

\subsection{Almost standard connections}\label{almost standard connections}

Let $\Gamma$ be a free abelian group of rank $n$. We fix an isomorphism of groups $\Gamma\simeq \Z^n$. Let $U\subset Hom(\Gamma, \C)$ be a domain of holomorphy.\footnote{This is a technical condition which simplifies our considerations.}

Let $\pi: E\to U\times Spf(\C[[t]])$ be a non-linear analytic fiber bundle over the formal disc endowed with the following data:

a) An isomorphism $\psi_0: \pi^{-1}(U\times \{0\})\simeq U\times Hom(\Gamma, \C^\ast)=U\times (\C^\ast)^n$ which commutes with $\pi$.

b) A non-linear flat connection $\nabla^{nl}$ with the pole of finite order at $t=0$. \footnote{It is assumed that $\nabla^{nl}$ is analytic along $U$ and formal along $Spf(\C[[t]])$. Same concerns trivializations, etc.}

Choose a trivialization $\psi: E\simeq (\C^\ast)^n\times (U\times Spf(\C[[t]]))$ of the non-linear fiber bundle $E$ which analytically extends $\psi_0$.

Consider the following class of non-linear connections $\nabla^{nl}$ on $E$ which in the trivialization $\psi$ satisfies the following:

1) {\bf Polar part property along the base}:

$$\nabla^{nl}_{\partial_t}=\partial_t+t^{-2}\sum_{1\le i\le n}Z_iz_i\partial_{z_i}+t^{-1}v_{-1}+...,$$
where $Z_i, 1\le i\le n$ are some complex numbers numbers which should be thought of as coordinates of a point $Z\in U\subset \C^n$,  the vector field $v_{-1}\in \OO^{an}(U)\widehat{\otimes}Vect_{(\C^\ast)^n}$ is an analytic vector field along the fiber which analytically depends on a point of $U$, dots mean terms in the formal expansion of non-negative degrees in $t$ with coefficients which belong to $ \OO^{an}(U)\widehat{\otimes}Vect_{(\C^\ast)^n}$ .

2) {\bf Polar part property along the fiber}:

$\nabla^{nl}_{\partial_{Z_j}}=\partial_{Z_j}-{1\over{t}}z_j\partial_{z_j}+..., 1\le j\le n$, where dots have the same meaning as in 1).

\begin{defn} Let us call connections $\nabla^{nl}$ satisfying the Polar properties 1) and 2) {\it almost standard}.

\end{defn}

Recall, that the Lie algebra $\g=Vect_{(\C^\ast)^n}$  of vector fields has $\Z^n$-grading such that degree $\gamma$ component $\g_\gamma$ is spanned over $\C$ by $z^{\gamma}\cdot z_i\partial_{z_i}, 1\le i\le n$.
In particular vector fields $z_i\partial_{z_i}, 1\le i\le n$ have degree zero. We will denote by $\g^{an}$ the analytic version of $\g$, i.e. $\OO^{an}(U)\widehat{\otimes}\g$. It inherits the grading from $\g$. The coefficients for $t^k, k{\ge 0}$ denoted by dots belong to $\g^{an}$.

Let us write $v_{-1}=\sum_{1\le i\le n}\delta_i z_i\partial_{z_i}+v$, where $v$ does not contain terms of degree $0$.
Let us write $\nabla^{nl}_{\partial_{Z_j}}=\partial_{Z_j}-{1\over{t}}z_j\partial_{z_j}+P_j$, where $P_j=\sum_{m\ge 1}t^mw_m^{(j)}$ and $w_m^{(j)}\in \g^{an}$.

\begin{prp}\label{independence of trivialization}
a) The conditions 1) and 2) are preserved after a change of  trivialization $\psi$. 

b) The degree zero parts of $v_{-1}$ and $w_0^{(j)}$ do not depend on the trivialization $\psi$.
In particular the functions $\delta_i\in \OO^{an}(U)$ do not depend on $\psi$.

\end{prp}

{\it Proof.}  Formal trivializations extending $\psi_0$ form a torsor over the pronilpotent Lie group $H$ such that $Lie(H)=\{h=\sum_{l\ge 1}t^lh_l|h_l\in \g^{an}\}$.The action of $H$ preserves the term with $t^{-2}$ in 1) and term with $t^{-1}$ in  2). This proves a).

Moreover, it preserves degree zero part of the term with $t^{-1}$ in 1). Indeed, the input of $h$ in the coefficient for $t^{-1}$ is $[h_1,\sum_{1\le i\le n}Z_iz_i\partial_{Z_i}]$. Since vector fields $z_i\partial_{z_i}$ preserve the degree and commute with terms of degree zero, we see that there is not input from the commutator to the coefficient for $t^{-1}$.
Hence $\delta=(\delta_i)_{1\le i\le n}$ is invariant with respect to the change of trivialization. 

In a similar way the coefficient $w_0^{(j)}$ for $t^0$ in 2) gets changed by $[h_1, w_0^{(j)}]$. Same argument as above shows that the degree $0$ part of $w_0^{(j)}$ is preserved for any $1\le j\le n$. 
The Proposition is proven.
$\blacksquare$.

Let us write $w_0^{(j)}=\sum_{1\le i\le n}a_{ij}z_i\partial_{z_i}$, where $a_{ij}\in \OO^{an}(U)$.
The flatness condition $[\nabla^{nl}_{\partial_{Z_{j_1}}}, \nabla^{nl}_{\partial_{Z_{j_2}}}]=0$ implies $\partial_{Z_{j_1}}a_{i,j_2}-\partial_{Z_{j_2}}a_{i,j_1}=0$. Thus we have a collection of closed $1$-forms $\alpha_i:=\alpha_i(Z)=\sum_{1\le j\le n}a_{ij}(Z)dZ_j, 1\le i\le n, Z\in U$.

Proposition \ref{independence of trivialization} implies that the class of connections described by 1) and 2) above has the following invariants:

a) $Z=(Z_i)_{1\le i\le n}$, where each $Z_i$ can be thought of as a restriction of the coordinate function on $Hom(\Gamma, \C)=\C^n$ to $U$.

b) $\delta=(\delta_i)_{1\le i\le n}$, where $\delta_i\in \OO^{an}(U)$.

c) Collection of closed analytic $1$-forms $\alpha_i, 1\le i\le n$ on $U$.

\subsection{Formal classification of almost standard connections}\label{nl HLT}

Aim of this subsection is to prove the following result.

\begin{thm}\label{nl HLT} Let $\nabla^{nl}$ be an almost standard connection on $E$. Assume that on $U$ we have $\delta=(\delta_i)_{1\le i\le n}=0$ and $\alpha_i=0, 1\le i\le n$. Then $\nabla^{nl}$ is gauge equivalent over formal series in $t$ to $\nabla^{nl,st}$. The corresponding gauge equivalence is unique.

\end{thm}

{\it Proof.}  The idea of the proof  is to interpret the existence of the desired gauge transformation and its uniqueness in terms of vanishing of the cohomology groups of a certain complex. 

First we notice that $\g^{an}=\g^{an}_0\oplus \g^{an}_{\ne 0}$, where subscripts in the direct sum refer to the degrees of terms belonging to the corresponding summand. Then $\g^{an}_0$ is a free $\OO^{an}(U)$-module with the basis $z_i\partial_{z_i}, 1\le i\le n$. For $\gamma\in \Gamma=\Z^n$ we use the notation  $v^\gamma$ in order to denote an element of degree $\gamma$.
Let us define the following $\OO^{an}(U)$-linear maps $\g^{an}\to \g^{an}$:
$$S_0: v^\gamma\mapsto Z(\gamma)v^\gamma, S_j: v^\gamma\mapsto \gamma_jv^\gamma, 1\le j\le n,$$
where $\gamma_j$ denote the $j$-th coordinate of $\gamma$ (recall that we fixed an isomorphism $\Gamma\simeq \Z^n$, hence the coordinates of $\gamma$ are well-defined). Notice that $S_0=[\sum_{1\le i\le n}Z_iz_i\partial_{z_i},\bullet]$ and $S_j=[z_j\partial_{z_j}, \bullet], 1\le j\le n$. Here we understand $Z_i, 1\le i\le n$ as coordinate functions on $U$. Notice that $S_aS_b=S_bS_a, 0\le a,b\le n$.

Define a $3$-term complex $0\to C^0\to C^1\to C^2\to 0$ of $\OO^{an}(U)$-modules such as follows:

$$C^0=\g^{an}, C^1=\g^{an}\oplus (\oplus_{1\le j\le n}\g^{an}), C^2=\oplus_{1\le i\le n}\g^{an},$$
where the linear map $d_{01}: C^0\to C^1$ is defined as $x\mapsto (S_0(x),\oplus_{1\le j\le n}(-S_j(x))$, and
$d_{12}: C^1\to C^2$ is defined as $(x_0,x_1,...,x_n)\mapsto \oplus_{1\le i\le n}S_i(x_0)+\oplus_{1\le i\le n}S_0(x_i)$.
One checks that $d_{12}\circ d_{01}=0$. 

Next we observe that $Ker(S_0)=\g^{an}_0$. Indeed, $\g^{an}_0\subset Ker(S_0)$. On the other hand if $v^\gamma\in Ker(S_0), \gamma\ne 0$
we see that $Z(\gamma)=0$ as a linear function on $U$. This is possible only if all coordinates of $\gamma$ are equal to zero, which is a contradiction. This  observation implies that $Ker(d_{01})=\g^{an}_0$. 

Let us now compute $Ker(d_{12})$. By definition it consists of the collections $(x_0,...,x_n)\in (\g^{an})^{n+1}$ such that 
$$\gamma_1x_0=-Z(\gamma)x_1, \gamma_2x_0=-Z(\gamma)x_2,...,\gamma_nx_0=-Z(\gamma)x_n.$$
We know that at least one of the coordinates of $\gamma$ is not equal to zero. We can assume for simplicity that $\gamma_1\ne 0$. Then we have $x_1=-{\gamma_1x_0\over {Z(\gamma)}}$ as a vector field on $(\C^\ast)^n$ which has coefficients analytic in $U$. This means that $x_0=x_0(\gamma)$ is divisible by $Z(\gamma)$.

\begin{lmm}The series $R=\sum_{\gamma\in \Gamma-\{0\}}{x_0(\gamma)\over{Z(\gamma)}}z^\gamma$ is a vector field on $(\C^\ast)^n$ which has analytic coefficients  in $U$.

\end{lmm}

In order to prove the Lemma, we cover $\R^n-\{0\}$ by $2n$ convex domains $V_i^\pm, 1\le i\le n$, such that $\Z^n-\{0\}\subset \cup_iV_i^{pm}$,
the domains $V_i^+$ and $V_i^-$ are symmetric with respect to $0\in \R^n$ and for each $\gamma\in \Z^n-\{0\}$ there exists index $i$ such that the coordinate $\gamma_i\ne 0$ as long as $\gamma\in V_i^\pm$. Then the series $R$ can be split into the finite sum of expressions of the type
$-\sum_{\gamma\in V_i^{\pm}}{x_i(\gamma)\over{\gamma_i}}z^\gamma$. The latter series converges because of the exponential estimates on the analytic in $U$ coefficients of $x_i(\gamma)$. Hence $R$ is convergent.

Let us now apply gauge transformations from the pronilpotent group $H$, or, equivalently from its Lie algebra $Lie(H)$ in order to transform almost standard $\nabla^{nl}$ to $\nabla^{nl,st}$. We have $\nabla^{nl}_{\partial_t}=\nabla^{nl,st}_{\partial_t}+t^{-1}v_{-1}+t^0v_0+...$
and $\nabla^{nl}_{\partial_{Z_j}}=\nabla^{nl,st}_{\partial_{Z_j}}+t^0w_0^{(j)}+t^1w_1^{(j)}+..., 1\le j\le n$. The coefficients $v_i$ and $w_i^{(j)}$ belong to $\g^{an}$. Moreover, by our assumptions we have $\delta_i=\alpha_i=0$, hence $v_{-1}, w_0^{(j)}\in \g^{an}_{\ne 0}$.

We would like to find an element $h=\sum_{i\ge 1}t^ih_i\in Lie(H)$ such that the gauge transformation given by $h$ removes both terms $t^{-1}v_{-1}$ and all terms $t^0w_0^{(j)}, 1\le j\le n$. Applying $[h, \bullet]$ to $\nabla^{nl,st}_{\partial_t}$ and taking the coefficient for $t^{-1}$ we arrive to the equation
$$v_{-1}=S_0(h_1).$$
Similarly, applying  $[h, \bullet]$ to $\nabla^{nl,st}_{\partial_{Z_j}}$ and taking the coefficient for $t^0$ we arrive for each $1\le j\le n$ to the equation
$$w_0^{(j)}=-S_j(h_1).$$

Consider the complex $C^{\bullet}_{\ne 0}$ obtained from $C^\bullet$ by keeping terms with $\g^{an}_{\ne 0}$ only.
Since the cohomology of $C^\bullet_{\ne 0}$ are trivial in the middle term, we conclude that  such $h_1$ does exist. \footnote{In order to prove surjectivity of $d_{01}$ we use the Lemma, which guarantees that the solution has analytic coefficients.} Proceeding by induction in  powers $t^l$, we can ``kill" all for all terms in $\nabla^{nl}_{\partial_t}=\nabla^{nl,st}_{\partial_t}+t^{-1}v_{-1}+t^0v_0+...$
and $\nabla^{nl}_{\partial_{Z_j}}=\nabla^{nl,st}_{\partial_{Z_j}}+t^0w_0^{(j)}+t^1w_1^{(j)}+..., 1\le j\le n$  denoted by dots their parts which do not belong to $\g^{an}_{\ne 0}$.

In order to finish the proof of the Theorem we have to show that by a suitable choice of $h\in Lie(H)$ we can also ``kill" parts of the above-mentioned terms which belong to $\g^{an}_0$. Thus we may assume that in both expressions  $\nabla^{nl}_{\partial_t}=\nabla^{nl,st}_{\partial_t}+t^{-1}v_{-1}+t^0v_0+...$
and $\nabla^{nl}_{\partial_{Z_j}}=\nabla^{nl,st}_{\partial_{Z_j}}+t^0w_0^{(j)}+t^1w_1^{(j)}+..., 1\le j\le n$ all the coefficients belong to $\g^{an}_0$. 
Recall that  by our assumptions  $\delta_i=\alpha_i=0, 1\le i\le n$, hence $v_{-1}=0, w_0^{(j)}=0$.
Let us consider first non-trivial terms in both expression, i.e.  $t^0v_0$ and $t^1w_1^{(j)}, 1\le j\le n$. As before we can find $h_1$ such that $S_0(h_1)=v_0$. Then from the flatness condition $[\nabla^{nl}_{\partial_t}, \nabla^{nl}_{\partial_{Z_j}}]=0, 1\le j\le n$ we see that the coefficients in  for $t^0$ in the commutators are equal to $w_1^{(j)}$. Hence $w_1^{(j)}=0$. We can do the same thing with terms $t^1v_1$ and $t^2w_2^{(j)}$. For that we can find the gauge transformation which starts with $t^2h_2$. Proceeding by induction we construct an element $\sum_{k\ge 1}t^kh_k\in Lie(H)$ which transforms $\nabla^{nl}$ to $\nabla^{nl,st}$. In other words all terms denoted by dots will be ``killed" by an appropriate gauge transformation $h$ such that $h_k\in \g^{an}_0$. Notice that on each step the corresponding $h_k$ is determined uniquely from the equations. This concludes the proof of the Theorem.
$\blacksquare$

\begin{rmk} \label{nl HLT and Liouville} 

a) Given an almost standard connection $\nabla^{nl}$ assume that the collection of complex numbers $Z=(Z_i)_{1\le i\le n}$ satisfies the following Liouville property: there exists $C>0$ such that for any $\gamma\ne 0$ we have $|Z(\gamma)|>e^{-c|\gamma|}$ (in particular $Z$
gives rise to an embedding of abelian groups $\Gamma\to \C$). Assume also that $\delta_i=0, \alpha_i=0, 1\le i\le n$. Then we can find a formal gauge transformation which transforms $\nabla^{nl}$ to $\nabla^{nl, st}$. This result is slightly different from our main theorem, since it is formulated for a fixed point $Z$ rather than for an open neighborhood.

b) The reader has noticed the similarity of Theorem \ref{nl HLT} with the formal classification of meromorphic connections (Hukuhara-Levelt-Turrittin theorem). In general the fact that $\nabla^{nl}$ can be transformed into the standard form leads to infinitely many conditions on the coefficients of the corresponding $\infty$-jets. Our result can be thought of as a ``non-linear HLT theorem with parameters". It suggests that the same purpose can be achieved by imposing finitely many conditions.

\end{rmk} 






\subsection{Axiomatic approach to analytic stability data}\label{data and axioms}

Now, having in mind the above considerations, let us discuss the axiomatic approach to the general  analytic stability data. Let $\Gamma$ be a free abelian group of rank $n$ as before.

{\bf Data and axioms:}

1) {\it Complex manifold $X$ endowed with a trivial local system of free abelian groups ${\Gamma}$, together with a morphism of complex manifolds $\phi: X\to  Hom({\Gamma}, \C)$, which is a local isomorphism.}



2) Let $E_0$ denote the trivial analytic fiber bundle $X\times {\bf T}_\Gamma\to X$, which we identify with the trivial analytic bundle over $X\times\{0\}\subset X\times \C$. The latter pair of spaces gives rise to a germ $Y$ of complex manifold $X\times \C$ at $X\times \{0\}$. The next piece of data is 
{\it a germ of analytic bundle $E\to Y$ which extends $E_0$.} 

More precisely we have a germ of complex  manifold $E$ at $E_0$, such that $dim_\C E=dim_\C E_0+1$ as well as morphism of germs $\pi:E\to Y$ which is surjective at the level of tangent spaces
and such that  the restriction of $E$ to $X\times \{0\}$ coincides with $E_0$. 

3) {\it A non-linear meromorphic flat connection $\nabla^{nl}$ on $E$ (meromorphic flat Ehresmann connection), such that the set of poles of $\nabla^{nl}$ coincides with  $E_0$.} Below in 7) we will explain that the pole in $t$ is of order $2$.

4) {\it Action of $\C^\ast$ on $X$ which under the isomorphism $\phi$ corresponds to the rescaling. We denote it by $x\mapsto \lambda x$.}
Notice  that the $d\phi$-image of corresponding vector field on $X$ is the vector field equal to $\xi$ at the point $\xi\in Hom(\Gamma, \C)$.
This action gives rise to the action of $\C^\ast$ on $X\times \C, (x, t)\mapsto (\lambda x, \lambda t)$. 

5) {\it Lift of the above $\C^\ast$-action  to $E$ such that the map $\pi$ above is $\C^\ast$-equivariant and moreover orbits of the lifted $\C^\ast$-action belong to the leaves of the foliation given by $\nabla^{nl}$.} The latter assumption implies that the above $\C^\ast$-action is unique.

6) Let $\partial_\xi$ denote the constant vector field on $Hom(\Gamma, \C) $ corresponding to the
element $\xi\in Hom(\Gamma, \C)$. Define the vector field $v_\xi$ on $X\times \C$ by $v_\xi=(t\partial_\xi, 0)$. Consider this vector field as the vector field on the germ $Y$. Let us lift $v_\xi$ to $E$ via $\nabla^{nl}$. A priori we obtain a meromorphic vector field.

{\it We require that the lift is holomorphic vector field on $E$.
Furthermore, at $t=0$ the lifted vector field is tangent to ${\bf T}_\Gamma$  at each point $x\in X=X\times \{0\}$, and we require that the corresponding vector field on ${\bf T}_\Gamma$ is the one defined by $\xi$ (understood as a constant vector field in logarithmic coordinates on ${\bf T}_\Gamma$). }


7) Consider the vector field $w=(0, t^2\partial_t)$ on $X\times \C$ and keep the same notation for the corresponding vector field on $Y$. 

{\it Similarly to 6) we require that the $\nabla^{nl}$-lift of $w$ is holomorphic on $E$ and its restriction to the fiber of $E_0$ at $x\in X=X\times \{0\}$ is a constant vector field in logarithmic coordinates corresponding to $\phi(x)$.}

8)  Let $Z_0\in X$. We keep the notation $\nabla^{nl}$ for the induced non-linear flat connection on the bundle on the product $U\times Spf(\C[[t]])$, where $U$ is an open in analytic topology neighborhood of $Z_0$. We require that in the notation of the Subsection \ref{almost standard connections} we have 
$\delta_i(Z)=0, 1\le i\le n$ and $\alpha_i(Z)=0, 1\le i\le n, Z\in U_{Z_0}$.

\begin{conj}\label{analytic stability and non-linear connection}
The data and axioms described in 1)-8) are in bijection with the set of continuous families of analytic stability data parametrized by  a connected component $X\subset Stab(Vect_\Gamma)$, where $\phi$ is the projection to the space of central charges.

\end{conj}

Let us discuss arguments in favor of the Conjecture \ref{analytic stability and non-linear connection}.

\begin{prp}\label{from stability to connection}
There is a natural map from the set of analytic stability data on $Vect_\Gamma$ to the set of data 1)-7).

\end{prp}

{\it Proof.} Let $X$ be a connected component in the topological space $Stab^{an}(Vect_\Gamma)$ of analytic stability data on $Vect_\Gamma$. Then we have the natural map $\phi$ to $Hom(\Gamma, \C)$. Furthermore, for $\sigma\in X$ we have the analytic fiber bundle ${\mathcal F}(\sigma)$  over a germ of $X\times \{0\}\subset X\times \C$ described previously (see e.g. Section \ref{resurgence conjecture}). As we explained earlier the restriction of ${\mathcal F}(\sigma)$ carries a non-linear flat connection induced by the $1$-parameter semigroup $t\mapsto exp(-Z(\beta)/t)$, where $\beta$ is  an element of a strict convex cone.  More precisely, for any finite cyclic collection of admissible rays $l_i, 1\le i\le m$ in $\C_t$ which bound admissible open sectors $V_i$ we have a trivial fiber bundle over each $V_i$ endowed with the standard non-linear flat connection $\nabla^{nl,st}$ described previously. We glue a new fiber bundle using automorphisms $g_{i,t}$. These automorphisms respect $\nabla^{nl,st}$. Since the collection of rays $l_i$ can be chosen arbitrarily, we have the analytic bundle with non-linear flat connection satisfying the conditions 1)-8). $\blacksquare$

Let us discuss the  inverse map, which would make the Conjecture \ref{analytic stability and non-linear connection} true.  For that, starting with the data 1)-8) we need to construct analytic stability data.
On the first step, having $Z\in Hom(\Gamma, \C)$ and $x\in X$ such that $\phi(x)=Z$, we need to construct stability data on $Vect_\Gamma$. On the second step we have to prove that they are analytic. 

We expect that the first step should follow from the formal non-linear version of the Hukuhara-Levelt-Turrittin theorem (HLT theorem for short) which ensures that when we consider the data 1)-8) over the formal series in $t$, then they should be isomorphic to the standard data, associated with the trivial stability data, as we discussed at the beginning of this Subsection. This result should be true as long as we consider central charges which belong to a small neighborhood $U_0$ of $Z_0$ such that $Z_0: \Gamma\to \C$ is an embedding and elements of $Z_0(\Gamma)$ do not admit good approximation by vectors with rational coordinate (Liouville property).

Having the result over formal series in $t$ we can try to upgrade it to the analytic version, using the data 1)-8) for constructing analogs of Stokes automorphisms such as follows. 


Let us now consider a finite collection of admissible rational rays $l_i, 1\le i\le m$ in $\C=\C_t$. They define the open rational cyclic cover of $\C^\ast$ by half-planes $\alpha_i$ orthogonal to $l_i$. Then we have a cyclically ordered collection of open admissible sectors $V_i=\alpha_i\cap \alpha_{i+1}$. 

We denote by $\{w_1^{(i)},...,w_n^{(i)}\}, w_l^{(i)}:=w_l=w_l(Z)$ a basis of flat sections of $\nabla^{nl}$ restricted to $U_0\times \alpha_i$. 
Since over the formal series everything is equivalent to the standard model, we may assume that over $U_0\times \alpha_i$ the non-linear flat connection  $\nabla^{nl}$ is equivalent to $\nabla^{nl,st}$.

Let us introduce the coordinates on the torus fiber over $t\in \C^\ast$ by the formulas $z_i=e^{-Z(e_i)/t}w_i, 1\le i\le n$, where $(e_i)_{1\le i\le n}$ is the standard basis in $\Gamma=\Z^n$.
Then the analytic transition functions which determine how $E$ is glued from its local trivializations on each $V_i$ can be chosen in the form
$$g_{V_i,t}(z_k)=z_k(1+\sum_{\beta\in \Z^n}c_{ik\beta}e^{-Z(\beta)/t}z^\beta), 1\le k\le n, t\in V_i$$
where $c_{ik\beta}\in \C$ (here we use the $\nabla^{nl}$-flatness condition). We need the Support Property which implies that the above sum
runs over $\beta \in C_i\cap \Z^n$, where $C_i$ is a strict convex cone in $\R^n$. The results of Section \ref{families and support property} give an estimate on the support of the arising family over $U_0$ of stability data. The condition $\beta\in C_i$ should also imply that $Re(Z(\beta)/t)>0$ as long as $t\in V_i$. 
Then the Taylor expansions as $t\to 0, t\in V_i$ of the  transition functions  give rise to stability data on $Vect_\Gamma$ which coincide with those constructed over the formal series on the previous step.

Since $E$ is analytic and $\nabla^{nl}$ is meromorphic it follows that each $g_{V_i,t}$ is analytic for fixed $t$. One can prove that similarly defined  $g_{l_i,t}$ is analytic.  Looking at the rational central charges in $U_0$, choosing the rays $l_i$ to be rational and applying Proposition \ref{analyticity rational} we will deduce that all stability data parametrized by $U_0$ are in fact analytic stability data.

\begin{rmk}

a) The data 1)-8) can be generalized to non-trivial local systems, thus giving a description of analytic wall-crossing structures.

b) The data 1)-8) resemble those appearing in the theory of non-commutative Hodge structures (see [KaKoPa]).

c) A generalization of analytic stability data which should be related to the notion of analytic wall-crossing structure arises when instead of the condition 8)  we require that $\delta_i(Z)=c_i, 1\le i\le n$, with some constants $c_i\in \C$.

d) There are other axiomatic approaches to analytic properties of families of stability data. E.g. in [Br] the author considers stability data with exponential bounds on the growth of DT-invariants (the case of Hamiltonian vector fields on the torus of characters). We cannot demonstrate that the notion of analytic of stability data introduced in our paper is different from the one with exponential bounds. The advantage of our approach is  the possibility to prove that analytic stability data form an open subset, hence analyticity is preserved in families. It is not clear how to derive this result from exponential bounds. 

\end{rmk}

\section{Final remarks: more examples and other classes of stability data}\label{final remarks}

\subsection{Holomorphic Floer Theory and Voros resurgence for WKB expansions}

Let $(M,\omega^{2,0})$ be a complex algebraic symplectic manifold, $L\subset M$ a smooth closed  algebraic Lagrangian subvariety with a chosen spin structure (i.e. a choice of $K_L^{1/2}$). We assume that $M$ and $L$ ``behave well at infinity".\footnote{This means that pseudoholomorphic discs for some almost complex structure compatible with $Re(\omega^{2,0})$ which have boundary on $L$ belong to a compact subset of $M$.}

Let $\Gamma:=H_2(M,L,\Z)$ and $Z:\Gamma\to \C$ be a homomorphism (central charge) given by $Z(\gamma)=\int_\gamma\omega^{2,0}$.
Let ${\bf T}$  denote  the torus of characters $Hom(\Gamma,\C^{\ast})$.
The torus ${\bf T}$ has an interpretation in terms of the Fukaya category ${\mathcal F}(M, Re(\omega^{2,0})+iB)$, where $B=Im(\omega^{2,0})$ is the $B$-field.\footnote{As the first approximation it parametrizes rank one complex local systems $\rho$ on $L$. 
More precisely ${\bf T}$ parametrizes deformations of the Fukaya category together with the pair $(L, \rho)$, modulo the action on this category of the group  $Hom(H_1(M, Z), \C^\ast)=H^1(M,\C^\ast)$.}

The stability data on $Vect_\Gamma$ are   determined by the Stokes automorphisms defined in terms of the virtual number of pseudo-holomorphic discs with the boundary on $L$. More precisely, there are Stokes rays with corresponding non-trivial automorphisms $g_\theta:=g_{Stokes}$ of an appropriate domain in the anlytification ${\bf T}^{an}$ of the algebraic torus ${\bf T}$. Stokes rays are defined by the condition $\theta=Arg(Z(\beta))$ for some non-trivial $\beta\in H_2(M,L,\Z)$. This is the situation when there is a 
pseudo-holomorphic disc (for a generic compatible with $Re(\omega^{2,0})$ almost complex structure on $M$) with the boundary on $L$. For a wheel of admissible sectors $V_i, i\in I$ existence of convex cones $C(V_i)$
satisfying the Support Property follows from the Gromov-type estimate
for the area of a disc with the class $\beta$ in terms of the length of its boundary. Let us {\it assume} uniform estimates
$|c_{i\gamma\beta}|\le e^{-const\,|\beta|}$ for the coefficients of the Stokes automorphisms $g_{Stokes}$. Under this assumption
the Conjecture \ref{analyticity and resurgence} holds. This amounts to  the resurgence properties of Gromov-Witten series for complex algebraic symplectic manifolds.  Since this and other examples in this paper are given just for illustration purposes, we postpone the detailed discussion until further publications on our project on Holomorphic Floer Theory (see [KoSo4]).

Let $M=T^\ast C$ be the cotangent bundle of a complex projective curve $C$, and $L\subset M$ to be the union of a fixed cotangent fiber $T_{x_0}^\ast C$ and a smooth spectral curve $\Sigma\subset M$. We can apply the above considerations to this case. \footnote{The fact that $L$ is not smooth does not create problems as long as $\Sigma$ has only Morse ramification points.} Then the above-described stability  data control the resurgence properties of the WKB expansion in the Planck constant $t=\hbar$ of the value of the wave  function $\psi$  at $x_0\in \C$ of the ordinary linear differential equation with the small parameter $\hbar$ and the spectral curve $\Sigma$ (Voros resurgence). An example is given by the wave function of the Schr\"odinger operator $P=-(\hbar \partial_x)^2+V(x)$ where $V(x)$ is a polynomial. The  Stokes automorphisms can be written in terms of the so-called Voros symbols. They have been known since 90's (see e.g. [DelPh] or [IwNak] for a more recent exposition). Our geometric considerations provide an alternative point of view on the subject. The story which involves more general holonomic $DQ$-modules is far beyond this paper and will be discussed in one of the subsequent papers on our big project ``Holomorphic Floer Theory" (see [KoSo4]).

Notice that if  we allow  $x_0$ vary  then we will  have more  walls and the corresponding wall-crossing structure is related to $2d-4d$ wall-crossing formulas described in Section 3.9. The general idea which relates the above geometry to differential equations (or, better, bundles with flat connections) is to consider $Hom(\Sigma, T_{x}^\ast C), x\in C$ in the Fukaya category of $T^\ast C$ (treated as $C^\infty$ $4$-dimensional symplectic manifold). Since Floer homology  are invariant with respect to Hamiltonian deformations,  the family  of these vector spaces gives rise  to a bundle with connection on $C$  via the Riemann-Hilbert correspondence.
 
\subsection{Resurgence of exponential integrals}

Let $X$ be a K\"ahler manifold of complex dimension $n$. The K\"ahler hermitian form gives rise to the K\"ahler metric and the symplectic form as its real and imaginary parts respectively.
Assume that $X$ carries a holomorphic volume form (i.e. it is Calabi-Yau).\footnote{In some of the following considerations it suffices to fix any holomorphic top degree form on $X$.}

We are interested in the exponential integral $I_C(t)=\int_Ce^{f/t}vol$ where $f$ is a given holomorphic function on $X$ with Morse critical points $x_1,..,x_k$, and $C$ is an integration cycle.

It is easy to see that $I_C(t)$ is an integer linear combination of the integrals over the cycles $th_{i,\theta+\pi}$  called {\it Lefschetz thimbles} $th_{i,\theta+\pi}$, where $\theta=Arg(t)$.
Each $th_{j,\theta+\pi}$ is the union of gradient lines with respect to the K\"ahler metric of the function $Re(e^{-i\theta}f)=Re(f/t)$ outcoming from the critical point $x_j\in X$. Let $S=\{z_i=f(x_i), 1\le i\le k\}$ be the set of critical values of $f$.  Then $f(th_{i,\theta+\pi})$ is the ray in the $t$-plane $\R^2$ with the vertex $z_i$ and the direction $\theta+\pi$.

Assume that elements of the set of critical values $S$ are in generic position in the sense that no straight line contains three points from $S$.  A {\it Stokes ray}  contains two different critical values  which we   order by their proximity to the vertex.

Thus we have a collection of exponential integrals for all  $t\in \C^{\ast}$ which do not belong to Stokes rays:
$$I_i(t)=\int_{th_{i,\theta+\pi}}e^{f/t}vol.$$

It  is easy to see that if we cross a Stokes ray
$s_{ij}:=s_{\theta_{ij}}$ containing critical values $z_i,z_j, i<j$, then the integral $I_i(t)$ changes such as follows:
$$I_i(t)\mapsto I_i(t)+n_{ij}I_j(t),$$
where $n_{ij}\in \Z$ is the number of gradient trajectories of the function $Re(e^{i(Arg(z_i-z_j))}f)$ joining critical points $x_i$ and $x_j$.

Let us modify the exponential integrals such as follows:
$$I_i^{mod}(t):=\left({1\over{2\pi t}}\right)^{n/2}e^{-z_i/t}I_i(t).$$
Then as $t\to 0$ the stationary phase expansion of $I_i(t)$ gives rise to a formal series for  $I_i^{mod}(t)$ which we denote  by
$c_{i,0}+c_{i,1}t+....\in \C[[t]],$
where $c_{i,0}\ne 0$.
The jump of the modified exponential integral across the Stokes ray $s_{ij}$ is given by 
$$I_i^{mod}\mapsto I_i^{mod}+n_{ij}I_j^{mod}(t)e^{-(z_i-z_j)/t}.$$ 

We define the analytic stability data such as follows.
Let $\Gamma=\Z^k, k=|S|$, and  $Z:\Gamma\to \C$ given by the formula
$Z(e_i)=z_i$, where $(e_i)_{1\le i\le k}$ is the standard basis of $\Gamma$.
With each Stokes ray $s_{ij}$ we associate the automorphisms of the torus ${\bf T}_\Gamma=Hom(\Gamma,\C^\ast)$ with coordinates $(x_1,...,x_k)$ given by the formula
$$A_{s_{ij}}: x_i\mapsto x_i(1+n_{ij}x^{\gamma_{ij}}),$$
where $\gamma_{ij}=e_j-e_i$. Notice that $x^{\gamma_{ij}}=x_jx_i^{-1}$.
Hence the above formula for the jump of $I_i^{mod}$ corresponds to the
formula
$$x_i\mapsto x_i(1+n_{ij}e^{-Z(\gamma_{ij})/t}x^{\gamma_{ij}}).$$
The Conjecture \ref{analyticity and resurgence} predicts that the series $I_i^{mod}(t)$ are resurgent. This is a well-known fact,
which gives one more confirmation of the Conjecture.\footnote{There are many ways to see the resurgence. Geometrically one can deduce it from the upper bounds on the volumes of sections of thimbles $f=const$ with respect to the volume form $vol/df$.}

More explicitly, the vector $\overline{I}^{mod}(t)=(I_1^{mod}(t),...,I_k^{mod}(t))$,  satisfies the Riemann-Hilbert problem on $\C$ with known jumps across the Stokes rays and known asymptotic expansion as $t\to 0$. We will discuss the Riemann-Hilbert problem in the next subsection.

\subsubsection{Discussion of the Riemann-Hilbert problem}\label{RH problem}

The Riemann-Hilbert problem mentioned above is of the following sort: for a sequence of $\C^k$-valued functions\footnote{Here $k$ is the rank of the relative Betti cohomology, which is under our assumptions is equal to the cardinality $|S|=k$.} $\Psi_1(t),...,\Psi_k(t)$ on $\C^{\ast}-\cup_{i,j} (Stokes\,rays\, s_{ij})$ each of which has a formal power asymptotic expansion in $\C[[t]]$ as $t\to 0$, and which satisfy the following jumping conditions along the Stokes rays $s_{ij}$:

$$ \Psi_j\mapsto \Psi_j,$$
$$\Psi_i\mapsto \Psi_i+n_{ij}e^{-{{{z_i-z_j}\over{t}}}}\Psi_j.$$

This collection $(\Psi_i)_{1\le i\le k}$ gives rise to a  holomorphic vector bundle as explained in Section \ref{resurgence conjecture}. Let us introduce coordinates $\Psi_{ij}$ of each $\Psi_i$, i.e. $\Psi_i(t)=(\Psi_{ij}(t))_{1\le j\le k}$.

For any given $R>0$ there is an explicit formula for $\Psi_{ij}(t)$ which gives absolutely convergent series in the disc $|t|<R$ (see [GaMoNe1], formula (C.12)):
$$\Psi_{ij}(t)=\delta_{ij}+$$
$$+\sum_{\substack{s\ge 1\\i=i_0\ne i_1\ne i_2\ne...\ne i_s=j}}{n_{i_0i_1}...n_{i_{s-1}i_s}\over{(2\pi i})^s}\cdot\int dt_1...dt_s{e^{z_{i_1}-z_{i_0}\over{t_1}}...e^{z_{i_{s}}-z_{i_{s-1}}\over{t_s}}\over{(t_1-t_2)...(t_{s-1}-t_s)(t_s-t)}},$$
where integration is taken over the parts of rays $t_l\in \R_{>0}\cdot(z_{i_{l}}-z_{i_{l-1}})$ given by the conditions $|t_l|<R$ for all $1\le l\le s$.

Analytic properties of this type of integrals were discussed in detail in several places, see e.g. [GaMoNe1], [BarbSt]. Since it is a minor issue with respect to the main topic of our paper, we give below just one result for illustration purposes.

\begin{prp}\label{convergence of integrals}

Each of the integrals as well as the whole above series is convergent.
\end{prp}

{\it Sketch of the proof.} The proof is pretty straightforward, so we give only a sketch of the proof.
Convergence of each integral follows from the rapid decay of exponential functions on the corresponding rays. Indeed let $min_i|z_i-z_{i+1}|:=A$. Then the integral is bounded by $C_1e^{-sA/R}$ for some positive constant $C_1$.

Convergence of the infinite sum is not automatic. It follows from the assumption that there are finitely many critical points of $f$. Then  $n_{ij}\le C$ for some $C>0$. Therefore the product of $n_{ij}$'s for each summand is less or equal than $C^s$. Taking $R$ to be sufficiently small we can make each summand of the series  smaller than $e^{-C_2s}$ for some $C_2>0$. This guarantees the convergence. $\blacksquare$

\begin{rmk}\label{infinite-dimensional integrals}
The finiteness condition does not hold in the infinite-dimensional case.
\end{rmk}

The collection of functions $\Psi_{ij}(t)$ can be combined into a $GL(k,\C)$-valued function $M(t)$, which is holomorphic in each Stokes sector and has there a well-defined limit $lim_{t\to 0}M(t)=id$.
It follows that the jumps of $M(t)$ over Stokes rays are given by integer upper-triangular matrices with $1$ on the main diagonal and the $ij$-th matrix element equal to $n_{ij}e^{-{(z_i-z_j)}/t}$.

Then the vector $M^{-1}(t)(\overline{I}(t))$ does not jump across Stokes rays. Hence it gives rise to a $\C^k$-valued holomorphic function in the dics $|t|<R$. Notice that the formal (as $t\to 0$) asymptotics $M_{form}$ belongs to $ GL(k,\C[[t]])$, while the formal asymptotics $\overline{I}_{form}$ belongs to $ \C^k[[t]]$. Both asymptotics can be computed by means of perturbation theory. The above discussion implies that $J_{form}=M_{form}^{-1}(\overline{I}_{form})$ has non-zero radius of convergence.
Equivalently, we can formulate the result of the above discussion in the following way.

\begin{prp} $\overline{I}(t)=M(t)(J_{form}(t))$, where $M(t)$ depends on the set $S$ and integer numbers $n_{ij}$ only, and $J_{form}(t)$ is a $\C^k$-valued formal power series in $t$, which is convergent in some disc. Its coefficients can be computed by formal procedure, which does not depend on the function $f$.

\end{prp}

\begin{rmk} We hope that a similar result about the resurgence holds in Quantum Field Theory in which the coupling constants are not renormalized, with exponential integrals being replaced by Feynman integrals.\footnote{Examples are given by Quantum Mechanics, Chern-Simons theory in $3d$, or $N=4$ Super Yang-Mills.} Then $\overline{I}_{form}$ can be computed by Feynman rules in a pure formal way. We claim that after applying some infinite rank operator to it, we obtain a convergent in the coupling constants series.

\end{rmk}

\subsection{Other classes of stability data}

There are several interesting classes of stability data on $Vect_\Gamma$. Let us mention some of them:

A) Polynomial stability data;

B) Rational stability data;

C) Algebraic stability data;

And we add to the list the class discussed in this paper:

D) Analytic stability data.

Definitions and tools for study these classes rely on the geometric interpretation of stability data presented in  [KoSo3] and Section 2 above. We will not discuss the details in this paper, which is exclusively devoted to the class D.
We just list some facts leaving the details for future publications:

i) There are natural embeddings $A\subset B\subset C\subset D$.

ii) In order  to define classes A-C one uses non-archimedean geometry. Namely, with generic stability data on $Vect_\Gamma$ and an additional choice of the wheel of cones,  one can associate a geometric object which is a {\it non-archimedean} analytic space over a valuation field $K$.  \footnote{We utilize V. Berkovich approach to non-archimedean geometry, see[Be].} In the case when  $K=\C((T))$ this $K$-analytic space is the generic fiber of the formal scheme over $\C$. The formal scheme is a formal neighborhood of a closed subset which is a wheel (or chain) of projective lines in a toric manifold as discussed in this paper. The corresponding analytic space contains a subset $Z^\ast(\R^2-\{0\})$.  We need only a germ of the non-archimedean analytic space at this  locally closed subset. This germ does not depend on the wheel of cones.

iii) The classes A-C are defined in terms of the field of meromorphic functions on this  $K$-analytic germ. They enjoy interesting properties which do not hold for more general analytic stability data. For example, the main  reason for introducing the class C is the theorem which says that   if the central charge is rational then the series $A_V(x^\gamma)$ are {\it algebraic}.
The subclass $B\subset C$ is singled out by the property that the field of meromorphic functions on the natural $\Z$-covering of the $K$-analytic germ has the transcendence degree equal to $n=rk\,\Gamma$.

iv) Most of the examples of analytic stability data coming from ``real life" are rational or even polynomial. Here are few of them:

1) Wild character varieties for local systems on complex curves give rise to open subsets in the non-archimedean spaces arising from polynomial stability data.

2) Class $\mathcal{P}$ of quivers (see [KoSo1], Section 8) and some of its generalizations (see [Lad]) is related  to rational stability data.

3) Conjecturally any canonical  transformation (a.k.a DT-transformation with canonical initial data, see [KoSo3], Section 3.3, [GrHaKeKo]) is rational.


4) Stability data with gaps on $Vect_\Gamma$  such that each $A_{V_i}$  (see loc.cit.) is given by the rational map on monomials $x^\mu, \mu\in \Gamma$ is rational.

v) Open  and closed properties:

1) The space of of polynomial stability data is open and closed  in $Stab(Vect_\Gamma)$.

2) The spaces of rational and algebraic stability data are open and closed in $Stab(Vect_\Gamma)$ provided
a certain Convexity Conjecture is satisfied. The latter roughly says that a meromorphic function on the
non-archimedean tube domain over an open set extends to a one over the convex hull of the set.

\section{Appendix: proof of the Theorem \ref{independence of cyclic cover}}

Below we are going to develop a purely combinatorial formalism of cyclic covers from Section 2. This will allow us to prove independence of certain quotient sets  of the combinatorial version of the cyclic covers. Then we explain how to reduce the geometric set up of the Theorem \ref{independence of cyclic cover} to the combinatorial one.

\subsection{Combinatorial reformulation of open cyclic covers}

\begin{defn} A {\it cyclic ordered set} (a.k.a a {\it cyclic set}) is as  non-empty finite set endowed with the transitive action of group $\mathbf Z$ (denoted by $s\mapsto s+1$).

\end{defn}

For any cyclic set $S$ we define the   the {\it dual} set $S^\vee$ as the set of ``holes" between consecutive elements of $S$ (i.e. ordered pairs $(i,i+1)$ for $i\in S$). The disjoint union $S\sqcup S^\vee$ is again naturally cyclically ordered:
$$ \dots\quad  i \quad(i,i+1)\quad i+1\quad (i+1,i+2) \quad  i+2\quad \dots$$
In what follows we will draw elements of $S$ as $\bullet$, and elements of $S^\vee$ as $-$.

\begin{defn} For a cyclic set $S$ denote by $Int(S)$ (intervals in $S$) the set of nonempty subsets of $S\sqcup S^\vee$ of type
 $$- \qquad  -\bullet - \qquad  -\bullet-\bullet- \qquad \dots$$
 
 \end{defn}

If $S$ has $k$ elements, then $Int(S)$ has $k^2$ elements. Each interval $I$ contains exactly $a$ elements of $S$ and $(a+1)$ elements of $S^\vee$ for some integer $a, \,0\le a\le k-1$, called the {\it length} of $I$.
The set $Int(S)$ is partially ordered (by inclusion), and is endowed with a partially defined associative multiplication
$$A B=C\text{ for intervals }A,B,C\in Int(S)$$
$$\text{ iff the last hole in }A \,is\, equal\, to\,\text{ the first hole in }B\text{, and }C=A\cup B$$
An ordered pair $(A,B)$ of intervals such that $AB$ is well-defined, is called an {\it adjacent} pair:
$$\begin{aligned} A=\quad &-\bullet  - \bullet -\phantom{ \bullet - \bullet -\bullet -\bullet -} \\
  B= \quad &   \phantom{\,\,\,- \bullet- \bullet}\! -\bullet - \bullet - \bullet -\bullet -\\
  AB=\quad &-\bullet  - \bullet - \bullet - \bullet -\bullet -\bullet - \end{aligned}$$

\

From now we will fix a ``base" cyclic set $H$. 

\

\begin{defn} An ordered pair of   intervals $(I_1,I_2)\in Int(H)\times Int(H)$ is  called {\it linked} if there exist (automatically unique)  intervals $A,B,C$ such that the interval $ABC\in Int(H)$ is defined, and
\begin{equation}I_1=AB, \quad I_2=BC\label{deflinked} \end{equation}

\end{defn}

Here is the picture:

$$\begin{aligned}I_1=\quad &-\bullet  - \bullet -\bullet - \bullet -\phantom{ \bullet -\bullet -} \\
  I_2= \quad &  \phantom{\,\,\, - \bullet}- \bullet -\bullet - \bullet - \bullet -\bullet -\end{aligned}$$

with $length(A)=1,length(B)=3,length(C)=2$.
Obviously, an ordered pair $(I_1,I_2)$ is adjacent iff it is linked and $B=I_1\cap I_2$ has length 0.

\

\begin{defn} A {\it combinatorial cyclic cover}  of $H$ is  a collection of  intervals $(I_s)_{s\in S}$ labeled by (another) cyclic set $S$ (in other words, a map $S\to Int(H)$), such that 
$$\forall s\in S:\quad (I_s,I_{s+1})\text{ is linked} $$
and such that the map
$$\left(\sqcup_{s\in S} (I_s-I_{s-1})\right)\cap H\to H$$
is a bijection (informally, this means that the cover is ``of degree 1").
A {\it disjoint} cover of $H$ is a combinatorial cyclic cover $(I_s)_{s\in S}$ such that 
$$\forall s\in S:\quad (I_s,I_{s+1})\text{ are adjacent} $$
\end{defn}
The simplest example of a disjoint cover is the canonical cover  $(\{-i-\})_{i\in H}$ by intervals of length 1.

\begin{rmk}
We use here the adjective ``combinatorial" in order to distinguish between cyclic covers of geometric origin from Section 2 and cyclic covers introduced in this Appendix. In fact both notion agree.
\end{rmk}

\begin{defn} A factorizable system of groups is a subset $SmallInt\subset Int(H)$ containing all intervals of lengths $0$ and $1$, and all subintervals of $I$ for any $I\in SmallInt$, and a functor $$(SmallInt,\subset)\to Groups$$ (written as $I\mapsto G_I$, so we have a group morphism 
$$(I_1\subset I_2)_*:G_{I_1}\to G_{I_2}$$ 
for every ordered pair $I_1\subset I_2\in SmallInt$), such that for any 
$$A,B,C\in SmallInt, AB=C$$ we have
 $$G_A\times G_B\to G_C, \quad (g_A,g_B)\mapsto (A\subset C)_*(g_A)\cdot_{G_C} (B\subset C)_*(g_B)\text{ is a bijection of sets}.$$
 
\end{defn}
Elements of $SmallInt$ are called {\it small intervals}, hence the notation.

 This definition implies that $G_I=\{id\}$ for any interval of length 0, and for any $I=-i_1-i_2-\dots -i_k-\in SmallInt$  we have a canonical bijection (and not a group homomorphism)
 $$ iso_I:G_{i_1}\times \dots \times G_{i_k}\xrightarrow{\sim\text{ as a set }}G_I,\quad G_i:=G_{\{i\}}\,\,\forall i\in H\,.$$

\begin{defn}  For any combinatorial cyclic cover $\varkappa=(I_s)_{s\in S}$ we define the set $X_{\varkappa}$ as
$$ \dots\times_{G_{I_{s-1}\cap I_{s}} } G_{I_s}\times_{G_{I_{s}\cap I_{s+1}}} G_{I_{s+1}}\times_{G_{I_{s+1}\cap I_{s+2}}}\dots$$

\end{defn}

Let ${\bf G}=\prod_{i\in H}G_i$. 
 Our goal is to construct for any combinatorial cyclic covers  $\varkappa$ a {\bf canonical}  bijection   
 $$iso_{\varkappa}: {\bf G}\simeq X_{\varkappa}$$ 

\

First of all, we define $iso_{\varkappa}$ for any {\it disjoint} combinatorial cyclic cover $\varkappa$. Indeed, if $\varkappa=(I_s)_{s\in S}$ is disjoint  then all the groups $G_{I_s \cap I_{s+1}}$ are trivial, therefore $X_\varkappa$ is canonically identified with the product $\prod_{s\in S} G_{I_s}$. Each individual factor $G_{I_s}$, where
$$I_s=-i_{1,s}-i_{2,s}-\dots -i_{k_s,s}-\,,$$
is canonically identified by $iso_I^{-1}$  with the product $\prod_{j=1}^{k_s} G_{i_{j,s}}$. Combining all these isomorphisms we obtain the bijection  $iso_{\varkappa}$.

\

In order to define $iso_{\varkappa}$ for an arbitrary combinatorial cyclic cover, we proceed as follows. First, pick a disjoint cover $\varkappa'=(J_s)_{s\in S}$ labeled by the same cyclic set $S$ of indices, and such  that $J_s\subset I_s$ for any $s\in S$. We call $\varkappa'$ a {\it disjoint subcover} of $\varkappa$. Then, by  functoriality of the assignement $I\rightsquigarrow G_I$,  we have canonical morphisms of groups
$$ G_{J_s}\to G_{I_s},\quad  G_{J_s \cap J_{s+1}}\to G_{I_s\cap I_{s+1}} \qquad \forall s\in S\,,$$
hence a map
 $$f_{\varkappa,\varkappa'}:X_{\varkappa'}\to X_{\varkappa}\,.$$
Combining this map with the bijection $iso_{\varkappa'}$ defined as above, we obtain a map 
$$\varphi_{\varkappa,\varkappa'}: {\bf G}\to X_\varkappa\,.$$

The construction of the canonical isomorphism $iso_{\varkappa}$ is based on the following result.

\

\begin{thm}

1) For any combinatorial cyclic cover $\varkappa$ there exists a disjoint subcover $\varkappa'$ as above.

2) The map $\varphi_{\varkappa,\varkappa'}: {\bf G}\to X_\varkappa$ does not depend on the choice of a disjoint subcover $\varkappa'$ (hence we can define $\varphi_\varkappa$ as $\varphi_{\varkappa,\varkappa'}$ for arbitrary disjoint subcover $\varkappa'$).

3) The map  $\varphi_\varkappa$ is a bijection. 

\end{thm}

Assuming the  theorem, we  define the isomorphism $iso_{\varkappa}:=\varphi_\varkappa$.

{\it Proof.} Statement  1) follows from the fact that
$$\varkappa^{min}:=(J^{min}_s)_{s\in S},\quad J_s^{min}:=(I_s-I_{s+1})\sqcup\{\text{the first hole in }I_{s+1}\}$$
is a disjoint subcover of $\varkappa$.

\ 

In order to prove  independence of $\varphi_{\varkappa,\varkappa'}$ of the choice of a disjoint subcover $\varkappa'$ of $\varkappa$, we observe the following. If a disjoint subcover  $\varkappa'=(J_s)_{s\in S}$ is not equal to $\varkappa^{min}$, then there exists $i_0\in H$ and $s_0\in S$ such that
$$ i_0\in J_{s_0},\quad  (i_0+1)\in J_{s_0+1},\quad i_0\in I_{s_0+1}$$
then we can form another disjoint subcover $\widetilde{\varkappa}'=(\widetilde{J}_s)_{s\in S}$ of $\varkappa$, by swapping $i_0\in H$ from term $s_0$ to term $(s_0+1)$:
$$  \widetilde{J}_{s_0}:=J_{s_0}-\{i_0,-\} , \quad  \widetilde{J}_{s_0+1}:=\{-,i_0\} \sqcup J_{s_0+1}, \qquad \widetilde{J}_s:=J_s \text{ for } s\ne s_0,s_0+1\,.$$
Here is the picture:
$$ \begin{aligned} J_{s_0}=\quad &-\bullet  - \bullet\, - \stackrel{i_0}{\bullet} - \phantom{\bullet -\bullet -\bullet -} \\
  J_{s_0+1}= \quad &   \phantom{- \bullet- \bullet -\bullet }\,\,\,\,\,- \stackrel{\!\!i_0+1\!\!\!}{\bullet}\! - \bullet -\bullet -\\
\, gives\,rise\,to\,\,\,\,\,\,\,\,\,\,\,\,\,\,\,\,\,\,\,\,\, \\
  \widetilde{J}_{s_0}=\quad &-\bullet  -\bullet -\phantom{  \bullet - \bullet -\bullet -\bullet - }\\
 \widetilde{J}_{s_0+1}=\quad &\phantom{-\bullet  - \bullet\,\, }\,\,-  \stackrel{i_0}{\bullet} - \!\stackrel{\!\!i_0+1\!\!\!}{\bullet}\! -\bullet -\bullet - \end{aligned}$$

It follows immediately from the factorization property that 
$$ \varphi_{\varkappa, \varkappa'}=\varphi_{\varkappa, \widetilde{\varkappa}'}\,.$$ 
By making modifications $\varkappa'\rightsquigarrow \widetilde{\varkappa}'$ we reach the final state in which no more modification as above is possible, i.e. when  we reach $\varkappa^{min}$. 

\

Now, after the proving first two statements,  we know that $\varphi_\varkappa$ is well-defined for any combinatorial cyclic cover $\varkappa$. 
In order to prove the third statement (which is the most cumbersome to prove), we will construct  a chain of modifications of combinatorial cyclic covers such that at each step we {\it do not change} the set  $X_\varkappa$ (i.e. the new quotient  set is {\it canonically} identified with the previous one), and also we do not change the map $\varphi_\varkappa$ (after making this identification). The final state is  a disjoint cover, which implies the bijectivity of $\varphi_\varkappa$. 

\

{\bf Step 1}: we eliminate repeating intervals. If $I_s=I_{s+1}$ for some $s\in S$, we remove one of repeating copies without changing the corresponding quotient set. Denote $I:=I_s=I_{s+1}$. We make the modification
$$\varkappa= \dots I_{s-1}, I,I, I_{s+2}\dots \rightsquigarrow \widetilde{\varkappa}= \dots I_{s-1}, I, I_{s+2}\dots, \qquad X_\varkappa\overset{?}{=}X_{\widetilde{\varkappa}}\,.$$
Indeed, we have 
$$\dots \underset{G_{I_{s-1}\cap I}}{\times} G_I \underset{G_I}{\times}  G_I \underset{G_{I \cap I_{s+2}}}{\times}\dots \stackrel{canonically}{\simeq}  \dots \underset{G_{I_{s-1}\cap I}}{\times} G_I \underset{G_{I \cap I_{s+2}}}{\times}\dots $$

\

{\bf Step 2}: we eliminate intervals of length 0. If $I_s=\{-\}$ has length 0, then the pair $(I_{s-1},I_{s+1})$ is linked, and we make the modification
$$\varkappa= \dots I_{s-1}, I_s, I_{s+1}\dots \rightsquigarrow \widetilde{\varkappa}= \dots I_{s-1}, I_{s+1}\dots\,.$$

\

{\bf Step 3}: we eliminate consecutive pairs $(I_s,I_{s+1})$ such that $I_s\not\subset I_{s+1}$ and $I_{s}\not\supset I_{s+1}$. Indeed, if it is the case then $I_s\cap I_{s+1}\ne I_s,I_{s+1}$. The modification is
$$\varkappa= \dots I_{s},  I_{s+1}\dots \rightsquigarrow \widetilde{\varkappa}= \dots I_{s}, I_s\cap I_{s+1}, I_{s+1}\dots, \qquad X_\varkappa\overset{?}{=}X_{\widetilde{\varkappa}}\,,$$
or, writing $I_s=AB, I_{s+1}=BC$ as in \eqref{deflinked},   
$$ \dots AB,BC\dots \rightsquigarrow \dots AB,B,BC\dots$$

This follows from the fact that $G_{AB​}\underset{G_B}{\times}​G_{BC​}$ as the set with $G_{AB}\times G_{BC}^{op}$-action,   coincides with $G_{AB}\underset{G_B}{\times} G_B\underset{G_B}{\times} G_{BC}$.

From now one we assume that for each $s\in S$ we have either 
$$ (I_s,I_{s+1})=(A,AB)\quad \text{ or } \quad  (I_s,I_{s+1})=(AB,B) \qquad\text{ for some }A,B\,.$$

\

{\bf Step 4}: we  eliminate consecutive inclusions in the same direction
$$ \dots A,AB,ABC\dots \rightsquigarrow \dots A, ABC\dots, \text{ or}$$
$$ \dots ABC,BC,C\dots \rightsquigarrow \dots ABC,C\dots \phantom{\text{ or}}$$

Here we use
 $$G_A\underset{G_A}{\times} G_{AB}\underset{G_{AB}}{\times} G_{ABC}=G_A\underset{G_A}{\times}  G_{ABC}, $$
  $$G_{ABC}\underset{G_{BC}}{\times} G_{BC}\underset{G_{C}}{\times} G_{C}=G_{ABC}\underset{G_C}{\times}  G_{C}\,.$$

Let us eliminate now the case 
$$   I_s\subsetneq I_{s+1}\supsetneq  I_{s+2}\,.$$
There are 3 distinct subcases (all intervals $A, C$ below have strictly positive length):
$${\bf I}:  \dots AB, ABC,BC\dots , $$
$${\bf II}:  \dots A, ABC,C\dots, $$
$${\bf III}:  \dots A, AC,C\dots\,.$$

\

 {\bf Step 5}: (for the subcase ${\bf I}$), replace 
 $$   \dots AB, ABC,BC\dots \rightsquigarrow \dots AB, B, BC\dots$$
Here the calculation is a bit non-trivial:
$$\begin{aligned}G_{AB}\times_{G_{AB}} G_{ABC}\times_{G_{BC}}G_{BC}=\\=G_{AB} \times_{G_{AB}} (G_{AB}\times_{G_{B}} G_{BC})\times_{G_{BC}} G_{BC}=\\
=(G_{AB} \times_{G_{AB}} G_{AB})\times_{G_{B}} (G_{BC}\times_{G_{BC}} G_{BC})=  \\
=G_{AB}\times_{G_{B}}  G_{BC}=G_{AB}\times_{G_B} G_B\times_{G_B} G_{BC}\end{aligned}$$
 After performing Step 5, we can arrive to the situation when we have to apply Step 4 again. As the total length of the cover only decreases, the process terminates.

\

 {\bf Step 6}: (for the subcase ${\bf II}$), replace 
 $$   \dots A, ABC,C\dots \rightsquigarrow \dots A, B, C\dots$$
 
\

  {\bf Step 7}: (for the subcase ${\bf III}$), replace 
 $$   \dots A, AC,C\dots \rightsquigarrow \dots A, C\dots$$

Now, when we cannot make any longer  one of modifications from the previous seven steps, we should have at least one consecutive pair $(I_s,I_{s+1})$ which is *adjacent*. All such pairs cut the cyclic cover into a cyclic system of covers of disjoint (in $H$ part) intervals. Moreover, each individual interval is of one of the following 4 forms:
$$ A \quad A,AB \quad AB,B \quad AB,B,BC$$

\ 

{\bf Step 8}: we reach a disjoint cover by using respectively 
$$ A \rightsquigarrow A\quad A,AB \rightsquigarrow A,B \quad AB,B \rightsquigarrow A,B \quad AB,B,BC\rightsquigarrow A,B,C$$ 
This finishes the proof. $\blacksquare$

\subsection{Independence of cyclic cover}

Let $f: S\to S^\prime$ be an embedding of cyclic sets, and $\varkappa=(I_s)_{s\in S}, \varkappa^\prime=(I^\prime_{s^\prime})_{s^\prime\in S^\prime}$ be the collection of intervals such that $I_s\subset I^\prime_{f(s)}, \forall s\in S$. Consider the induced map $\varphi: \prod_{s\in S}G_{I_s}\to \prod_{s^\prime\in S^\prime}G_{I^\prime_{s^\prime}}$.

\begin{thm} \label{combinatorial isomorphism} This map descends to the map $X_{\varkappa}\to X_{\varkappa^\prime}$ and induces a bijection of these two sets.

\end{thm}
{\it Proof.} Proof of the fact that the map gives rise to the map of quotient sets is similar to the one of the Proposition \ref{functor}.
More difficult is to prove the second statement.
For that we notice that the disjoint cover for $\varkappa$ induces via the embedding  $S\to S^\prime$ a disjoint cover for $\varkappa^\prime$. They induce identifications of $X_{\varkappa}$ and $ X_{\varkappa^\prime}$ with ${\bf G}=\prod_{i\in H}G_i$. Moreover these identifications are compatible with the map $\varphi$.
This finishes the proof.
 $\blacksquare$

Let us explain how this result implies Theorem \ref{independence of cyclic cover}. 
First, let us show that it holds for  cyclic covers of any type which do not contain open half-planes in case if we use the groups $G_{V_i}^{old}$ instead of $G_{V_i}$.

Let  $(V_i)_{i\in J}$  be such a cover. The boundary rays of all sectors $V_i$ decompose $\R^2-\{0\}$ into a disjoint union of rays and open sectors. This set of rays and sectors has a natural clockwise cyclic order. By definition this is our cyclic ordered set $H$. Then our cyclic cover gives rise to a combinatorial disjoint cover by some intervals $I_j, j\in H$. The set of short intervals is defined as the set of all  subintervals of all intervals $I_j, j\in H$. Finally as groups $G_{I_s}$ we take groups $G^{old}$ for the corresponding sectors. Applying Theorem \ref{combinatorial isomorphism} we get the result in this case.

Next step is to reduce the Theorem \ref{independence of cyclic cover} to the case we have just considered. For that we start with an open cyclic cover $(V_i)_{i\in J}$ which can include open half-planes, and we will use groups $G_{V_i}$ rather than $G_{V_i}^{old}$.

Let us choose $\epsilon>0$ such that $2\epsilon$ is smaller than the angle between any two  boundary rays of the sectors $V_i, i\in J$. Consider closed admissible sectors $W_{i,\epsilon}\subset V_i$   such that the angle between left (resp. right) boundary ray of $W_{i,\epsilon}$ and left (resp. right) boundary ray of $V_i$ is equal to $\epsilon$. Then we have the groups $G_{W_{i,\epsilon}}^{old}=G_{W_{i,\epsilon}}$. For each $i\in J$ the inductive limit of these groups as $\epsilon\to 0$ is equal to $G_{V_i}$. 
The sectors $W_{i,\epsilon}, i\in J$ form a closed cyclic cover, so the above considerations can be applied. Then we have the quotient set $\MM_{(W_{i,\epsilon})_{i\in J}}^{old}=\MM_{(W_{i,\epsilon})_{i\in J}}$. Moreover, the inductive limit of $\MM_{(W_{i,\epsilon})_{i\in J}}$
as $\epsilon\to 0$ coincides with $\MM_{(V_i)_{i\in J}}$.  This concludes the proof of the Theorem \ref{independence of cyclic cover}.

\begin{cor} There is a natural bijection of the set $\MM_{(W_{i,\epsilon})_{i\in J}}$ with $Stab(\g)$.
\end{cor}
Furthermore, as we explained in Section 2 this bijection  induces a homeomorphism of the corresponding topological spaces.

\subsection{Declarations}

1. There is no conflict of interests.

2. Data sharing not applicable to this article as no datasets were generated or analysed during the current study.

\vspace{3mm}

{\bf References}

\vspace{2mm}





[Ar] D. Arinkin, Moduli of connections with a small parameter on a curve, arXiv:math/0409373.

\vspace{2mm}

[BarbSt] A. Barbieri, J. Stoppa, Frobenius type and CV-structures for Donaldson-Thomas theory and a convergence property,  arXiv:1512.01176.

\vspace{2mm}

 [Be]  V. Berkovich, Spectral theory and analytic geometry over non-archimedean fields,
Mathematical Surveys and Monographs, vol. 33, American Mathematical Society, Providence, RI, 1990, 169 pp.

\vspace{2mm}






[Br]  T. Bridgeland, Riemann-Hilbert problems from Donaldson-Thomas theory,
arXiv:1611.03697.

\vspace{2mm}

[CorSh] C. Cordova, S-H. Shao, Shur indices, BPS particles and Argyres-Douglas theories, arXiv:1506.00265.

\vspace{2mm}

[DelPh] E. Delabaere, F. Pham, Resurgent methods in semi-classical analysis, Annales IHP, Section A, 71:1, 1999, 1-94.
\vspace{2mm}






[Dou] A. Douady, Le probl\`eme des modules pour les sous-espaces analytiques compacts d'un espace analytique donn\'e, Ann. Inst. Fourier, 16:1, 1966, 1-95.




\vspace{2mm}
[Ec] J. \'Ecalle, Les Fonctions r\'esurgentes, Pub. Math. Orsay, 1985.

\vspace{2mm}

[FoGo1] V. Fock, A. Goncharov, Cluster ensembles, quantization and the dilogarithm, Invent.
Math. 175 (2009), no. 2, 223-286, see also  arXiv:math/0311245.


\vspace{2mm}

[GaMoNe1] D. Gaiotto, G. Moore, A. Neitzke,  Four-dimensional wall-crossing via three-dimensional field theory, arXiv:0807.4723.

\vspace{2mm}
[GaMoNe2] D. Gaiotto, G. Moore, A. Neitzke,  Wall-crossing, Hitchin Systems, and the WKB Approximation, arXiv:0907.3987.

\vspace{2mm}

[GaMoNe3] D. Gaiotto, G. Moore, A. Neitzke,  Wall-crossing in coupled 2d-4d systems, arXiv:1103.2598.
\vspace{2mm}


[GarGuM] S. Garoufalidis, J. Gu, M. Marino, The resurgent structure of quantum knot invariants, arXiv:2007.10190.
\vspace{2mm}

[GrHaKeKo]  M. Gross, P. Hacking, S. Keel, M. Kontsevich, Canonical bases for cluster algebras, arXiv:1411.1394.

\vspace{2mm}


[GuMaPu] S. Gukov, M. Marino, P. Putrov, Resurgence in complex Chern-Simons theory, arXiv:1605.07615.
\vspace{2mm}

[IwNak] K. Iwaki, T. Nakanishi, Exact WKB  analysis and cluster algebras, arXiv: 1401.7094.

\vspace{2mm}

[KaKoPa] L. Katzarkov, M. Kontsevich, T. Pantev, Hodge theoretic aspects of mirror symmetry, arXiv:0806.0107. 
\vspace{2mm}

[KoSo1] M. Kontsevich, Y. Soibelman, Stability structures, motivic Donaldson-Thomas invariants and cluster transformations, arXiv:0811.2435.

\vspace{2mm}

[KoSo2] M. Kontsevich, Y. Soibelman, Affine structures and non-archimedean analytic spaces, math.AG/0406564.


\vspace{2mm}








\vspace{2mm}

[KoSo3] M. Kontsevich, Y. Soibelman, Wall-crossing structures in Donaldson-Thomas invariants, integrable systems and Mirror Symmetry, arXiv:1303.3253.

\vspace{2mm}

[KoSo4] M. Kontsevich, Y. Soibelman, Holomorphic Floer Theory, in preparation.

\vspace{2mm}

[Lad] S. Ladkani, Mutation classes of certain quivers with potentials as derived equivalence,  arXiv:1102.4108.


\vspace{2mm}
[MitSau] C. Mitschi, D. Sauzin, Divergent series, summability and resurgence I, Lecture Notes in Math. vol.  2153, Springer, 2016.

\vspace{2mm}




[NiSi] T. Nishinou, B. Siebert, Toric degenerations of toric varieties and tropical curves, Duke Math. J., 135:1, 1-51.

\vspace{2mm}





[Vo] A. Voros, The return of the quartic oscillator (the complex WKB method), Annales Institut H. Poincar\'e, 29:3, 1983, 211-338.

\vspace{2mm}

Addresses:

M.K.: IHES, 35 route de Chartres, F-91440, France, {maxim@ihes.fr}

Y.S.: Department of Mathematics, KSU, Manhattan, KS 66506, USA, {soibel@math.ksu.edu}

\end{document}